%% file: main.tex
\documentclass[final,leqno,onefignum,onetabnum]{siamltex}

\usepackage{fullpage}
\usepackage{comment}
\usepackage{upgreek}
\usepackage{cite}
\usepackage[mathscr]{eucal}
\usepackage{color}
\usepackage[usenames,dvipsnames,svgnames,table]{xcolor}
\usepackage{amsmath,amssymb,amsopn,mathtools}
\usepackage{graphicx}
\usepackage{hyperref}       
\hypersetup{
	colorlinks=true,   
	citecolor=blue,     
	filecolor=blue,     
	linkcolor=red,    
	urlcolor=blue      
}

\usepackage{relsize}
\numberwithin{equation}{section}
\usepackage{enumitem}

\newlist{Assumption}{enumerate}{1}
\setlist[Assumption]{label=A\arabic*}

\definecolor{Blue}{rgb}{0,0,1}
\definecolor{Red}{rgb}{1,0,0}
\definecolor{Green}{rgb}{0,1,0}
\definecolor{Cyan}{rgb}{0,0.72,0.92}
\definecolor{Amethyst}{rgb}{0.6,0.4,0.8}
\definecolor{Bronze}{rgb}{0.8,0.5,0.2}
\definecolor{Violet}{rgb}{0.54,0.17,0.89}

\newcommand{\rev}[1]{#1}
\newcommand{\revv}[1]{#1}

\usepackage{amsmath,amssymb}
\usepackage{epsfig, algorithmic, algorithm}
\usepackage{multirow, booktabs}
\usepackage{siunitx}
\usepackage{xr}
\mathtoolsset{showonlyrefs=true}

\newlist{steps}{enumerate}{1}
\setlist[steps, 1]{label = Step \arabic*:}

\title{Theory and numerics of subspace approximation of eigenvalue problems}

\author{
  Siu Wun Cheung\thanks{Center for Applied Scientific Computing, Lawrence
  Livermore National Laboratory, Livermore, CA 94550 (cheung26@llnl.gov)}
  \and
  Youngsoo Choi\thanks{Center for Applied Scientific Computing, Lawrence
  Livermore National Laboratory, Livermore, CA 94550 (choi15@llnl.gov)}
  \and
  Seung Whan Chung\thanks{Center for Applied Scientific Computing, Lawrence
  Livermore National Laboratory, Livermore, CA 94550 (chung28@llnl.gov)}
  \and
  Jean-Luc Fattebert\thanks{Computational Sciences and Engineering Division, 
  Oak Ridge National Laboratory, Oak Ridge, TN 37830 (fattebertj@ornl.gov)}
  \and
  Coleman Kendrick\thanks{Center for Applied Scientific Computing, Lawrence
  Livermore National Laboratory, Livermore, CA 94550 (kendrick6@llnl.gov)}
  \and
  Daniel Osei-Kuffuor\thanks{Center for Applied Scientific Computing, Lawrence
  Livermore National Laboratory, Livermore, CA 94550 (oseikuffuor1@llnl.gov)}
}

\begin{document}
\setlength{\abovedisplayskip}{3pt}
\setlength{\belowdisplayskip}{3pt} 
\setlength{\abovedisplayshortskip}{3pt} 
\setlength{\belowdisplayshortskip}{3pt}

\maketitle

\begin{abstract}
Large-scale eigenvalue problems arise in various fields of science and engineering and demand computationally efficient solutions.
In this study, we investigate the subspace approximation for parametric linear eigenvalue problems, aiming to mitigate the computational burden associated with high-fidelity systems.
We provide general error estimates under non-simple eigenvalue conditions, \rev{establishing} \revv{some} \rev{theoretical foundations for understanding the convergence behavior of subspace approximations}.
Numerical examples, \rev{including problems with} one-dimensional to three-dimensional \rev{spatial domain and one-dimensional to two-dimensional parameter domain}, are presented to demonstrate the efficacy of reduced basis method in handling parametric variations in boundary conditions and coefficient fields to achieve significant computational savings while maintaining high accuracy, making them promising tools for practical applications in large-scale eigenvalue computations.
\end{abstract}

\begin{keywords} 
reduced order model, eigenvalue problems
\end{keywords}

\input{introduction}

\input{formulation}

\input{ROM}
\input{analysis}
\input{numerical}

\input{conclusion}

\bibliographystyle{unsrt}
\bibliography{references}

\appendix
\input{appendix}

\end{document}

%% file: introduction.tex
\section{Introduction}

Eigenvalue problems are of significant importance in various scientific and engineering disciplines, 
as they provide critical insights into the behavior and properties of complex systems. 
Applications range from structural dynamics, where eigenvalue analysis is used for modal analysis and vibration analysis \cite{rao2019vibration}, 
to quantum mechanics, where the Kohn-Sham equation is solved as an eigenvalue problem to determine energy levels and wave-functions \cite{CANCES20033}.
In many cases, the analytical solution of such problems can not be obtained, 
and numerical methods can be used to approximate the solutions. 
However, the complexity and size of these problems can lead to significant computational demands, 
especially as the fidelity and resolution of models increase. 
It may take a long time to solve a problem even with high performance computing and 
efficient iterative methods such as the Locally Optimal Block Preconditioned Conjugate Gradient (LOBPCG) \cite{knyazev2001toward}. 
Consequently, reducing computational cost of eigenvalue problems without sacrificing much accuracy 
has become a crucial goal in scientific computing and engineering simulations. 

One promising approach to achieve this goal is through reduced order modeling (ROM). 
Reduced order models provide efficient and accurate approximations of large-scale systems by utilizing low-dimensional representations of the solution manifold, thus significantly reducing the dimensionality of the problem. 
The low-dimensional representations can be constructed by compressing high-fidelity snapshot data through linear techniques such as proper orthogonal decomposition (POD) \cite{berkooz1993proper}, balanced truncation \cite{safonov1989schur}, and the reduced basis method \cite{rozza2008reduced}, as well as nonlinear approaches like autoencoders (AE) \cite{lee2020model,maulik2021reduced,kim2022fast}.
Reduced systems are typically obtained by projecting the original large-scale systems onto these low-dimensional structures. 
Projection-based ROM techniques directly incorporate the reduced solution representations into the governing equations and numerical discretization frameworks, ensuring that the ROMs honor the physical principles underlying the original problem. 
This enhances their reliability and accuracy while requiring less data to achieve similar levels of precision. 
Additionally, hyper-reduction techniques can be employed to further reduce the computational complexity associated with evaluating the nonlinear terms in the governing equations \cite{chaturantabut2010nonlinear, drmac2018discrete, lauzon2024s}. 
Projection-based ROMs have demonstrated success across a wide range of time-dependent nonlinear problems, including wave equations \cite{fares2011reduced, cheng2016reduced, cheung2021explicit}, the Burgers equation \cite{choi2019space, choi2020sns, carlberg2018conservative}, the Euler equations \cite{copeland2022reduced, cheung2023local}, porous media flow \cite{ghasemi2015localized, cheung2020constraint}, and Boltzmann transport problems \cite{choi2021space}. 
Detailed surveys on classical projection-based ROM techniques are available in \cite{gugercin2004survey, benner2015survey}.

In contrast to the extensive research on projection-based ROMs for time-dependent nonlinear problems 
\rev{and their widely spread applications across various scientific and engineering fields, 
our focus here is on the distinct field of parametric eigenvalue problems.}
Earlier works in this area include \cite{maday1999general, machiels2000output} 
on reduced-basis output bound method for symmetric eigenvalue problems. 
While their approach is effective for the first eigenpair, it can be restrictive in applications where \rev{non-simple eigenvalues, 
meaning eigenvalues with an algebraic multiplicity greater than one, lead to multiple eigenpairs of interest}.
\rev{More recently, the use of subspace methods for approximating parametric eigenvalue problems 
has gained significant traction. In \cite{kangal2018subspace}, the authors address the efficient min-max characterization 
of the $k$-th eigenvalue over a parametric domain. 
The development of certified greedy strategies for the approximation of the smallest eigenvalue, 
as seen in \cite{sirkovic2016subspace}, has also been instrumental. 
Recent advancements also include reduced basis approximation and a posteriori error estimates 
for parametrized elliptic eigenvalue problems \cite{fumagalli2016reduced}, 
tackling the challenges posed by clusters and intersections in parametric eigenvalue problems \cite{boffi2024reduced}, and 
the coupling of reduced basis methods with perturbation theory for eigenvalue problems is explored in \cite{garrigue2024reduced}. 
Furthermore, recent advancements in the uniform approximation of the smallest eigenvalue and 
singular value over a continuum parametric domain are presented in \cite{manucci2024uniform}, 
and the application of subspace methods for eigenvalue problems in pseudospectra computation is discussed in \cite{sirkovic2019reduced}. 
Applications to quantum spin systems are explored in \cite{herbst2022surrogate, brehmer2023reduced}. 
For non-symmetric cases, \cite{taumhas2024reduced} explores the reduced basis method for non-symmetric eigenvalue problems, 
with applications in neutron diffusion equations.}
In \cite{horger2017simultaneous}, the authors extended the reduced basis method to affinely parameterized elliptic eigenvalue problems, 
aiming to approximate several of the smallest eigenvalues simultaneously. 
The authors developed a posteriori error estimators for eigenvalues and presented various greedy strategies. 
In \cite{bertrand2023reduced}, the authors considered a reduced order method for approximating eigenfunctions of the Laplace problem 
with Dirichlet boundary conditions in a finite element setting. 
They used a time continuation technique to transform the problem into time-dependent and 
adopted a proper orthogonal decomposition (POD) approach for approximation of subsequent eigenmodes. 
They also presented theoretical results for choosing the optimal POD basis dimension.

In this study, we consider subspace approximation for solving parametric linear eigenvalue problems. 
We present a-priori error estimates of multiple eigenvalues and eigenvectors 
in a general setting of non-simple eigenvalues, 
providing rigorous theoretical foundations for the accuracy of the subspace approximation. 
We include several numerical examples arising from the finite element discretization of second-order elliptic differential equations, 
including Laplace problem, Schr\"{o}dinger equation, and heterogeneous diffusion problem. 
These examples range from one-dimensional to three-dimensional, 
and incorporate various parametric boundary conditions and coefficient fields.  
We demonstrate the capability of reduced basis method to handle complex parametric dependencies efficiently, 
significantly reducing computational costs while maintaining high accuracy.

\subsection{Organization}
In Section~\ref{sec:fom}, we describe the general setting of linear eigenvalue problems. 
Next, we introduce the subspace approximation in Section~\ref{sec:rom} and 
present an error analysis in Section~\ref{sec:analysis}. 
We provide various numerical examples of using the reduced basis methods for solving 
parametric linear eigenvalue problems in Section~\ref{sec:numerical}. 
Finally, a conclusion is given. 

\subsection{Notations}

Throughout this paper, we use boldfaced lowercase letters to denote scalar vectors, and boldfaced uppercase letters to denote matrices. 
We also follow the following notations: 

\begin{itemize}
\item $\mathbb{R}_{+}$ -- the set of all non-negative real numbers,
\item $\mathbb{R}_{++}$ -- the set of all positive real numbers,
\item $\mathbb{R}^{n}$ -- the space of all $n$-tuple real-valued vectors,
\item $\mathbb{R}^{m \times n}$ -- the space of all $m \times n$ real-valued matrices,
\item $\mathcal{R}(\mathbf{B})$ -- the column space of a matrix $\mathbf{B}$, 
\item $\mathbb{S}^{n}$ -- the set of all symmetric matrices in $\mathbb{R}^{n \times n}$,
\item $\mathbb{S}^{n}_{++}$ -- the set of all symmetric and positive definite matrices in $\mathbb{R}^{n \times n}$,
\item $\text{St}(m,n)$ -- the set of all orthogonal matrices in $\mathbb{R}^{m \times n}$, and 
\item $\text{Gr}_m(\mathcal{V})$ -- the set of all $m$-dimensional subspaces of a vector space $\mathcal{V}$. 
\end{itemize}

%% file: formulation.tex
\section{Problem statement}
\label{sec:fom}

Suppose $\mathbf{A} \in \mathbb{S}^{n}$ is a symmetric matrix, and 
$\mathbf{M} \in \mathbb{S}^{n}_{++}$ is a symmetric positive definite matrix. 
The problem of finding the non-trivial solutions of the equation 
\begin{equation}
    \label{eq:fom-evp}
    \mathbf{A} \boldsymbol{\phi} = \lambda \mathbf{M} \boldsymbol{\phi},  
\end{equation}
is called the generalized eigenvalue problem \cite{parlett1998symmetric}. 
Under this setting, 
$\boldsymbol{\phi} \in \mathbb{R}^{n} \setminus \{\mathbf{0}\}$ satisfying 
is called a generalized eigenvector of $\mathbf{A}$ 
with respect to $\mathbf{M}$, and 
$\lambda \in \mathbb{R}$ is the associated generalized eigenvalue of $\mathbf{A}$ 
with respect to $\mathbf{M}$. 
By the spectral theorem of matrices, there exists an $\mathbf{M}$-orthonormal basis 
$\{ \boldsymbol{\phi}_k \}_{k=1}^{n} \subset \mathbb{R}^{n}$ 
consisting of generalized eigenvectors of $\mathbf{A}$ with respect to $\mathbf{M}$, 
and the associated eigenvalues of $\mathbf{A}$ with respect to $\mathbf{M}$ are real 
and assumed to be arranged in ascending order, i.e. 
\begin{equation}
\label{eq:fom-eig}
    \lambda_1 \leq \lambda_2 \leq \cdots \leq 
    \lambda_k \leq \lambda_{k+1} \leq \cdots \leq \lambda_n.
\end{equation}
In a compact manner, we write the generalized eigenvalue decomposition as 
\begin{equation}
    \label{eq:fom-evd}
    \mathbf{A} \boldsymbol{\Phi} = \mathbf{M} \boldsymbol{\Phi} \boldsymbol{\Lambda},
\end{equation}
where $\boldsymbol{\Lambda} = \text{diag}(\lambda_1, \lambda_2, \cdots, \lambda_n)$ 
is the diagonal matrix consisting of the generalized eigenvalues of $\mathbf{A}$ with respect to $\mathbf{M}$, 
and $\boldsymbol{\Phi} = \left[ \boldsymbol{\phi}_1, \boldsymbol{\phi}_2, \cdots, \boldsymbol{\phi}_n \right] \in \mathbb{R}^{n \times n}$ 
is the matrix consisting of the generalized eigenvectors of $\mathbf{A}$ with respect to $\mathbf{M}$ 
satisfying $\boldsymbol{\Phi}^\top \mathbf{M} \boldsymbol{\Phi} = \mathbf{I}$. 

Furthermore, we assume that the matrices $\mathbf{A}(\boldsymbol{\mu})$ and $\mathbf{M}(\boldsymbol{\mu})$
depend on a set of problem parameters $\boldsymbol{\mu} \in \mathsf{D}$, 
\rev{where the parameter domain $\mathsf{D}$ is a compact subset of $\mathbb{R}^{d_{\mu}}$}.
As a result, the eigendecomposition
$(\boldsymbol{\Lambda}(\boldsymbol{\mu}), \boldsymbol{\Phi}(\boldsymbol{\mu}))$
\rev{is} also parametrized by $\boldsymbol{\mu}$.
\rev{
We remark that it is a common assumption in subspace approximation methods 
that the matrices exhibit an affine-parametric dependence. 
This means that $\mathbf{A}(\boldsymbol{\mu})$ and $\mathbf{M}(\boldsymbol{\mu})$ 
can be expressed as a linear combination of parameter-independent operators:
\begin{equation}
\begin{split}
    \mathbf{A}(\boldsymbol{\mu}) & = \sum_{q=1}^{N_A} \theta_A^q(\boldsymbol{\mu}) \mathbf{A}_q, \\
    \mathbf{M}(\boldsymbol{\mu}) & = \sum_{q=1}^{N_M} \theta_M^q(\boldsymbol{\mu}) \mathbf{M}_q,
\end{split}
\label{eq:affine_param}
\end{equation}
where $\mathbf{A}_q$ and $\mathbf{M}_q$ are parameter-independent matrices,
and $\theta_A^q(\boldsymbol{\mu})$ and $\theta_M^q(\boldsymbol{\mu})$ are scalar } 
\revv{analytic}
\rev{ functions 
that depend only on the parameter $\boldsymbol{\mu}$. 
While this assumption is satisfied by many parameterized physical problems after discretization, 
some other problem formulations, such as those involving atomic well potentials with parameterized positions, 
do not inherently satisfy this assumption.
}

%% file: ROM.tex
\section{Subspace approximation}
\label{sec:rom}

In this section, we present the projected eigenvalue problem onto
the column space of an orthogonal matrix $\mathbf{Q} \in \text{St}(n,r)$ with $r \leq n$,
which is referred to as the basis matrix.
Let \rev{$\widehat{\mathbf{A}} := \mathbf{Q}^\top \mathbf{A} \mathbf{Q} \in \mathbb{S}^{r}$}
and \rev{$\widehat{\mathbf{M}} := \mathbf{Q}^\top \mathbf{M} \mathbf{Q} \in \mathbb{S}_{++}^{r}$}
be the projection of $\mathbf{A}$ and $\mathbf{M}$ onto the column space $\mathcal{R}(\mathbf{Q})$ of $\mathbf{Q}$ respectively.
\rev{For conciseness in notation, we temporarily drop the explicit dependence on the parameter $\boldsymbol{\mu}$ from $\mathbf{A}$, $\mathbf{M}$, and their projected counterparts ($\widehat{\mathbf{A}}$, $\widehat{\mathbf{M}}$). It is understood that these operators, and consequently the eigenvalues and eigenvectors, inherently depend on $\boldsymbol{\mu}$, as introduced in Section 2.}
Analogously, the projected generalized eigenvalue problem is given by
\begin{equation}
    \label{eq:rom-evp}
    \widehat{\mathbf{A}} \widehat{\boldsymbol{\phi}} = \widetilde{\lambda} \widehat{\mathbf{M}} \widehat{\boldsymbol{\phi}},
\end{equation}
where $\widehat{\boldsymbol{\phi}} \in \mathbb{R}^{r} \setminus \{\mathbf{0}\}$
is called a generalized eigenvector of $\widehat{\mathbf{A}}$
with respect to $\widehat{\mathbf{M}}$, and
$\widetilde{\lambda} \in \mathbb{R}$ is the associated generalized eigenvalue of $\widehat{\mathbf{A}}$ with respect to $\widehat{\mathbf{M}}$.
Again, by the spectral theorem of symmetric matrices, there exists an $\widehat{\mathbf{M}}$-orthonormal basis
$\{ \widehat{\boldsymbol{\phi}}_k\}_{k=1}^{r} \subset \mathbb{R}^{r}$
consisting of generalized eigenvectors of $\widehat{\mathbf{A}}$ with respect to $\widehat{\mathbf{M}}$,
and the associated eigenvalues $\{\widetilde{\lambda}_k\}_{k=1}^{r}$ of $\widehat{\mathbf{A}}$ with respect to $\widehat{\mathbf{M}}$ are real
and assumed to be arranged in ascending order, i.e.
\begin{equation}
    \widetilde{\lambda}_1 \leq \widetilde{\lambda}_2 \leq \cdots \leq \widetilde{\lambda}_k \leq \widetilde{\lambda}_{k+1} \leq \cdots \leq \widetilde{\lambda}_r.
\end{equation}
The projected generalized eigenvalue decomposition can be written as
\begin{equation}
\label{eq:rom-evd}
    \widehat{\mathbf{A}} \widehat{\boldsymbol{\Phi}} = \widehat{\mathbf{M}} \widehat{\boldsymbol{\Phi}} \widetilde{\boldsymbol{\Lambda}},
\end{equation}
where $\widetilde{\boldsymbol{\Lambda}} = \text{diag}(\widetilde{\lambda}_1, \widetilde{\lambda}_2, \cdots, \widetilde{\lambda}_r)$
is the diagonal matrix consisting of the generalized eigenvalues of $\widehat{\mathbf{A}}$ with respect to $\widehat{\mathbf{M}}$, and
$\widehat{\boldsymbol{\Phi}} = \left[ \widehat{\boldsymbol{\phi}}_1, \widehat{\boldsymbol{\phi}}_2, \cdots, \widehat{\boldsymbol{\phi}}_r \right] \in \mathbb{R}^{r \times r}$
is the matrix consisting of the generalized eigenvectors of $\widehat{\mathbf{A}}$ with respect to $\widehat{\mathbf{M}}$
satisfying $\widehat{\boldsymbol{\Phi}}^\top \widehat{\mathbf{M}} \widehat{\boldsymbol{\Phi}} = \mathbf{I}$.
The subspace approximation of the first $r$ eigenvectors is then given by
$\widetilde{\boldsymbol{\Phi}} = \mathbf{Q} \widehat{\boldsymbol{\Phi}} \in \mathbb{R}^{n \times r}$,
which forms a basis for $\mathcal{R}(\mathbf{Q})$ and
satisfies $\widetilde{\boldsymbol{\Phi}}^\top \mathbf{M} \widetilde{\boldsymbol{\Phi}} = \mathbf{I}$.

%% file: analysis.tex
\section{Analysis}
\label{sec:analysis}

In this section, we provide an error analysis for the projected eigenvalue problem \eqref{eq:rom-evp}, 
in terms of a bound on the ratio between $\lambda_k$ and $\widetilde{\lambda}_k$, 
and the best approximation error of $\boldsymbol{\phi}_k$ by the projected eigenspace of the same eigenvalue 
in the $\mathbf{M}$-induced norm, which is given by 
$\| \mathbf{x} \|_{\mathbf{M}} = \sqrt{\mathbf{x}^\top \mathbf{M} \mathbf{x}}$, 
for any $\mathbf{x} \in \mathbb{R}^n$. 
This norm is associated with the $\mathbf{M}$-induced inner product, 
$(\mathbf{x}, \mathbf{y}) \in \mathbb{R}^n \times \mathbb{R}^n \mapsto \mathbf{x}^\top \mathbf{M} \mathbf{y} \in \mathbb{R}$, 
which satisfies the Cauchy-Schwarz inequality: 
\begin{equation}
    \mathbf{x}^\top \mathbf{M} \mathbf{y} \leq \| \mathbf{x} \|_{\mathbf{M}} \| \mathbf{y} \|_{\mathbf{M}} \text{ for } \mathbf{x}, \mathbf{y} \in \mathbb{R}^n.
\end{equation}

To simplify the notations of eigenpair indices,
given an index set $\mathcal{S}$ with cardinality $s = \vert \mathcal{S} \vert$, 
we denote the matrix consisting the generalized eigenvectors of $\mathbf{A}$ with respect to $\mathbf{M}$ indexed with $\mathcal{S}$ by 
$\boldsymbol{\Phi}^{(\mathcal{S})} = \left[ \boldsymbol{\phi}_k \right]_{k \in \mathcal{S}} \in \mathbb{R}^{n \times s}$.
In particular, for two positive integers $m_1 < m_2$, 
we denote the range by $\mathbb{N}(m_1,m_2) = \{m_1, m_1+1, \cdots, m_2\}$,
and $\boldsymbol{\Phi}^{(m_1,m_2)} = \boldsymbol{\Phi}^{(\mathbb{N}(m_1,m_2))}$. 
Furthermore, for any integer $k > 1$, 
we denote $\mathbb{N}(k) = \mathbb{N}(1,k)$ and $\boldsymbol{\Phi}^{(k)} = \boldsymbol{\Phi}^{(1,k)}$. 

We start with some technical results on the generalized Rayleigh quotient with respect to $\mathbf{A}$ and $\mathbf{M}$,
which is defined as: for $\mathbf{x} \in \mathbb{R}^{n} \setminus \{\mathbf{0}\}$,  
\begin{equation}
    R_{\mathbf{A}, \mathbf{M}}(\mathbf{x}) = \dfrac{\mathbf{x}^\top \mathbf{A} \mathbf{x}}{\mathbf{x}^\top \mathbf{M} \mathbf{x}}.
\end{equation}

\begin{lemma}
\label{lemma:rayleigh-quotient-fom}
    Let $m_1 < m_2$ be two positive integers. There hold 
    \begin{equation}
        \begin{split}
            \min_{\mathbf{x} \in \mathcal{R}(\boldsymbol{\Phi}^{(m_1,m_2)}) \setminus \{\mathbf{0}\}}
            R_{\mathbf{A}, \mathbf{M}}(\mathbf{x}) & = \lambda_{m_1} = R_{\mathbf{A}, \mathbf{M}}(\boldsymbol{\phi}_{m_1}), \\
            \max_{\mathbf{x} \in \mathcal{R}(\boldsymbol{\Phi}^{(m_1,m_2)}) \setminus \{\mathbf{0}\}}
            R_{\mathbf{A}, \mathbf{M}}(\mathbf{x}) & = \lambda_{m_2} = R_{\mathbf{A}, \mathbf{M}}(\boldsymbol{\phi}_{m_2}).
        \end{split}
    \end{equation}
    \begin{proof}
    Let $\mathbf{x} \in \mathcal{R}(\boldsymbol{\Phi}^{(m_1,m_2)})$. 
    By the properties of $\mathbf{M}$-orthonormal basis, 
    \begin{equation}
        \mathbf{x} = \sum_{k=m_1}^{m_2} (\boldsymbol{\phi}_k^\top \mathbf{M} \mathbf{x}) \boldsymbol{\phi}_k.
    \end{equation}
    By the properties of generalized eigenvectors, 
    \begin{equation}
        \mathbf{A} \mathbf{x} = \sum_{k=m_1}^{m_2} (\lambda_k \boldsymbol{\phi}_k^\top \mathbf{M} \mathbf{x}) \mathbf{M} \boldsymbol{\phi}_k.
    \end{equation}
    Again using the properties of $\mathbf{M}$-orthonormal basis,
    \begin{equation}
        R_{\mathbf{A}, \mathbf{M}}(\mathbf{x}) 
        = \dfrac{\mathbf{x}^\top \mathbf{A} \mathbf{x}}{\mathbf{x}^\top \mathbf{M} \mathbf{x}}
        = \dfrac{\sum_{k=m_1}^{m_2} \lambda_k (\boldsymbol{\phi}_k^\top \mathbf{M} \mathbf{x})^2}{\sum_{k=m_1}^{m_2} (\boldsymbol{\phi}_k^\top \mathbf{M} \mathbf{x})^2}.
    \end{equation}
    By the ordering of the generalized eigenvalues $\lambda_{m_1} \leq \lambda_k \leq \lambda_{m_2}$ for $k \in \mathbb{N}(m_1,m_2)$, we infer that 
    $\lambda_{m_1} \leq R_{\mathbf{A}, \mathbf{M}}(\mathbf{x}) \leq \lambda_{m_2}$.
    \end{proof}
\end{lemma}

The generalized eigenvalues of $\mathbf{A}$ with respect to $\mathbf{M}$ is characterized by the generalized Rayleigh quotient through the following min-max theorem.
\begin{theorem}[Min-max]
\label{thm:min-max}
    For $1 \leq k \leq n$, there holds 
    \begin{equation}
        \lambda_k 
        = \min_{\mathcal{X} \in \text{Gr}_k(\mathbb{R}^n)} 
        \max_{\mathbf{x} \in \mathcal{X} \setminus \{\mathbf{0}\} }  R_{\mathbf{A}, \mathbf{M}}(\mathbf{x}). 
    \end{equation}
    \begin{proof}
        Let $\mathcal{X} \in \text{Gr}_k(\mathbb{R}^n)$. 
        By rank-nullity theorem, $\text{dim}(\mathcal{X} \cap \mathcal{R}(\boldsymbol{\Phi}^{(k-1)})^\perp) \geq 1$. 
        Take $\mathbf{x}_k \in \mathcal{X} \cap \mathcal{R}(\boldsymbol{\Phi}^{(k-1)})^\perp \setminus \{\mathbf{0}\} = \mathcal{X} \cap \mathcal{R}(\boldsymbol{\Phi}^{k,n}) \setminus \{\mathbf{0}\}$. 
        By Lemma~\ref{lemma:rayleigh-quotient-fom}, 
        \begin{equation}
            \lambda_m \leq R_{\mathbf{A}, \mathbf{M}}(\mathbf{x}_k) \leq \max_{\mathbf{x} \in \mathcal{X} \setminus \{\mathbf{0}\} }  R_{\mathbf{A}, \mathbf{M}}(\mathbf{x}).
        \end{equation}
        Since $\mathcal{X} \in \text{Gr}_k(\mathbb{R}^n)$ is arbitrary, we obtain 
        \begin{equation}
            \lambda_k 
            \leq \inf_{\mathcal{X} \in \text{Gr}_k(\mathbb{R}^n)} 
            \max_{\mathbf{x} \in \mathcal{X} \setminus \{\mathbf{0}\} }  
            R_{\mathbf{A}, \mathbf{M}}(\mathbf{x}). 
        \end{equation}
        On the other hand, since $\mathcal{R}(\boldsymbol{\Phi}^{(k)}) \in \text{Gr}_k(\mathbb{R}^n)$, again by Lemma~\ref{lemma:rayleigh-quotient-fom}, we have 
        \begin{equation}
            \inf_{\mathcal{X} \in \text{Gr}_k(\mathbb{R}^n)} 
            \max_{\mathbf{x} \in \mathcal{X} \setminus \{\mathbf{0}\} }  
            R_{\mathbf{A}, \mathbf{M}}(\mathbf{x}) 
            \leq \max_{\mathbf{x} \in \mathcal{R}(\boldsymbol{\Phi}^{(k)}) \setminus \{\mathbf{0}\} }  
            R_{\mathbf{A}, \mathbf{M}}(\mathbf{x}) = \lambda_k.
        \end{equation}
        This completes the proof.
    \end{proof}
\end{theorem}

With the simplified index notations extended to $\widetilde{\boldsymbol{\Phi}}$ and $\widehat{\boldsymbol{\Phi}}$, 
and the definition of generalized Rayleigh-quotient extended to $(\widehat{\mathbf{A}}, \widehat{\mathbf{M}})$, 
we have an analogous result for the approximation of generalized eigenvalues $\widetilde{\lambda}_k$. 

\begin{corollary}
\label{corollary:rayleigh-quotient-rom}
    Let $m_1 < m_2$ be two positive integers. There hold 
    \begin{equation}
        \begin{split}
            \min_{\mathbf{x} \in \mathcal{R}(\widetilde{\boldsymbol{\Phi}}^{(m_1,m_2)}) \setminus \{\mathbf{0}\}}
            R_{\mathbf{A}, \mathbf{M}}(\mathbf{x}) & = \widetilde{\lambda}_{m_1} = R_{\mathbf{A}, \mathbf{M}}(\widetilde{\boldsymbol{\phi}}_{m_1}), \\
            \max_{\mathbf{x} \in \mathcal{R}(\widetilde{\boldsymbol{\Phi}}^{(m_1,m_2)}) \setminus \{\mathbf{0}\}}
            R_{\mathbf{A}, \mathbf{M}}(\mathbf{x}) & = \widetilde{\lambda}_{m_2} = R_{\mathbf{A}, \mathbf{M}}(\widetilde{\boldsymbol{\phi}}_{m_2}).
        \end{split}
    \end{equation}
    \begin{proof}
    Let $\mathbf{x} \in \mathcal{R}(\widetilde{\boldsymbol{\Phi}}^{(m_1,m_2)}) \setminus \{ \mathbf{0} \}$. 
    Take $\widehat{\mathbf{x}} = \mathbf{Q}^\top \mathbf{x} \in \mathcal{R}(\widehat{\boldsymbol{\Phi}}^{(m_1,m_2)}) \setminus \{ \mathbf{0} \}$. 
    Then we have $\mathbf{x} = \mathbf{Q} \widehat{\mathbf{x}}$, and
    \begin{equation}
        R_{\mathbf{A}, \mathbf{M}}(\mathbf{x})
        = \dfrac{(\mathbf{Q}\widehat{\mathbf{x}})^\top \mathbf{A} (\mathbf{Q}\widehat{\mathbf{x}})}{(\mathbf{Q}\widehat{\mathbf{x}})^\top \mathbf{M} (\mathbf{Q}\widehat{\mathbf{x}})}
        = \dfrac{\widehat{\mathbf{x}}^\top \widehat{\mathbf{A}} \widehat{\mathbf{x}}}{\widehat{\mathbf{x}}^\top \widehat{\mathbf{M}} \widehat{\mathbf{x}}}
        = R_{\widehat{\mathbf{A}},\widehat{\mathbf{M}}}(\widehat{\mathbf{x}}).
    \end{equation}
    As a direct consequence of Lemma~\ref{lemma:rayleigh-quotient-fom}, we obtain
    \begin{equation}
        \widetilde{\lambda}_{m_1} \leq R_{\widehat{\mathbf{A}},\widehat{\mathbf{M}}}(\widehat{\mathbf{x}}) \leq \widetilde{\lambda}_{m_2}.
    \end{equation}
    \end{proof}
\end{corollary}

In what follows, we assume $\mathbf{A} \in \mathbb{S}^n_{++}$ without loss of generality, 
or otherwise we can replace $\mathbf{A}$ by $\mathbf{A} + t \mathbf{M}$ with $t > -\lambda_1$.
This allows us to define the $\mathbf{A}$-induced norm by 
$\| \mathbf{x} \|_{\mathbf{A}} = \sqrt{\mathbf{x}^\top \mathbf{A} \mathbf{x}}$. 
This norm is associated with the $\mathbf{A}$-induced inner product, 
$(\mathbf{x}, \mathbf{y}) \in \mathbb{R}^n \times \mathbb{R}^n \mapsto \mathbf{x}^\top \mathbf{A} \mathbf{y} \in \mathbb{R}$, 
which also satisfies the Cauchy-Schwarz inequality: 
\begin{equation}
    \mathbf{x}^\top \mathbf{A} \mathbf{y} \leq \| \mathbf{x} \|_{\mathbf{A}} \| \mathbf{y} \|_{\mathbf{A}} \text{ for } \mathbf{x}, \mathbf{y} \in \mathbb{R}^n.
\end{equation}

To facilitate the analysis, we introduce the oblique projection onto $\mathcal{R}(\mathbf{Q})$ with the matrix $\mathbf{A}$, 
defined by $\mathbf{P}_{\mathbf{A}} = \mathbf{Q} \widehat{\mathbf{A}}^{-1} \mathbf{Q}^\top \mathbf{A} \in \mathbb{R}^{n \times n}$.
\rev{
Since $\mathbf{A} \in \mathbb{S}^n_{++}$ and $\mathbf{Q} \in \text{St}(n,r)$ has full column rank, 
the projected matrix $\widehat{\mathbf{A}} = \mathbf{Q}^\top \mathbf{A} \mathbf{Q}$ is also positive definite. 
This ensures that the oblique projection $\mathbf{P}_{\mathbf{A}}$ is well-defined.
}
The following lemma states that $\mathbf{P}_{\mathbf{A}}$ is non-expansive in the $\mathbf{A}$-induced norm.
\begin{lemma}
\label{lemma:A-projection-estimate}
    For any $\mathbf{x} \in \mathbb{R}^n$, 
    $\| \mathbf{P}_{\mathbf{A}} \mathbf{x} \|_{\mathbf{A}} \leq \| \mathbf{x} \|_{\mathbf{A}}$. 
    \begin{proof}
    It is straightforward to check that 
    \begin{equation}
        \mathbf{Q}^\top \mathbf{A} (\mathbf{I} - \mathbf{P}_{\mathbf{A}}) = \mathbf{0}.
    \label{eq:A-galerkin-orthogonality}
    \end{equation}
    As a direct consequence, 
    we have $(\mathbf{P}_{\mathbf{A}} \mathbf{x})^\top \mathbf{A} (\mathbf{I} - \mathbf{P}_{\mathbf{A}}) \mathbf{x} = \mathbf{0}$, and hence
    \begin{equation}
    \begin{split}
        \| \mathbf{P}_{\mathbf{A}} \mathbf{x} \|_{\mathbf{A}}^2
        & = (\mathbf{P}_{\mathbf{A}} \mathbf{x})^\top \mathbf{A} (\mathbf{P}_{\mathbf{A}} \mathbf{x}) \\
        & = (\mathbf{P}_{\mathbf{A}} \mathbf{x})^\top \mathbf{A} \mathbf{x} \\
        & \leq \| \mathbf{P}_{\mathbf{A}} \mathbf{x} \|_{\mathbf{A}} \| \mathbf{x} \|_{\mathbf{A}}, 
    \end{split}
    \end{equation}
    where we have used the Cauchy-Schwarz inequality for the $\mathbf{A}$-induced inner product.
    \end{proof}
\end{lemma}

Similarly, given an index set $\mathcal{S}$, we define the oblique projection onto 
$\mathcal{R}(\widetilde{\boldsymbol{\Phi}}^{(\mathcal{S})})$
with the matrix $\mathbf{M}$ respectively by 
\begin{equation}
\widetilde{\mathbf{P}}_{\mathbf{M}}^{(\mathcal{S})} = \widetilde{\boldsymbol{\Phi}}^{(\mathcal{S})} 
\left( \widetilde{\boldsymbol{\Phi}}^{(\mathcal{S})} \right)^\top 
\mathbf{M} \in \mathbb{R}^{n \times n}.
\end{equation}
The following lemma states that $\widetilde{\mathbf{P}}_{\mathbf{M}}^{(\mathcal{S})}$ maps to the best approximation in
$\mathcal{R}(\widetilde{\boldsymbol{\Phi}}^{(\mathcal{S})})$ with respect to the $\mathbf{M}$-induced norm.
\begin{lemma}
\label{lemma:M-projection-estimate}
    For any $\mathbf{x} \in \mathbb{R}^n$, 
    $\| (\mathbf{I} - \widetilde{\mathbf{P}}_{\mathbf{M}}^{(\mathcal{S})}) \mathbf{x} \|_{\mathbf{M}} = 
    \min_{\mathbf{y} \in \mathcal{R}(\widetilde{\boldsymbol{\Phi}}^{(\mathcal{S})})} \| \mathbf{x} - \mathbf{y} \|_{\mathbf{M}}$.
    \begin{proof}
    It is straightforward to check that 
    \begin{equation}
        \left( \widetilde{\mathbf{P}}_{\mathbf{M}}^{(\mathcal{S})} \right)^\top \mathbf{M} (\mathbf{I} - \widetilde{\mathbf{P}}_{\mathbf{M}}^{(\mathcal{S})}) = \mathbf{0}.
    \label{eq:M-galerkin-orthogonality}
    \end{equation}
    As a direct consequence, for any $\mathbf{y} \in \mathcal{R}(\widetilde{\boldsymbol{\Phi}}^{(\mathcal{S})})$, we have 
    $(\mathbf{y} - \widetilde{\mathbf{P}}_{\mathbf{M}}^{(\mathcal{S})} \mathbf{x})^\top \mathbf{M} (\mathbf{I} - \widetilde{\mathbf{P}}_{\mathbf{M}}^{(\mathcal{S})}) \mathbf{x} = \mathbf{0}$, and therefore 
    \begin{equation}
    \begin{split}
        \| (\mathbf{I} - \widetilde{\mathbf{P}}_{\mathbf{M}}^{(\mathcal{S})}) \mathbf{x} \|_{\mathbf{M}}^2
        & = (\mathbf{x} - \widetilde{\mathbf{P}}_{\mathbf{M}}^{(\mathcal{S})} \mathbf{x})^\top \mathbf{M} (\mathbf{I} - \widetilde{\mathbf{P}}_{\mathbf{M}}^{(\mathcal{S})}) \mathbf{x} \\
        & = (\mathbf{x} - \mathbf{y})^\top \mathbf{M} (\mathbf{I} - \widetilde{\mathbf{P}}_{\mathbf{M}}^{(\mathcal{S})}) \mathbf{x} \\
        & \leq \| \mathbf{x} - \mathbf{y} \|_{\mathbf{M}} \| (\mathbf{I} - \widetilde{\mathbf{P}}_{\mathbf{M}}^{(\mathcal{S})}) \mathbf{x} \|_{\mathbf{M}},
    \end{split}
    \end{equation}
    where we have used the Cauchy-Schwarz inequality for the $\mathbf{M}$-induced inner product. 
    Since $\mathbf{y} \in \mathcal{R}(\widetilde{\boldsymbol{\Phi}}^{(\mathcal{S})})$ is arbitrary, we arrive at the desired result.
    \end{proof}
\end{lemma}

We are now ready to present an error analysis of the generalized eigenvalues. 
We assume that the matrix $\mathbf{P}_{\mathbf{A}}\boldsymbol{\Phi}^{(r)} \in \mathbb{R}^{n \times r}$ 
is of full rank $r$, which assure that, for $1 \leq k \leq r$,
\begin{equation}
    \kappa_k = \sup_{\mathbf{y} \in \mathcal{R}(\boldsymbol{\Phi}^{(k)})} \dfrac{\| \mathbf{y} \|_{\mathbf{M}}}{\| \mathbf{P}_{\mathbf{A}} \mathbf{y} \|_{\mathbf{M}}} < \infty.
    \label{eq:condition_number}
\end{equation}
\rev{The assumption that $\mathbf{P}_{\mathbf{A}}\boldsymbol{\Phi}^{(r)}$ is of full rank $r$ 
implies that the oblique projection of the true eigenvector subspace onto the reduced basis $\mathcal{R}(\mathbf{Q})$ is not rank-deficient. 
It ensures the reduced operator $\widehat{\mathbf{A}}$ is effectively capturing the components of the full-order eigenvectors 
$\boldsymbol{\Phi}^{(r)}$ that are relevant for the approximation.
However, this assumption is not always guaranteed. 
For a certain parameter $\boldsymbol{\mu} \in \textsf{D}$, if a vector $\mathbf{y} \in \mathcal{R}(\boldsymbol{\Phi}^{(r)})$ 
spanned by the first $r$ eigenvectors is orthogonal to the column space $\mathcal{R}(\mathbf{Q})$ of the basis matrix $\mathbf{Q}$, 
then the projection $\mathbf{P}_{\mathbf{A}}\mathbf{y}$ would be zero, leading to $\mathbf{P}_{\mathbf{A}}\boldsymbol{\Phi}^{(r)}$ becoming rank-deficient.
This scenario can occur in practice, especially if the basis $\mathbf{Q}$ is not sufficiently rich or optimally constructed 
to capture the behavior of the true eigenvectors across the parameter domain.
In practice, the construction of the basis $\mathbf{Q}$ should ideally ensure that the chosen basis vectors are well-aligned with the solution manifold, 
thereby minimizing the likelihood of such near-orthogonality and maintaining the full rank of $\mathbf{P}_{\mathbf{A}}\boldsymbol{\Phi}^{(r)}$.}

\begin{theorem}
\label{thm:eigval-error}
    For $1 \leq k \leq r$, there holds 
    \begin{equation}
        \lambda_k \leq \widetilde{\lambda}_k \leq \kappa_k^{2} \lambda_k. 
    \end{equation}
    \begin{proof}
        First, by Lemma~\ref{corollary:rayleigh-quotient-rom}, we note that 
        \begin{equation}
            \widetilde{\lambda}_k 
            = \max_{\mathbf{x} \in \mathcal{R}(\widetilde{\Phi}_k) \setminus \{ \mathbf{0} \}}
            R_{\mathbf{A}, \mathbf{M}}(\mathbf{x}).
        \end{equation}
        Since $\mathcal{R}(\widetilde{\boldsymbol{\Phi}}^{(k)}) \in \text{Gr}_m(\mathbb{R}^n)$, we conclude from Theorem~\ref{thm:min-max} that
        \begin{equation}
            \widetilde{\lambda}_k \geq 
            \min_{\mathcal{X} \in \text{Gr}_k(\mathbb{R}^n)} 
            \max_{\mathbf{x} \in \mathcal{X} \setminus \{\mathbf{0}\} }  R_{\mathbf{A}, \mathbf{M}}(\mathbf{x})
            = \lambda_k.
        \end{equation}
        On the other hand, by rank-nullity theorem, $\text{dim}(\mathcal{R}(\mathbf{P}_{\mathbf{A}}\boldsymbol{\Phi}^{(k)}) \cap \mathcal{R}(\widetilde{\boldsymbol{\Phi}}^{(k-1)})^\perp) \geq 1$. 
        Take $\mathbf{x}_k \in \mathcal{R}(\mathbf{P}_{\mathbf{A}}\boldsymbol{\Phi}^{(k)}) \cap \mathcal{R}(\widetilde{\boldsymbol{\Phi}}^{(k-1)})^\perp \setminus \{\mathbf{0}\} = \mathcal{R}(\mathbf{P}_{\mathbf{A}}\boldsymbol{\Phi}^{(k)}) \cap \mathcal{R}(\widetilde{\boldsymbol{\Phi}}^{(k,r)}) \setminus \{\mathbf{0}\}$. 
        By Lemma~\ref{corollary:rayleigh-quotient-rom}, $\widetilde{\lambda}_k \leq R_{\mathbf{A}, \mathbf{M}}(\mathbf{x}_k)$.
        Take $\mathbf{y}_k \in \mathcal{R}(\boldsymbol{\Phi}^{(k)})$ such that 
        $\mathbf{x}_k = \mathbf{P}_{\mathbf{A}} \mathbf{y}_k$. By Lemma~\ref{lemma:A-projection-estimate}, we have 
        \begin{equation}
        \begin{split}
            R_{\mathbf{A}, \mathbf{M}}(\mathbf{x}_k)
            & = \dfrac{\| \mathbf{P}_{\mathbf{A}} \mathbf{y}_k \|_{\mathbf{A}}^2}{\| \mathbf{P}_{\mathbf{A}} \mathbf{y}_k \|_{\mathbf{M}}^2} \\
            & \leq \dfrac{\| \mathbf{y}_k \|_{\mathbf{A}}^2}{\| \mathbf{P}_{\mathbf{A}} \mathbf{y}_k \|_{\mathbf{M}}^2} \\
            & = \dfrac{\| \mathbf{y}_k \|_{\mathbf{M}}^2}{\| \mathbf{P}_{\mathbf{A}} \mathbf{y}_k \|_{\mathbf{M}}^2} R_{\mathbf{A}, \mathbf{M}}(\mathbf{y}_k) \\
            & \leq \kappa_k^2 R_{\mathbf{A}, \mathbf{M}}(\mathbf{y}_k).
        \end{split}
        \end{equation}
        Again, by Lemma~\ref{lemma:rayleigh-quotient-fom}, since $\mathbf{y}_k \in \mathcal{R}(\boldsymbol{\Phi}^{(k)})$, we have 
        \rev{$$R_{\mathbf{A}, \mathbf{M}}(\mathbf{y}_k) \leq \max_{\mathbf{x} \in \mathcal{R}(\boldsymbol{\Phi}^{(k)}) \setminus \{\mathbf{0}\}}
            R_{\mathbf{A}, \mathbf{M}}(\mathbf{x}) = \lambda_k.$$}
        This completes the proof.
        \end{proof}
\end{theorem}

Next, we present an error analysis for the eigenspaces spanned by the generalized eigenvectors. 
\rev{It is important to note that the analysis focuses on the accuracy of these eigenspaces instead of individual eigenvectors.}
\revv{The main reason is that eigenvectors are not uniquely defined, even after normalization. For simple eigenvalues, they remain ambiguous up to a sign. For non-simple eigenvalues, the ambiguity is even greater, extending to arbitrary linear combinations within the eigenspace. In contrast, the corresponding eigenspace itself is uniquely defined. }
\rev{Moerover, individual eigenvectors can be highly sensitive to small perturbations. 
As a result, the distance between an eigenspace and its projection provides a more robust and physically meaningful measure of approximation quality.}
Let the distinct generalized eigenvalues of $\mathbf{A}$ with respect to $\mathbf{M}$ 
be arranged in ascending order, i.e. 
\begin{equation}
    \nu_1 < \nu_2 < \cdots < 
    \nu_j < \nu_{j+1} < \cdots < \nu_p,
\label{eq:fom-eig-distinct}
\end{equation}
and $\gamma_j$ be the geometric multiplicity of $\nu_j$ for $1 \leq j \leq p$. 
To classify the index in \eqref{eq:fom-eig} into eigenspaces, we denote 
$\Gamma_0 = 0$ and $\Gamma_j = \sum_{i=1}^{j} \gamma_i$ for $1 \leq j \leq p$. 
Then $\Gamma_p = n$. 
For simplicity, we assume that
$\lambda_r < \lambda_{r+1}$, which implies
$\Gamma_\ell = r$ for some $\ell \in \mathbb{N}(p)$.
For $1 \leq j \leq p$, we denote 
$\mathcal{S}_j = \mathbb{N}(\Gamma_{j-1} + 1, \Gamma_{j})$, 
which implies $\lambda_k = \nu_j$ for $k \in \mathcal{S}_j$, and 
the generalized eigenspace associated with $\nu_j$ is given by 
\begin{equation}
\mathcal{E}_j = \mathcal{R}(\boldsymbol{\Phi}^{(\mathcal{S}_j)}) 
= \{ \mathbf{x} \in \mathbf{R}^n \setminus \{0\} : \mathbf{A} \mathbf{x} = \nu_j \mathbf{M}\mathbf{x} \}.
\end{equation}
We also denote the complement of $\mathcal{S}_j$ in $\mathbb{N}(r)$ by 
$\mathcal{S}_j^{\text{c}} = \mathbb{N}(r) \setminus \mathcal{S}_j$.
We assume that for $k \in \mathcal{S}_j^{\text{c}}$,
$\nu_j \notin [\lambda_k, \kappa_k^2 \lambda_k]$, 
which assures that $\nu_j \neq \widetilde{\lambda}_k$ thanks to Theorem~\ref{thm:eigval-error}, 
and hence 
\begin{equation}
    \tau_j = \dfrac{\nu_j}{\min_{k \in \mathcal{S}_j^{\text{c}}} \vert \nu_j - \widetilde{\lambda}_{k} \vert} < \infty.
\end{equation}
\rev{The assumption that $\nu_j \notin [\lambda_k, \kappa_k^2 \lambda_k]$ is generally non-generic, 
since $\kappa_k$ can be large if the basis matrix $\mathbf{Q}$ does not provide a 
sufficiently good approximation space for the $k$-th eigenvector. 
When $\kappa_k$ is large, the interval $[\lambda_k, \kappa_k^2 \lambda_k]$ can become very wide, 
increasing the likelihood that $\nu_j$ falls within it. 
This scenario, where the assumption is violated, can occur in practice, especially if the basis $\mathbf{Q}$ 
is not sufficiently rich or optimally constructed to ensure small $\kappa_k$ values. 
For a robust subspace approximation, the basis $\mathbf{Q}$ must adequately capture the eigenmodes of interest, 
thereby keeping $\kappa_k$ small and helping to satisfy this condition.}

\begin{theorem}
\label{thm:eigvec-error}
    For $1 \leq j \leq \ell$, for any $\mathbf{x} \in \mathcal{E}_j$, we have 
    \begin{equation}
        \left\| \left(\mathbf{I} - \widetilde{\mathbf{P}}_{\mathbf{M}}^{\mathcal{S}_j} \right) \mathbf{x} \right\|_{\mathbf{M}} 
        \leq \left(1 + \tau_j \right) \| (\mathbf{I} - \mathbf{P}_{\mathbf{A}}) \mathbf{x} \|_{\mathbf{M}}.
    \end{equation}
    \begin{proof}
        Fix $1 \leq j \leq \ell$. For any $m \in \mathcal{S}_j$ and $1 \leq k \leq r$, by the properties of generalized eigenvector $\boldsymbol{\phi}_m$, we have 
        \begin{equation}
            \widetilde{\boldsymbol{\phi}}_k^\top \mathbf{A} \boldsymbol{\phi}_m
            = \lambda_m \widetilde{\boldsymbol{\phi}}_k^\top \mathbf{M} \boldsymbol{\phi}_m
            = \rev{\nu}_j \widetilde{\boldsymbol{\phi}}_k^\top \mathbf{M} \boldsymbol{\phi}_m.
        \end{equation}
        On the other hand, by the equality \eqref{eq:A-galerkin-orthogonality} and the properties of generalized eigenvector $\widehat{\boldsymbol{\phi}}_k$,
        \begin{equation}
            \widetilde{\boldsymbol{\phi}}_k^\top \mathbf{A} \boldsymbol{\phi}_m
            = \widetilde{\boldsymbol{\phi}}_k^\top \mathbf{A} \mathbf{P}_{\mathbf{A}}\boldsymbol{\phi}_m \\
            = \widetilde{\lambda}_k \widetilde{\boldsymbol{\phi}}_k^\top \mathbf{M} \mathbf{P}_{\mathbf{A}}\boldsymbol{\phi}_m.
        \end{equation}
        Let $\mathbf{x} \in \mathcal{E}_j$. Then, for any $k \in \mathcal{S}_j^{\text{c}}$, we have 
        \begin{equation}
            \widetilde{\boldsymbol{\phi}}_k^\top \mathbf{M} \mathbf{P}_{\mathbf{A}} \mathbf{x}
            = \dfrac{\nu_j}{\widetilde{\lambda}_k - \nu_j} \widetilde{\boldsymbol{\phi}}_k^\top \mathbf{M} (\mathbf{I} - \mathbf{P}_\mathbf{A}) \mathbf{x}.
        \end{equation}
        By the properties of $\mathbf{M}$-orthonormal basis, we infer that 
        \begin{equation}
            \left\| \left(\mathbf{I} - \widetilde{\mathbf{P}}_{\mathbf{M}}^{\mathcal{S}_j} \right) \mathbf{P}_{\mathbf{A}} \mathbf{x} \right\|_{\mathbf{M}}^2
            \leq \tau_j^2 \| (\mathbf{I} - \mathbf{P}_{\mathbf{A}})  \mathbf{x} \|_{\mathbf{M}}^2.
        \end{equation}
        By a triangle inequality for the $\mathbf{M}$-induced norm, we get  
        \begin{equation}
            \left\| \left(\mathbf{I} - \widetilde{\mathbf{P}}_{\mathbf{M}}^{\mathcal{S}_j} \mathbf{P}_{\mathbf{A}} \right) \mathbf{x} \right\|_{\mathbf{M}} 
        \leq (1 + \tau_j) \left\| \left(\mathbf{I} - \mathbf{P}_{\mathbf{A}}\right) \mathbf{x} \right\|_{\mathbf{M}}.  
        \end{equation}
        Finally, we arrive at the desired result by applying Lemma~\ref{lemma:M-projection-estimate} on $\mathcal{S}_j$.
    \end{proof}
\end{theorem}

%% file: numerical.tex
\section{Numerical examples}
\label{sec:numerical}

In this section, we consider the finite element discretization of 
a general spectral problem on elliptic differential operators 
in an open and bounded domain $\Omega \subset \mathbb{R}^d$: 
\begin{equation}
\begin{split}
-\text{div}(\sigma(x; \rev{\boldsymbol{\mu}}) \nabla \phi(x)) + \rho(x; \rev{\boldsymbol{\mu}}) \phi(x) & = \lambda \phi(x) \text { for all } x \in \Omega, \\
\alpha(x; \rev{\boldsymbol{\mu}}) \phi(x) + \beta(x; \rev{\boldsymbol{\mu}}) \mathbf{n} \cdot \nabla \phi(x) & = 0 \text{ for all } x \in \partial \Omega. 
\end{split}
\label{eq:elliptic_evp}
\end{equation}
Here $\rev{\boldsymbol{\mu}} \in \textsf{D}$ is a physical parameter for parametrizing the coefficient fields, 
including the conductivity $\sigma: \Omega \to \mathbb{R}_{++}$, 
the potential $\rho: \Omega \to \mathbb{R}_{+}$, 
and the boundary coefficients $\alpha, \beta: \partial \Omega \to \mathbb{R}_{+}$. 
It is assumed that $\alpha(x) > 0$ or $\beta(x) > 0$ for $x \in \partial \Omega$, and we denote 
$\Gamma_D = \{ x \in \partial \Omega: \beta(x) = 0 \}$.

Let $\mathcal{V}_h$ be a conforming finite element space for the boundary value problem, i.e. 
\begin{equation}
\mathcal{V}_h \subset \{ v \in H^1(\Omega): v \vert_{\Gamma_D} = 0 \}, 
\end{equation}
and $\{ v_i \}_{i=1}^{n}$ be a basis for $\mathcal{V}_h$. 
The finite element discretization of \eqref{eq:elliptic_evp} reads: find $(\lambda, \phi_h) \in \mathbb{R} \times \mathcal{V}_h$ such that 
\begin{equation}
\int_\Omega \left(a \nabla \phi_h \cdot \nabla v + q \phi_h v \right) \, \text{d}x
+ \int_{\partial\Omega \setminus \Gamma_D} \dfrac{\alpha}{\beta} \phi_h v \, \text{d}s 
= \lambda \int_\Omega \phi_h v  \, \text{d}x, \text{ for all } v \in \mathcal{V}_h, 
\end{equation}
which translates to the algebraic problem \eqref{eq:fom-evp} with $\boldsymbol{\phi} \in \mathbb{R}^n$ 
being the coefficients of $\phi_h$ with respect to $\{ v_i \}_{i=1}^{n}$, and 
\begin{equation}
\label{eq:fem}
\begin{split}
\mathbf{A}_{ij} &
= \int_\Omega \left(a\nabla v_i \cdot \nabla v_j + q v_i v_j \right) \, \text{d}x
+ \int_{\partial\Omega \setminus \Gamma_D} \dfrac{\alpha}{\beta} v_i v_j \, \text{d}s, \\
\mathbf{M}_{ij} &
=  \int_\Omega v_i v_j  \, \text{d}x. 
\end{split}
\end{equation}
We refer to the resultant eigenvalue problem \eqref{eq:fom-evp} as full order model (FOM).

For each example, the basis matrix $\mathbf{Q}$ will be constructed from corresponding sample snapshot data on 
a set of sample parameters $\textsf{D}_{\text{train}} = \{ \rev{\boldsymbol{\mu}}_\ell \}_{\ell = 1}^{s} \subset \textsf{D}$. 
More precisely, the eigenvectors $\{ \boldsymbol{\phi}_k(\rev{\boldsymbol{\mu}}_\ell) \}_{k=1}^{p}$ are computed for $1 \leq \ell \leq s$, 
and assembled to form a snapshot matrix of size $n \times r$, where $r = ps$. 
The basis matrix $\mathbf{Q}$ can be obtained from QR factorization of the snapshot matrix, 
and is used to formulate the projected system \eqref{eq:rom-evp} at generic 
\rev{testing parameters $\boldsymbol{\mu} \in \textsf{D}_{\text{test}} \subset \textsf{D}$}. 
We refer to the projected eigenvalue problem \eqref{eq:rom-evp} as reduced order model (ROM).

\rev{
In the first two examples, we present numerical results to verify the analysis in Section~\ref{sec:analysis}. 
At selected parameters $\boldsymbol{\mu} \in \textsf{D}$, for $1 \leq k \leq r$, we will show 
the ROM approximation of the $k$-th eigenvalue $\widetilde{\lambda}_k$, 
with the bounds $\lambda_k$ and $\kappa_k^2 \lambda_k$ in Theorem~\ref{thm:eigval-error}, 
and the ROM approximation error of the $k$-th eigenvector $\varepsilon_k$ 
with the bound $(1+\tau_j) \delta_k$ in Theorem~\ref{thm:eigvec-error}, where $\nu_j = \lambda_k$ and
\begin{equation}
\begin{split}
\varepsilon_k & = \left\| \left(\mathbf{I} - \widetilde{\mathbf{P}}_{\mathbf{M}}^{\mathcal{S}_j} \right) \boldsymbol{\phi}_k \right\|_{\mathbf{M}}, \\
\delta_k & = \left\| \left(\mathbf{I} - \widetilde{\mathbf{P}}_{\mathbf{A}} \right) \boldsymbol{\phi}_k \right\|_{\mathbf{M}}.
\end{split}
\end{equation}
We remark that the approximations for $p+1 \leq k \leq r$ exceed our need of approximating the first $p$ eigenpairs, 
but the results are presented for completeness. 
The verification of the eigenvalue upper bound relies on the evaluation of $\kappa_k$, 
for which the precise procedure is outlined in Appendix~\ref{sec:compute_kappa}.
While our theoretical analysis in Theorem~\ref{thm:eigvec-error} focuses on the error of arbitrary elements in eigenspaces, 
we present results for individual eigenvectors approximation for computational convenience and more direct insight into the approximation of specific eigenmodes. 
We remark the accurate approximations of eigenvectors do not directly imply the distance between the projectors onto the FOM eigenspace and ROM eigenspace is small. However, the metric still offers valuable information on the approximation quality of the computed eigenmodes.
In all the examples, we evaluate the accuracy of the ROM approximation 
at generic testing parameters $\rev{\boldsymbol{\mu}} \in \textsf{D}_{\text{test}}$ by 
the error in the $k$-th eigenvalue $\widetilde{\lambda}_k - \lambda_k$ 
the error in the $k$-th eigenvector $\varepsilon_k$. 
}

In the first two examples, the simulations are performed using MATLAB on Apple M1 Pro. 
In the last three examples, the simulations are performed using the finite element methods library MFEM \cite{mfem} 
on Dane in Livermore Computing Center \footnote{High performance computing at LLNL, https://hpc.llnl.gov/hardware/compute-platforms/dane},
on Intel Sapphire Rapids CPUs with 256 GB memory, and peak TFLOPS of 10,723.0.
A visualization tool, VisIt \cite{ahern2013visit}, is used to visualize the coefficient and solution fields. 

\subsection{1D parametric boundary-value Laplacian}

In this example, the computational domain is taken as the unit interval 
$\Omega = (0, 1) \subset \mathbb{R}$, and the parameter domain is 
$\textsf{D} = \rev{[0,10]}$ 
\rev{which reduces the parameter $\boldsymbol{\mu}$ to a scalar $\mu$}.
For $\mu \in \textsf{D}$, 
the coefficient fields are given by 
\begin{equation}
\begin{split}
\sigma(x; \mu) & = 1 \quad \text{for } x \in \Omega, \\
\rho(x; \mu) & = 0 \quad \text{for } x \in \Omega, 
\end{split}
\end{equation}
i.e., the Laplacian operator is under consideration, 
and the boundary conditions are taken as 
\begin{equation}
\begin{split}
\alpha(0; \mu) = 1, \quad & \beta(0; \mu) = 0, \\
\alpha(1; \mu) = \mu, \quad & \beta(1; \mu) = 1,
\end{split}
\end{equation}
i.e., the Dirichlet boundary condition is applied at $x = 0$, 
and the Neumann and Robin boundary conditions are applied at $x = 1$ 
for $\mu = 0$ and for $\mu > 0$, respectively. 
It is a classical result that the eigenfunctions are of the form $\phi(x) = \sin(\lambda x)$, 
where the spectral values $\lambda$ are solutions of 
$\sqrt{\lambda} + \mu \tan(\sqrt{\lambda}) = 0$.

\begin{figure}[ht!]
\centering
\includegraphics[width=0.48\linewidth]{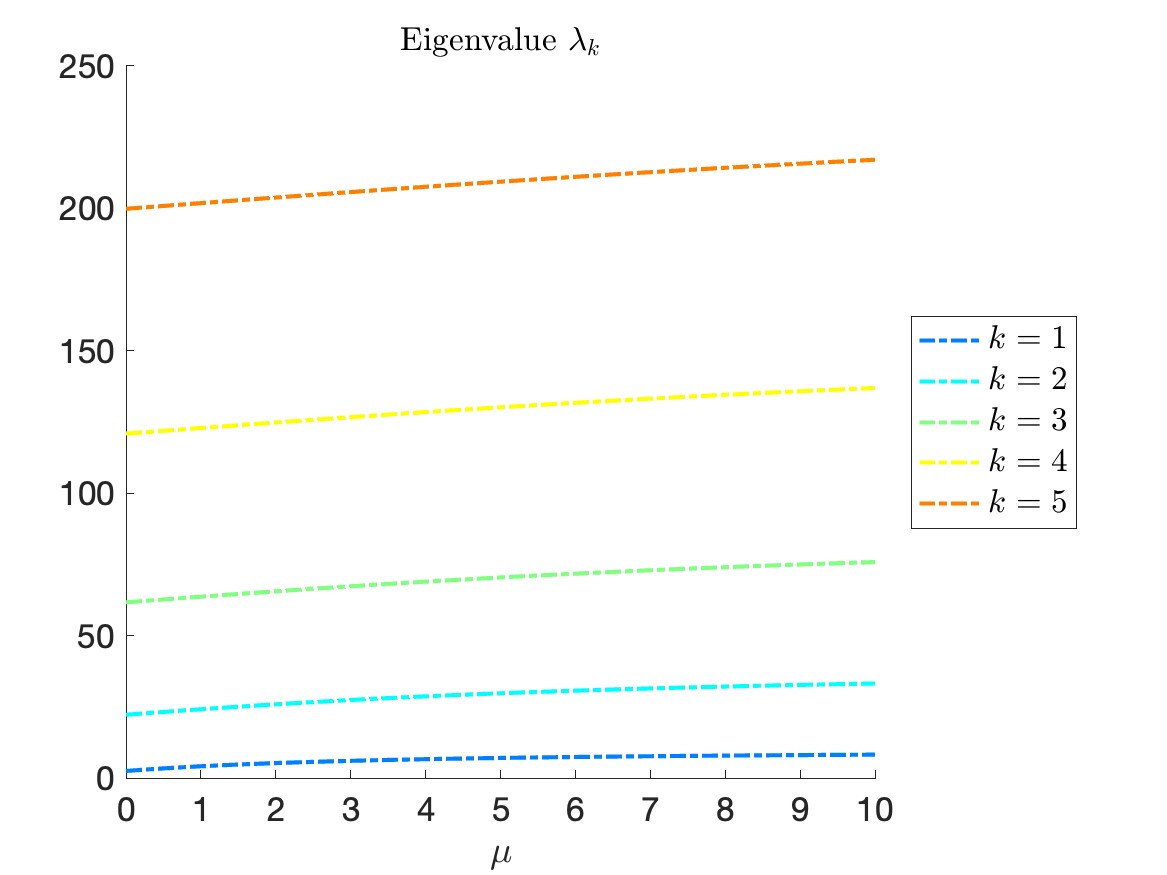}
\caption{First 5 eigenvalues at different parameters \rev{$\mu \in \textsf{D}$}  in the example of parametric boundary value Laplacian.}
\label{fig:laplace_robin_eigenvalue}
\end{figure}

The domain $\Omega$ is divided into a uniform mesh of size $h = 1/2000$, 
and $P^1$ Lagrange finite elements is used in the finite element discretization, 
which results in a system size of $n = 1999$. 
We are interested in the first $p = 5$ eigenpairs.
The eigenvalue problem \eqref{eq:fom-evp} is solved by the built-in function \texttt{eigs} in MATLAB.
Figure~\ref{fig:laplace_robin_eigenvalue} shows the first 5 eigenvalues $\{\lambda_k(\mu)\}_{k=1}^5$ 
of \eqref{eq:fom-evp} at different parameters $\mu \in [0, 10]$, 
which is verified to be the 5 smallest positive approximated solutions of $\sqrt{\lambda} + \mu \tan(\sqrt{\lambda}) = 0$. 

\begin{figure}[ht!]
\centering
\includegraphics[width=0.48\linewidth]{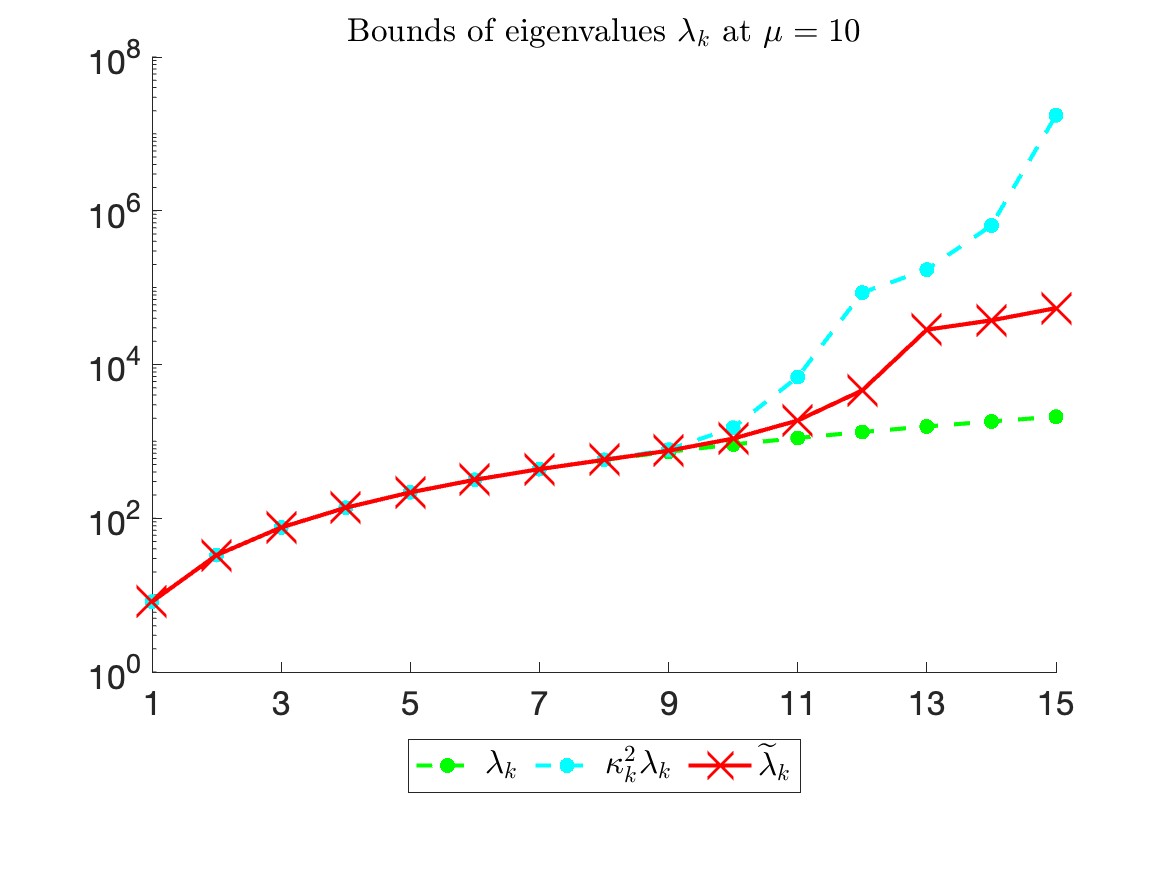}
\includegraphics[width=0.48\linewidth]{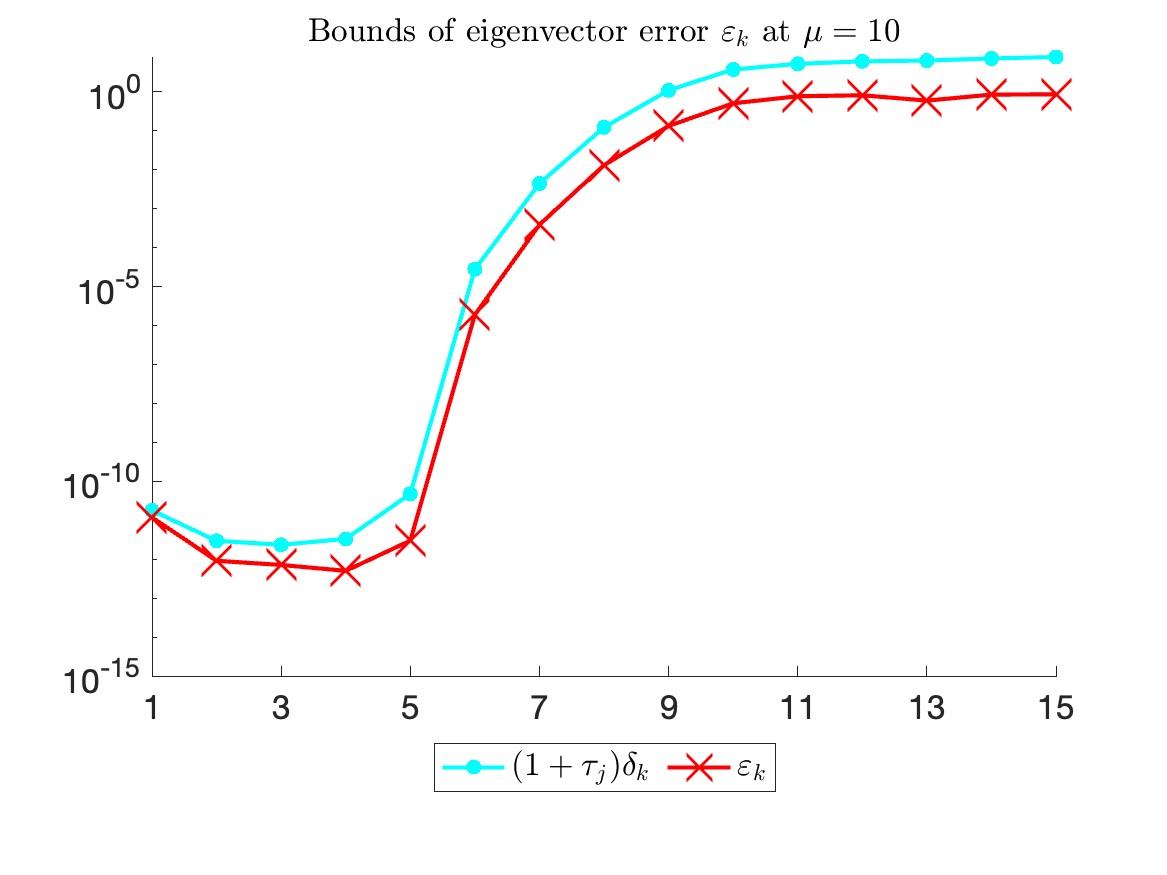}
\caption{Bounds of all 15 eigenvalues (left) and eigenvectors (right) at the parameter  
$\mu = 10$ in the example of parametric boundary value Laplacian.}
\label{fig:laplace_robin_error_bound}
\end{figure}

\begin{figure}[ht!]
\centering
\includegraphics[width=0.48\linewidth]{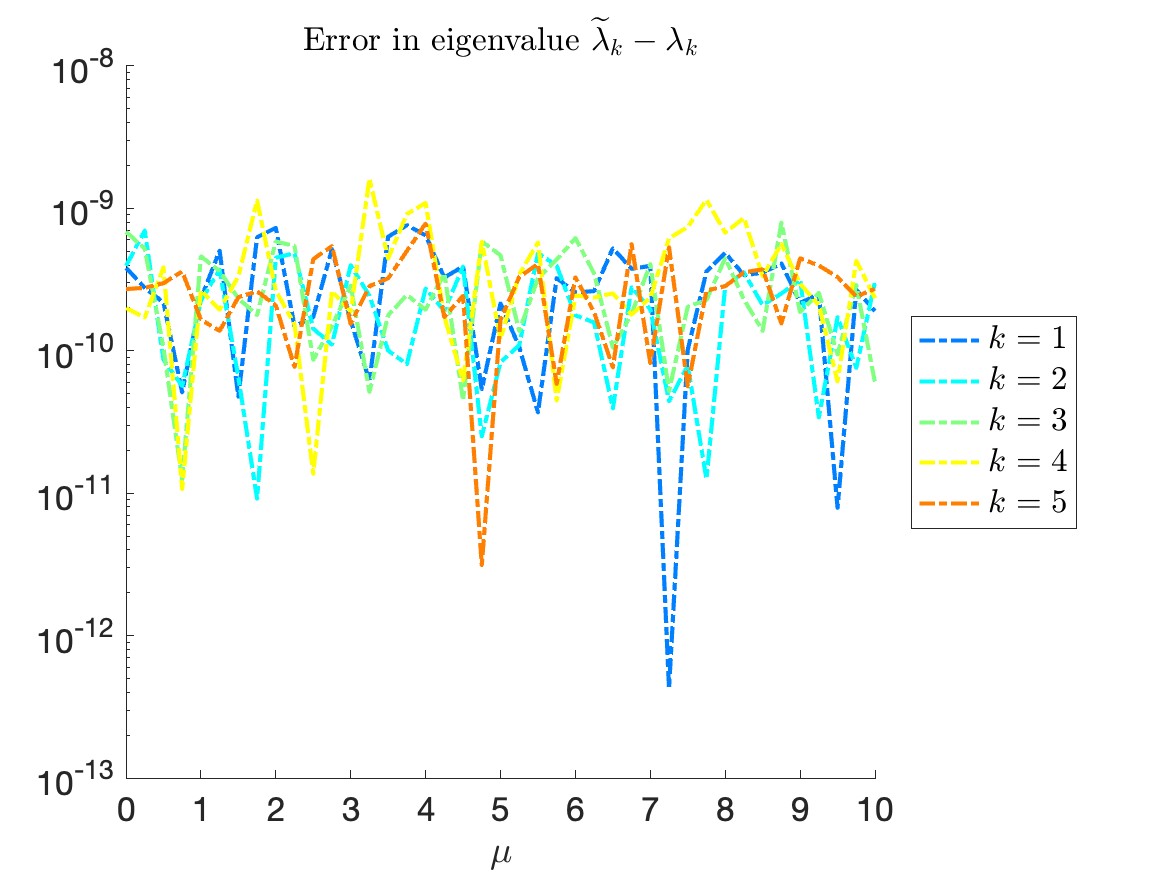}
\includegraphics[width=0.48\linewidth]{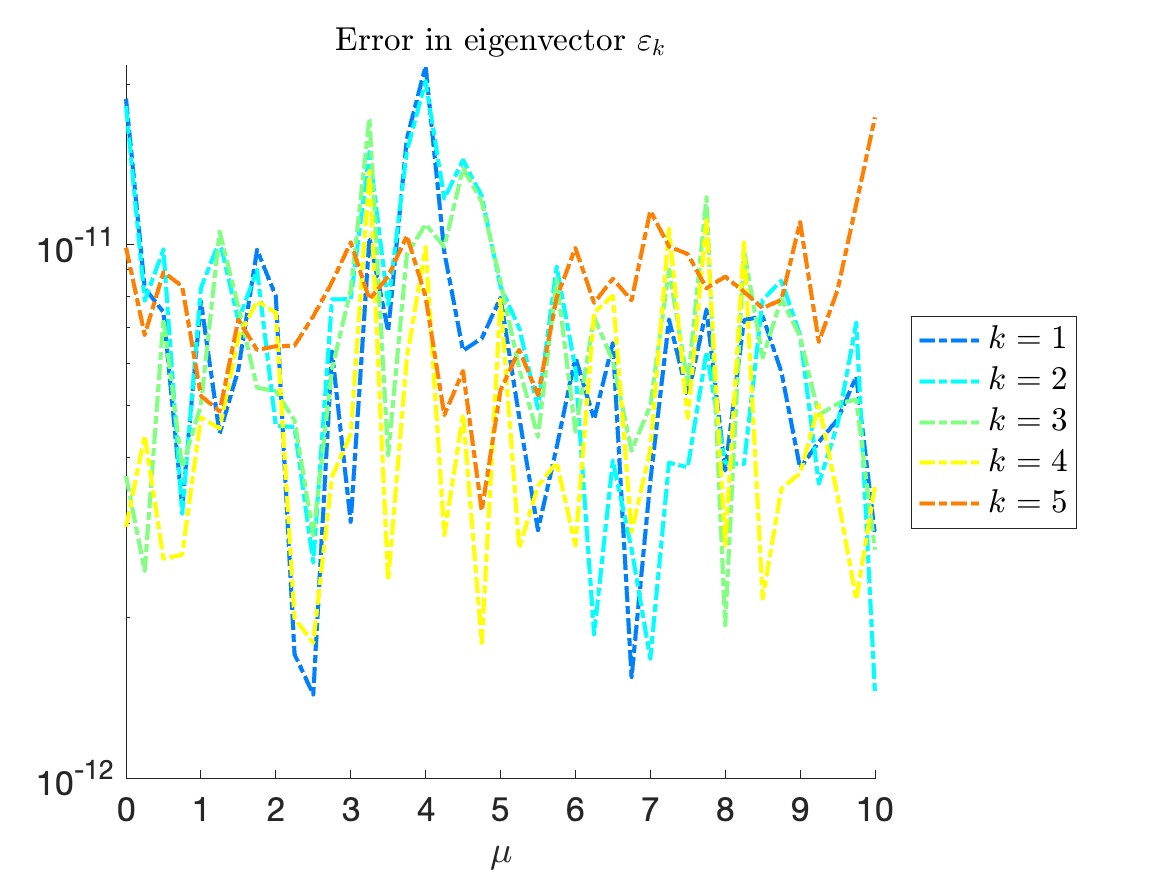}
\caption{Error in first 5 eigenvalues (left) and eigenvectors (right) at different testing parameters \rev{$\mu \in \textsf{D}_{\text{test}}$} 
in the example of parametric boundary value Laplacian.}
\label{fig:laplace_robin_error}
\end{figure}

The eigenvectors $\{ \boldsymbol{\phi}_k(\mu) \}_{k=1}^5$ at the sample parameters $\mu \in \textsf{D}_{\text{train}} = \{ 1, 5, 9 \}$ are used as snapshots 
to construct the basis matrix $\mathbf{Q}$ of column size $r = 15$.  
The projected system \eqref{eq:rom-evp} is solved with the built-in function \texttt{eig} in MATLAB.
Since all eigenvalues are \rev{simple}, direct comparison can be made to the FOM and ROM eigenvectors of the same index.
\rev{
To verify the theoretical results in Section~\ref{sec:analysis}, 
Figure~\ref{fig:laplace_robin_error_bound} illustrates the behavior of the ROM approximated $k$-th eigenvalue $\widetilde{\lambda}_k$ and 
the ROM approximation error of the $k$-th eigenvector $\varepsilon_k$, for $1 \leq k \leq r = 15$ at the parameter $\mu = 10$. 
The figure explicitly demonstrates that these quantities are bounded by 
$\lambda_k \leq \widetilde{\lambda}_k \leq \kappa_k^2 \lambda_k$ and 
$\varepsilon_k \leq (1+\tau_j) \delta_k$ as shown in 
Theorem~\ref{thm:eigval-error} and Theorem~\ref{thm:eigvec-error}, respectively. 
This verification confirms the validity of our theoretical error estimates.
We note that the observed deterioration in accuracy for $p+1 \leq k \leq r$
is attributed to the construction of the basis matrix $\mathbf{Q}$
from only the first $p=5$ eigenvectors at $\mu \in \textsf{D}_{\text{train}}$.
}
Figure~\ref{fig:laplace_robin_error} shows the error in the first 5 eigenvalues and eigenvectors at different parameters 
\rev{$\mu \in \textsf{D}_{\text{test}} = \{ 0.25t : t \in [0, 40] \cap \mathbb{Z} \}$}. 
It can be seen that both approximations are extremely accurate uniformly in $\mu \in [0,10]$. 
The error in the eigenvalues remain below $10^{-8}$, and 
the error in the eigenvectors remain below $10^{-10}$. 

\subsection{1D parametric harmonic oscillator}

\begin{figure}[ht!]
\centering
\includegraphics[width=0.48\linewidth]{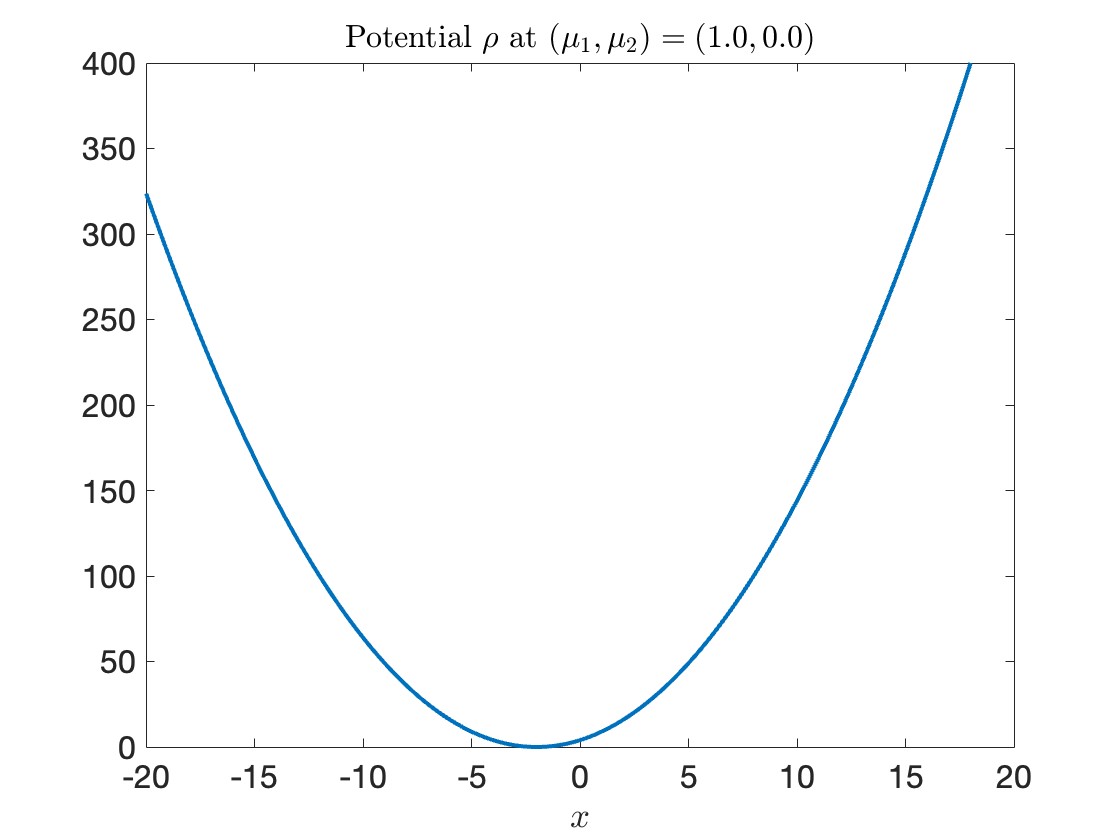}
\includegraphics[width=0.48\linewidth]{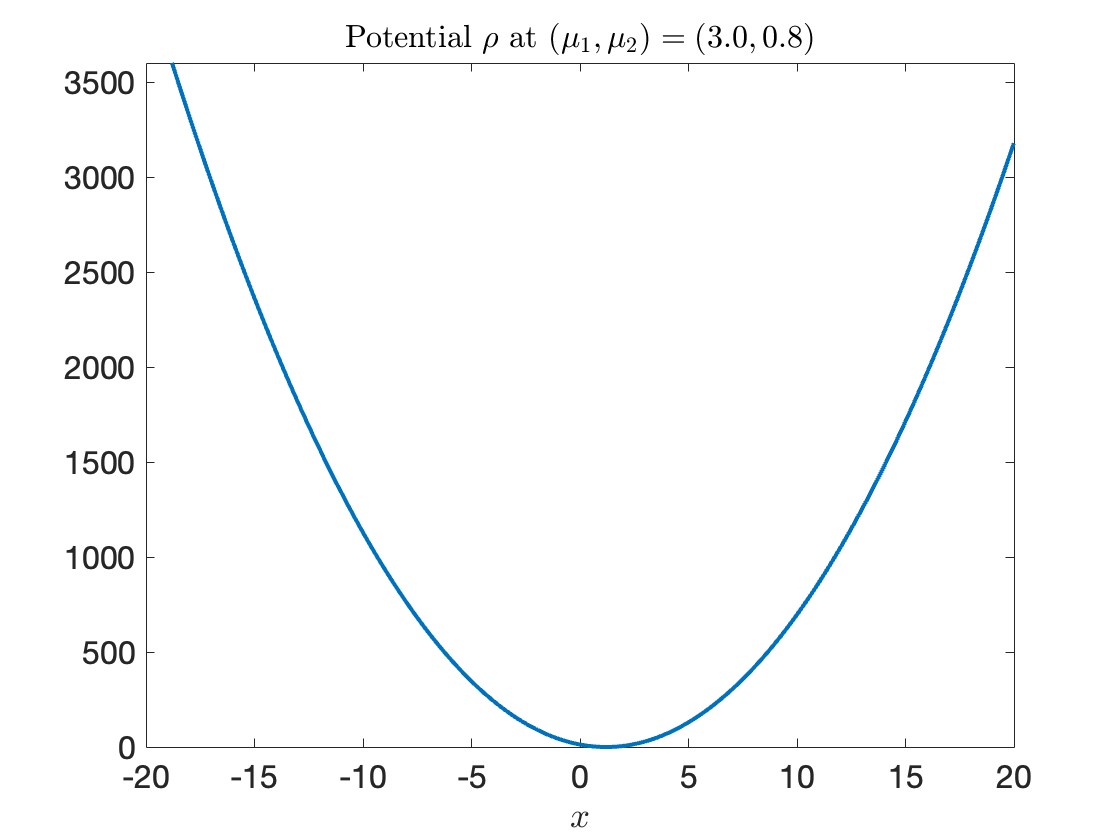}
\caption{Potential $\rho$ at the parameters $\boldsymbol{\mu} = (1, 0)$ (left) and $\boldsymbol{\mu} = (3, 0.8)$ (right) in the example of parametric harmonic oscillator.}
\label{fig:schrodinger_harmonic_potential}
\end{figure}

In this example, the computational domain is taken as the interval 
$\Omega = (-20, 20) \subset \mathbb{R}$, and the parameter domain is 
\rev{$\textsf{D} = [1,3] \times [0,1]$}.
For \rev{$\boldsymbol{\mu} = (\mu_1, \mu_2) \in \textsf{D}$}, 
the coefficient fields are given by 
\begin{equation}
\begin{aligned}
\sigma(x; \boldsymbol{\mu}) & = 1                   && \text{for } x \in \Omega, \\
\rho(x; \boldsymbol{\mu}) & = \mu_1^2(x + 2 - 4\mu_2)^2   && \text{for } x \in \Omega,
\end{aligned}
\end{equation}
which corresponds to the Schr\"{o}dinger operator with parametric harmonic oscillator potential,
\rev{with the parameter $\mu_1$ controls the magnitude of the potential $\rho$ 
and $\mu_2$ controls the equilibrium position.}
The potential function $\rho$ with the parameters \rev{$\boldsymbol{\mu} = (1, 0)$} and \rev{$\boldsymbol{\mu} = (3, 0.8)$} 
are depicted in Figure~\ref{fig:schrodinger_harmonic_potential}. 
The Neumann boundary condition is prescribed, i.e., 
\begin{equation}
\alpha(x; \boldsymbol{\mu}) = 0, \quad \beta(x; \boldsymbol{\mu}) = 1, \quad \text{for } x \in \partial\Omega = \{-20, 20\}. 
\end{equation}
This boundary condition represents a truncation of the unbounded problem on the real line, 
with the equilibrium position in the harmonic oscillator potential shifting along $[-2,2]$ as \rev{$\mu_2$} varies in $[0,1]$. 
As a result, the eigenfunctions $\phi(x)$ also shifts as \rev{$\mu_2$} varies in $[0,1]$. 
On the other hand, the spectral values are parametric-independent. It can be easily shown that $k$-th analytic spectral values of the unbounded problem on the real line is \rev{$(2k-1) \mu_1$} \cite{zwiebach2022mastering}.

\begin{figure}[ht!]
\centering
\includegraphics[width=0.48\linewidth]{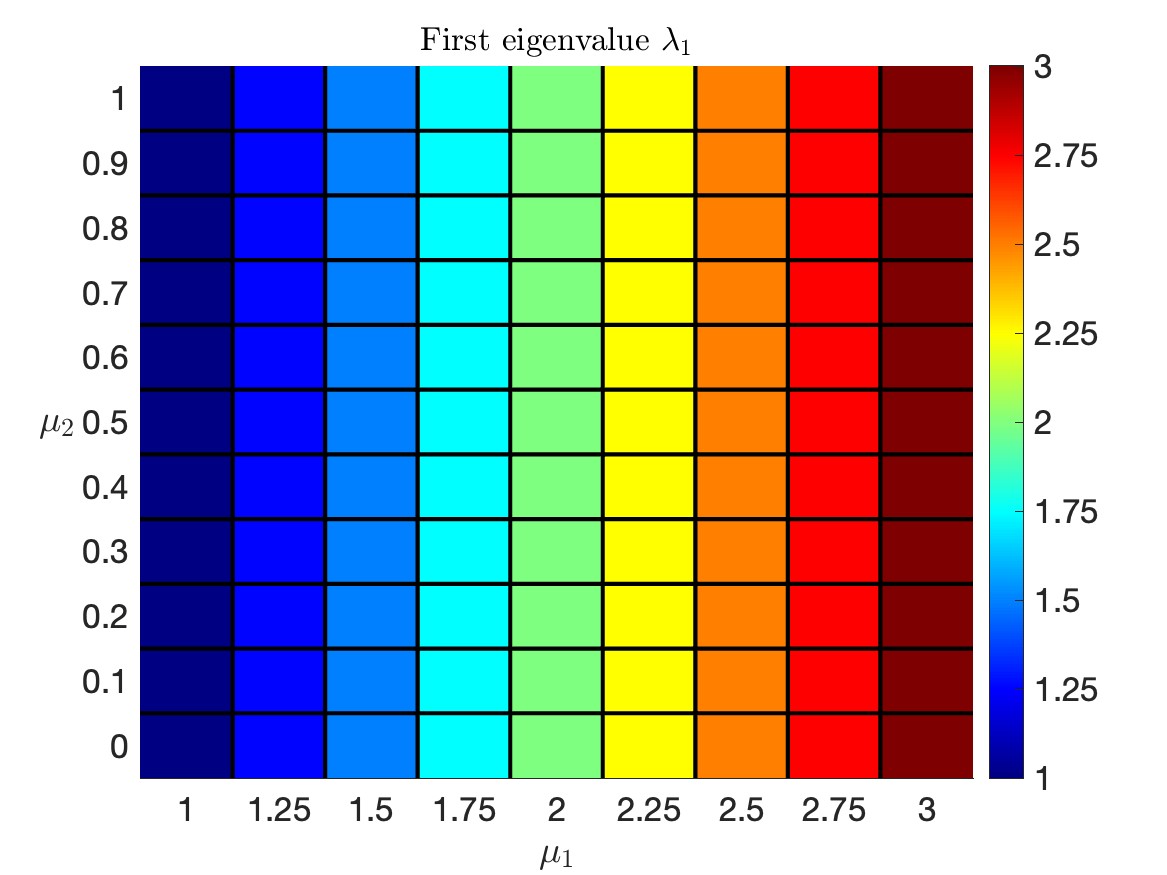}
\includegraphics[width=0.48\linewidth]{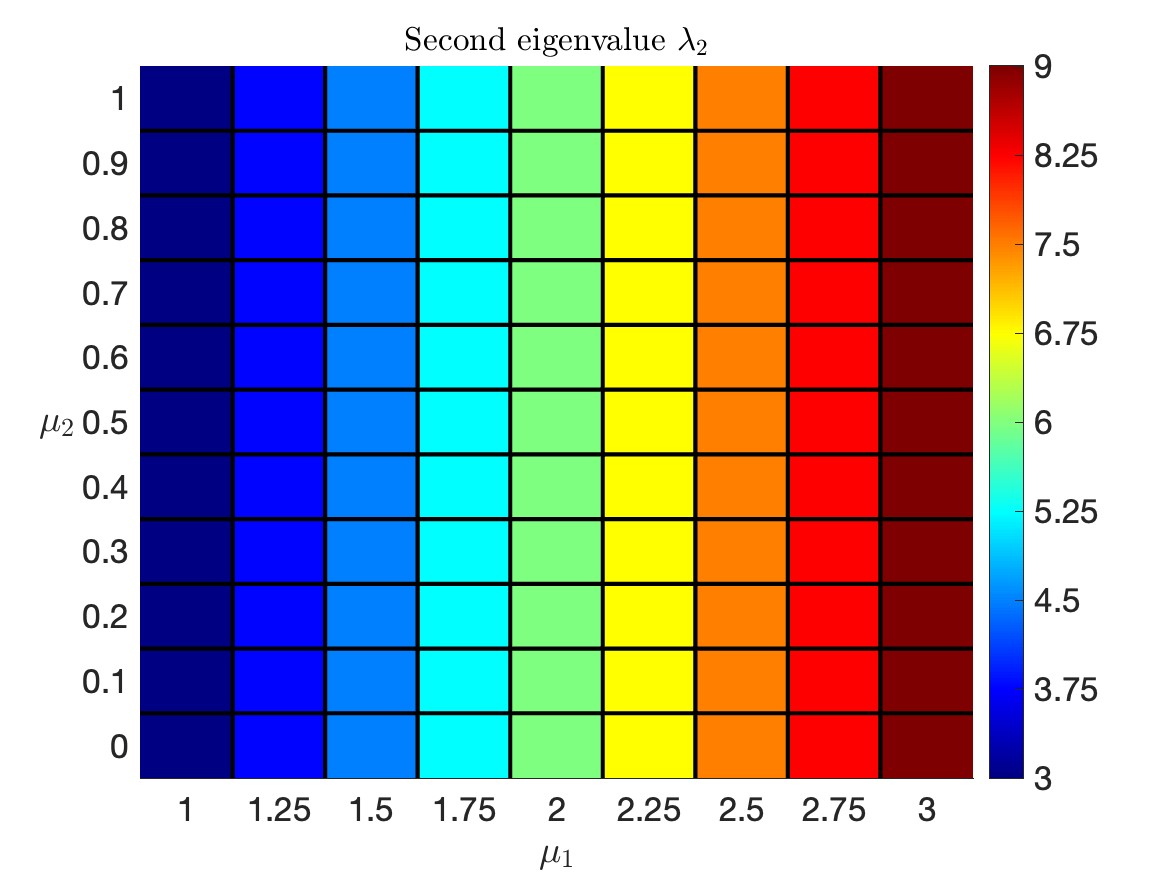} \\
\includegraphics[width=0.48\linewidth]{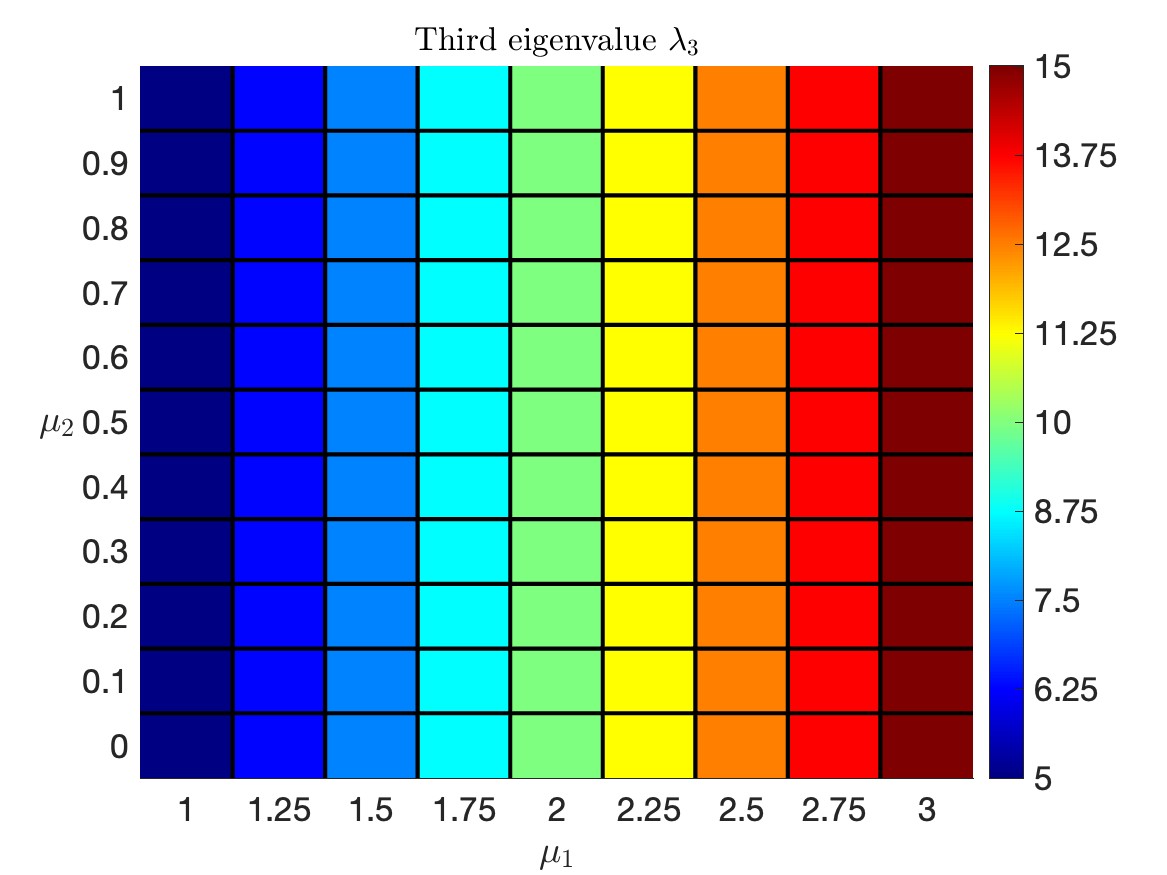}
\includegraphics[width=0.48\linewidth]{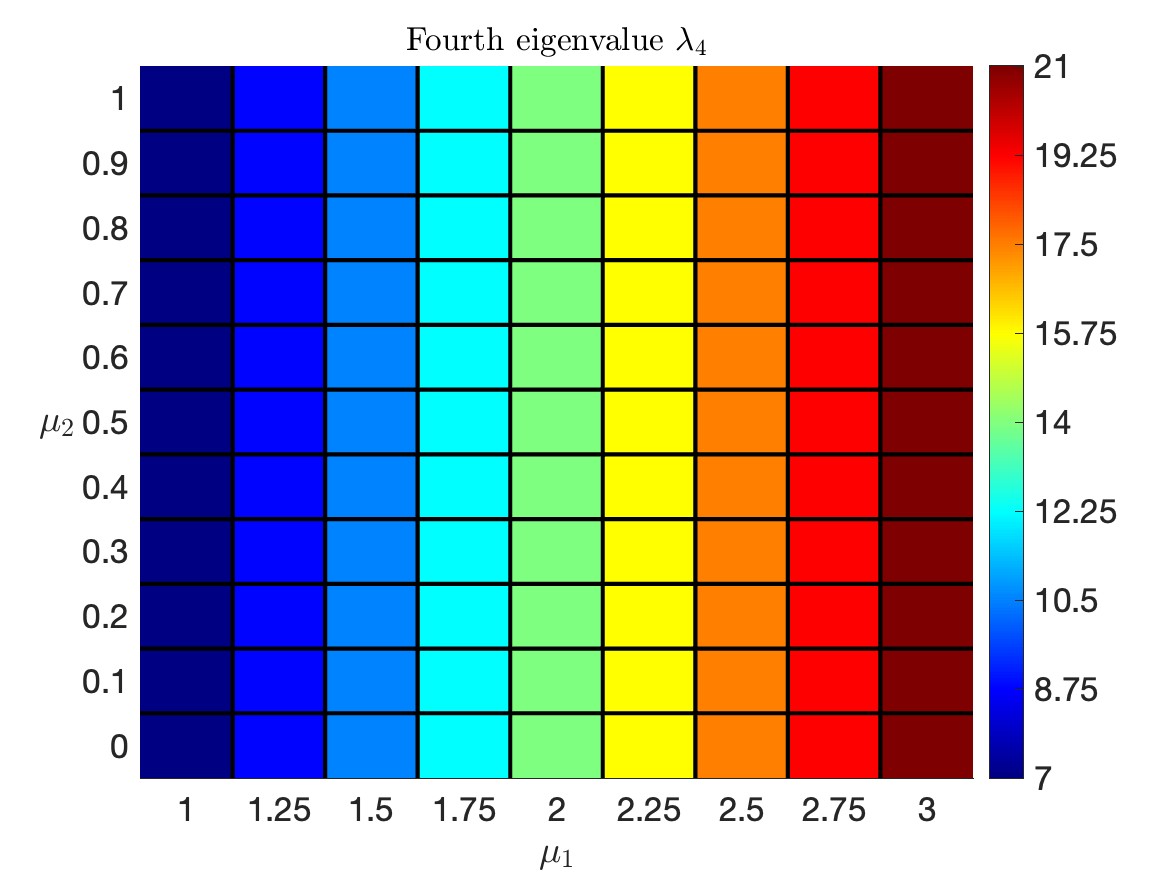}
\caption{First 4 eigenvalues at different parameters $\boldsymbol{\mu} \in \textsf{D}$ in the example of parametric harmonic oscillator.}
\label{fig:schrodinger_harmonic_eigenvalue_coarse}
\end{figure}

The domain $\Omega$ is divided into a uniform mesh of size $h = 1/50$, 
and $P^1$ Lagrange finite elements is used in the finite element discretization, 
which results in a system size of $n = 2001$
We are interested in the first \rev{$p = 4$} eigenpairs.
The eigenvalue problem \eqref{eq:fom-evp} is solved by the built-in function \texttt{eigs} in MATLAB.
Figure~\ref{fig:schrodinger_harmonic_eigenvalue_coarse} shows the first \rev{4 eigenvalues $\{\lambda_k(\boldsymbol{\mu})\}_{k=1}^4$}
of \eqref{eq:fom-evp} at different parameters $\boldsymbol{\mu} \in \textsf{D}$. 
As the analytic spectral values of the unbounded problem are parametric-independent and the computational domain is sufficiently large, 
the FOM eigenvalues accurately approximate the corresponding analytic spectral values 
uniformly for $\boldsymbol{\mu} \in \textsf{D}$. 

\begin{figure}[ht!]
\centering
\includegraphics[width=0.48\linewidth]{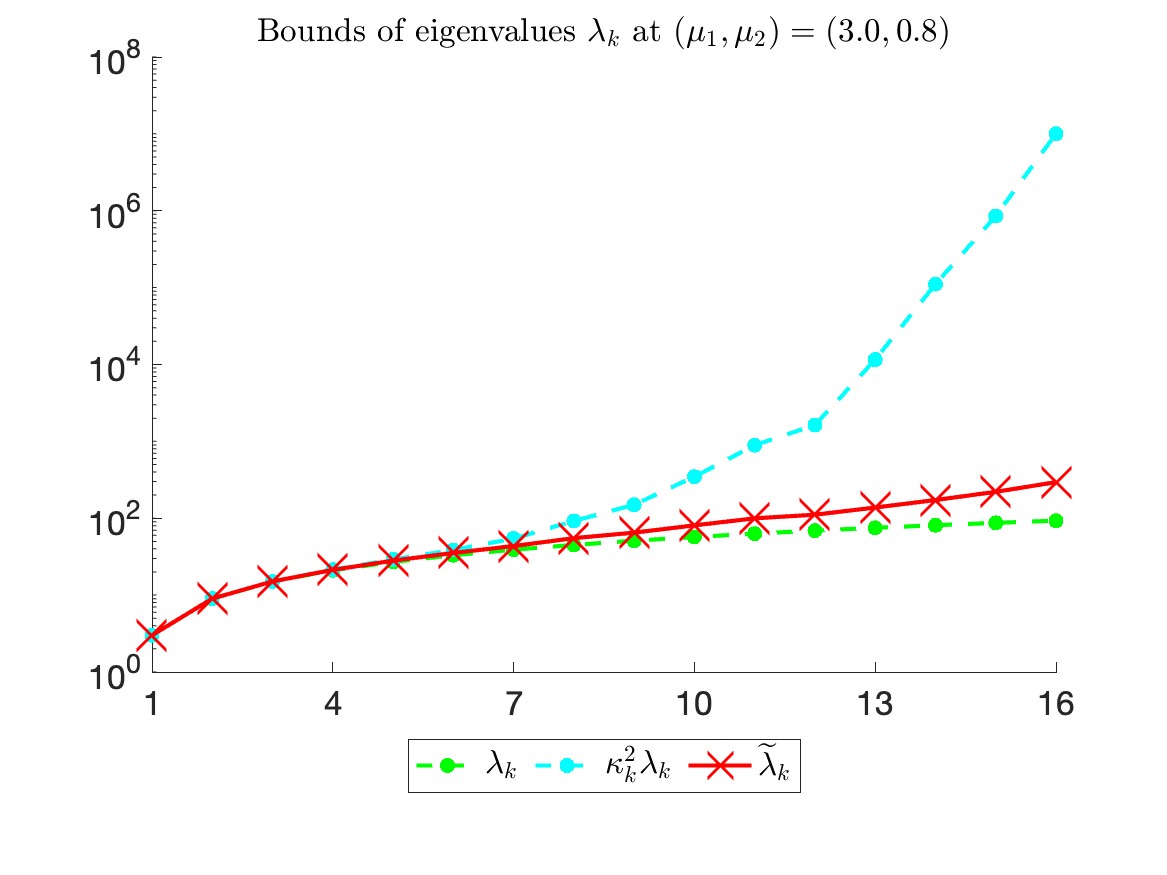}
\includegraphics[width=0.48\linewidth]{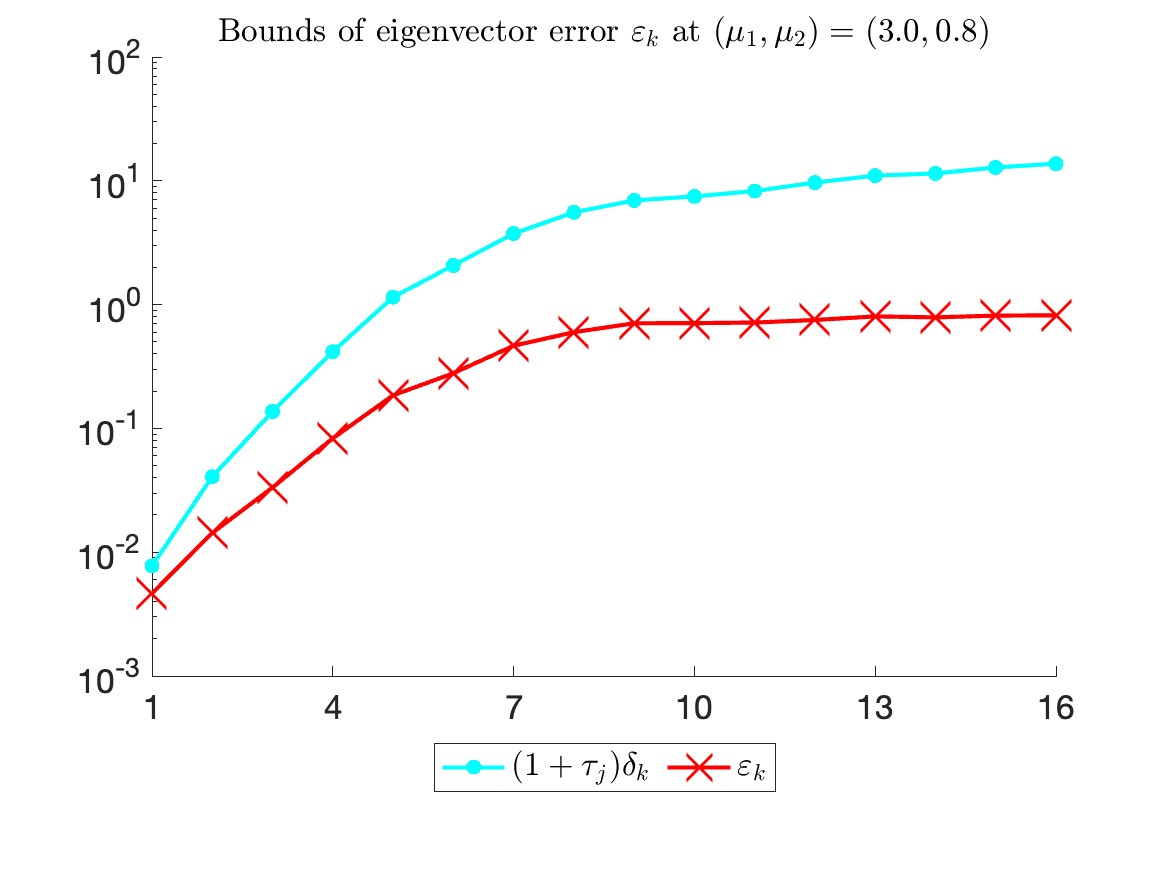}
\caption{Bounds of all 16 eigenvalues (left) and eigenvectors (right) at the parameter  
$\boldsymbol{\mu} = (3, 0.8)$ 
with coarse sampling of training parameters $\textsf{D}_{\text{train}} = \{ 1, 3 \} \times \{ 0, 1 \}$
in the example of parametric harmonic oscillator.}
\label{fig:schrodinger_harmonic_error_bound_coarse}
\end{figure}

\begin{figure}[ht!]
\centering
\includegraphics[width=0.48\linewidth]{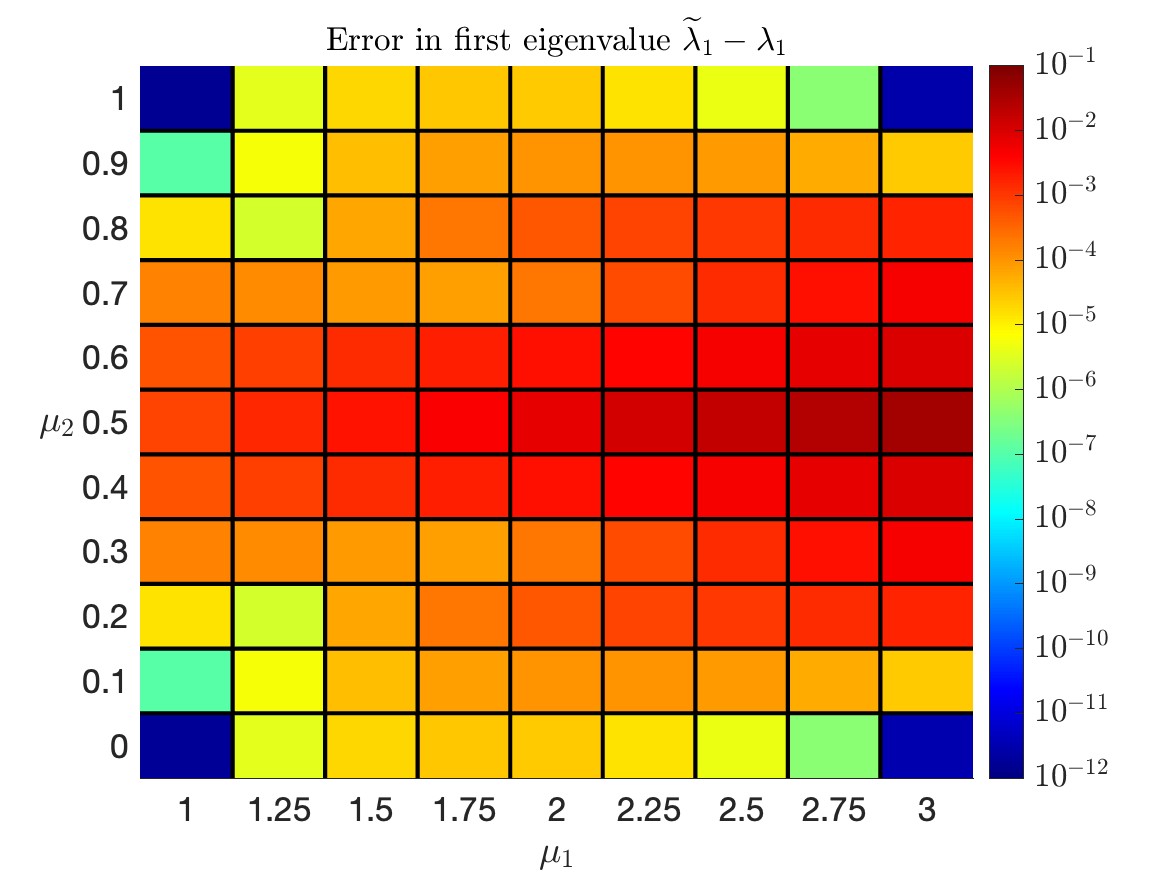}
\includegraphics[width=0.48\linewidth]{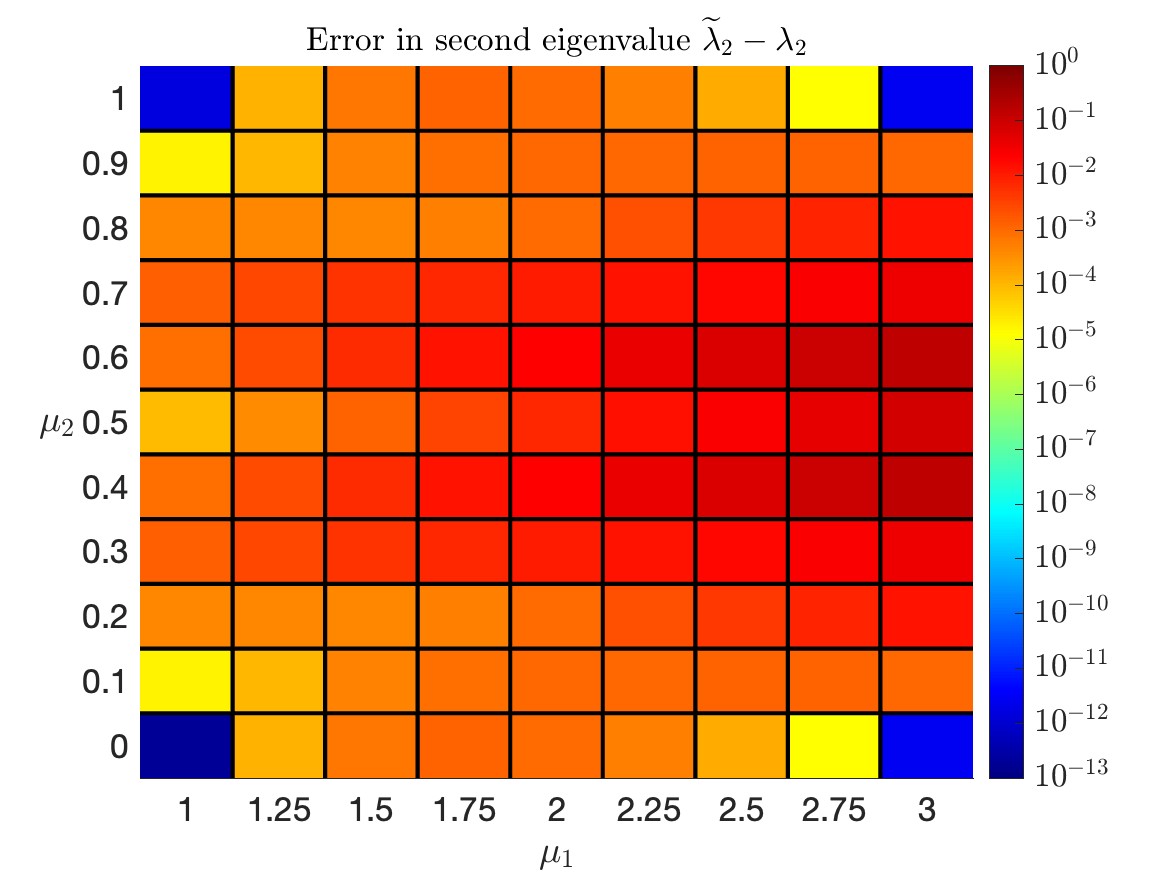}\\
\includegraphics[width=0.48\linewidth]{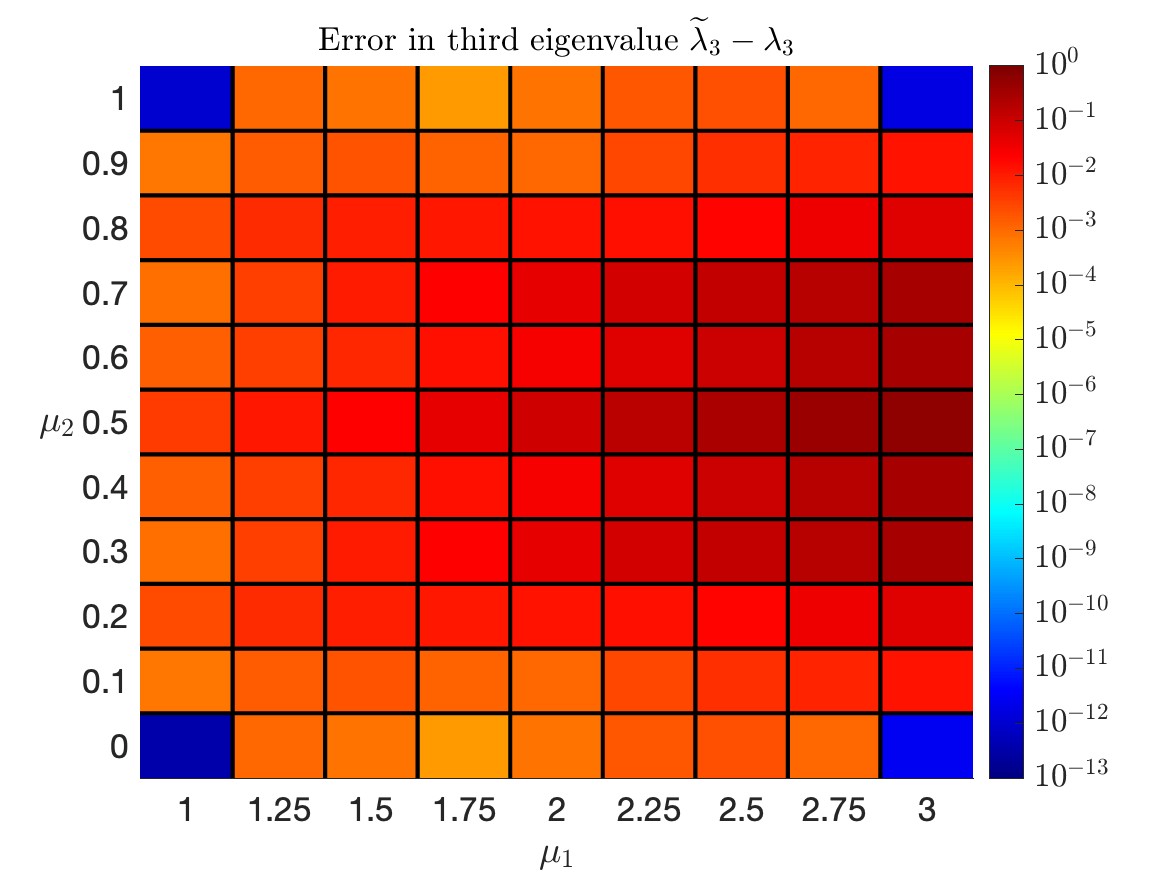}
\includegraphics[width=0.48\linewidth]{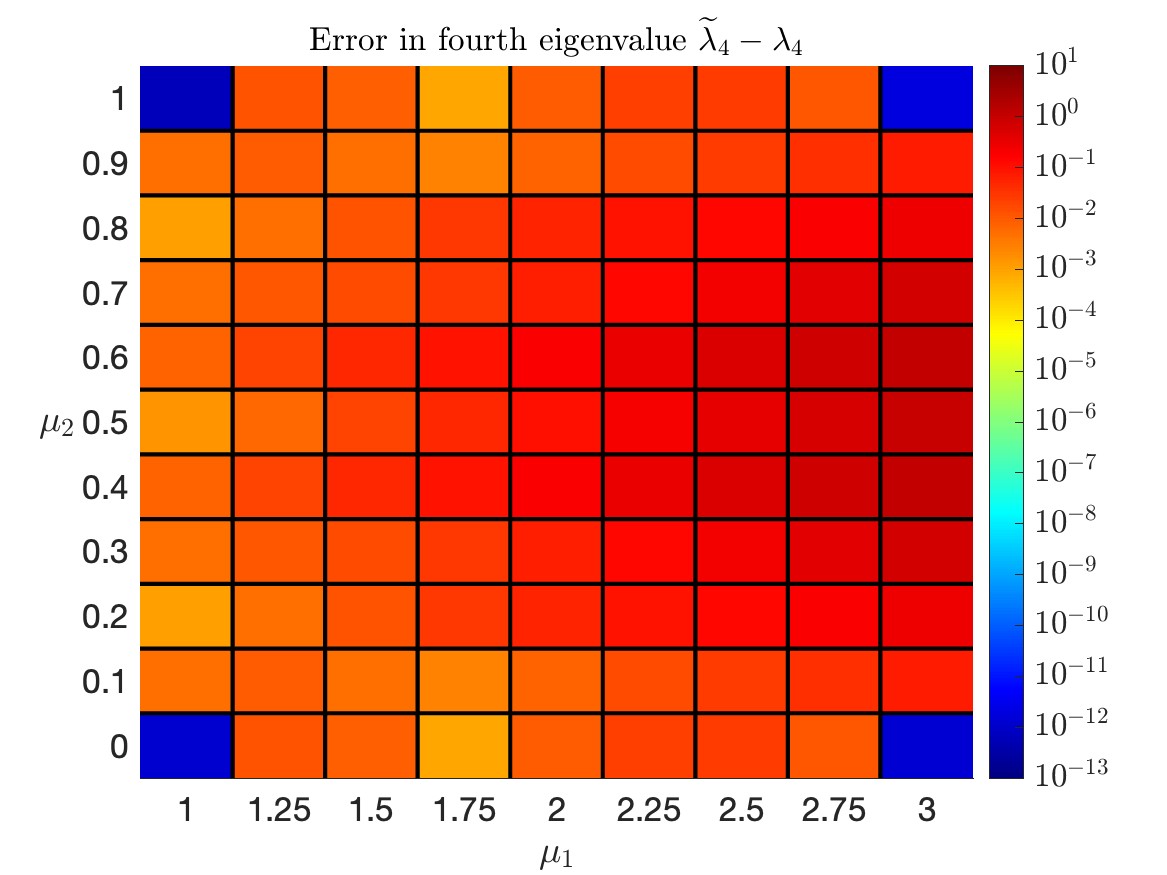}\\
\caption{Error in first 4 eigenvalues at different testing parameters 
\rev{$\boldsymbol{\mu} \in \textsf{D}_{\text{test}}$} 
with coarse sampling of training parameters $\textsf{D}_{\text{train}} = \{ 1, 3 \} \times \{ 0, 1 \}$
in the example of parametric harmonic oscillator.}
\label{fig:schrodinger_harmonic_eigenvalue_error_coarse}
\end{figure}

\begin{figure}[ht!]
\centering
\includegraphics[width=0.48\linewidth]{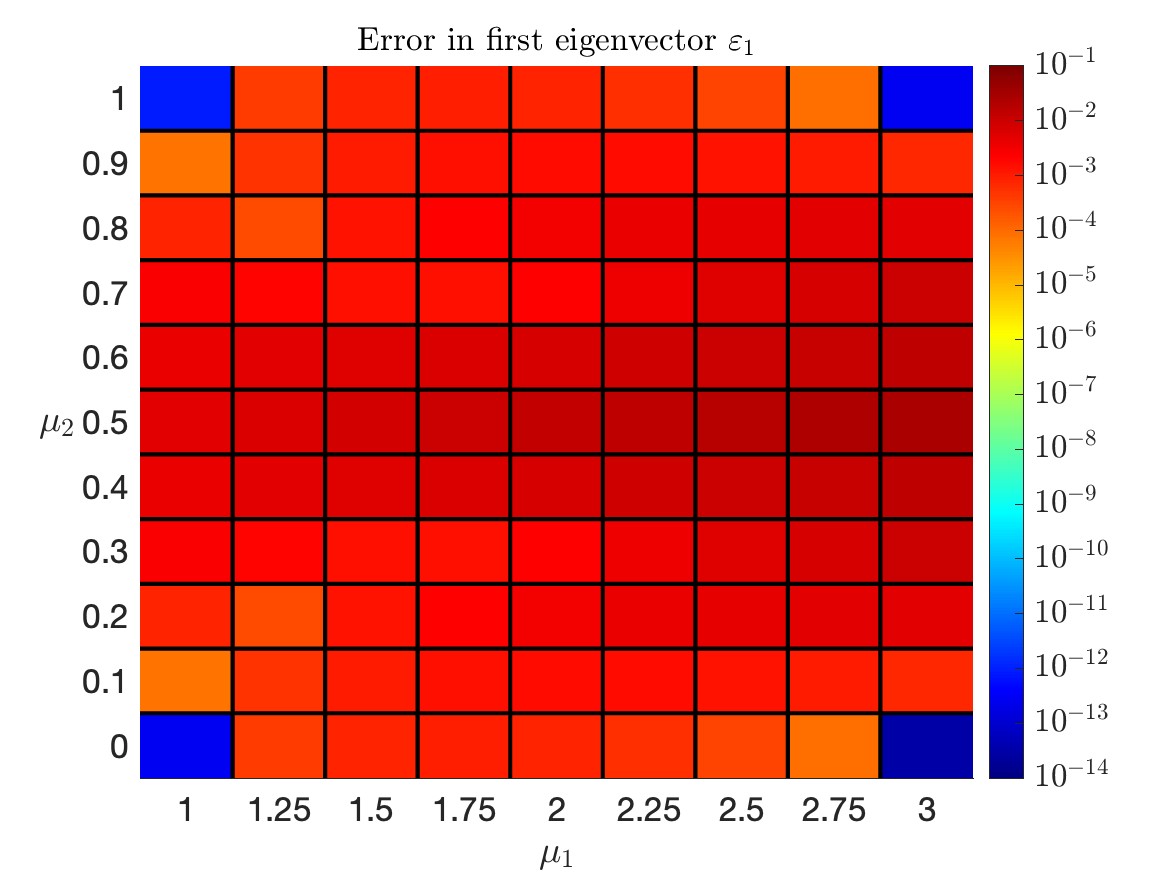}
\includegraphics[width=0.48\linewidth]{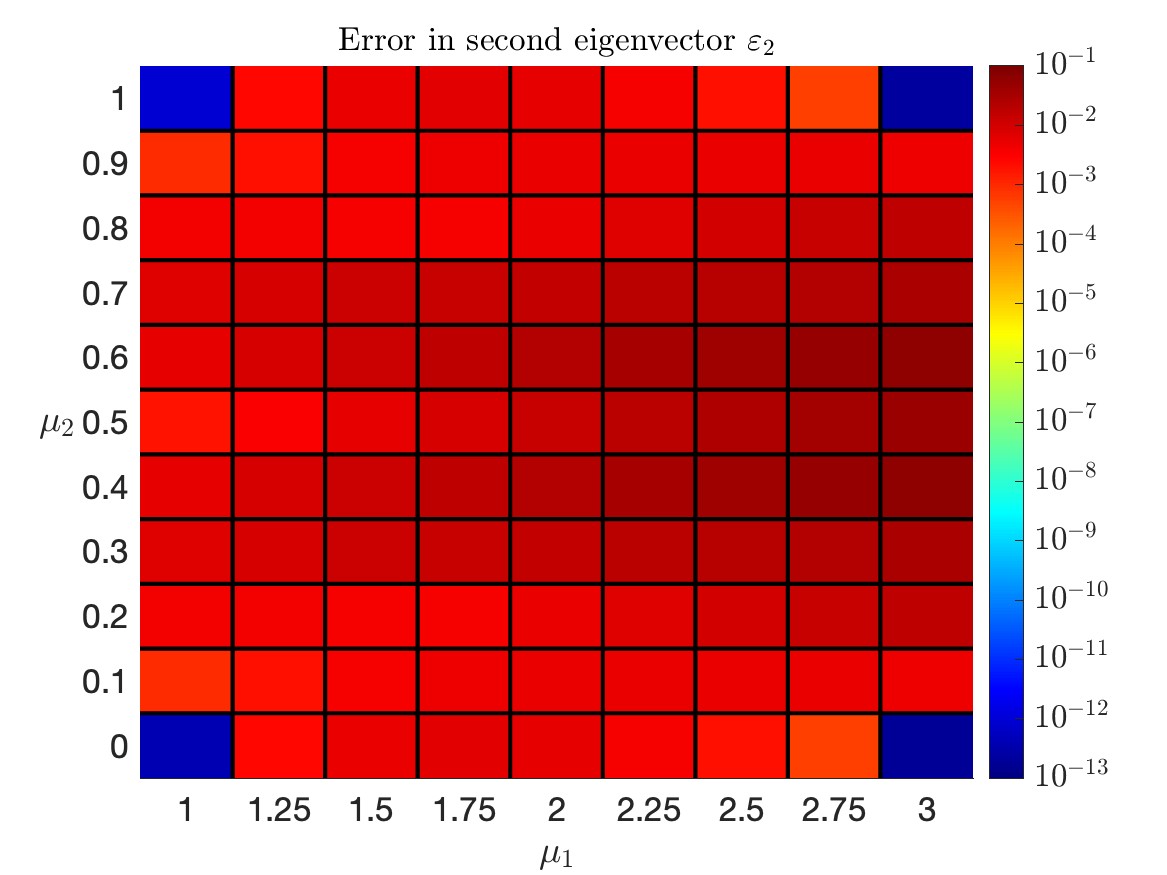}\\
\includegraphics[width=0.48\linewidth]{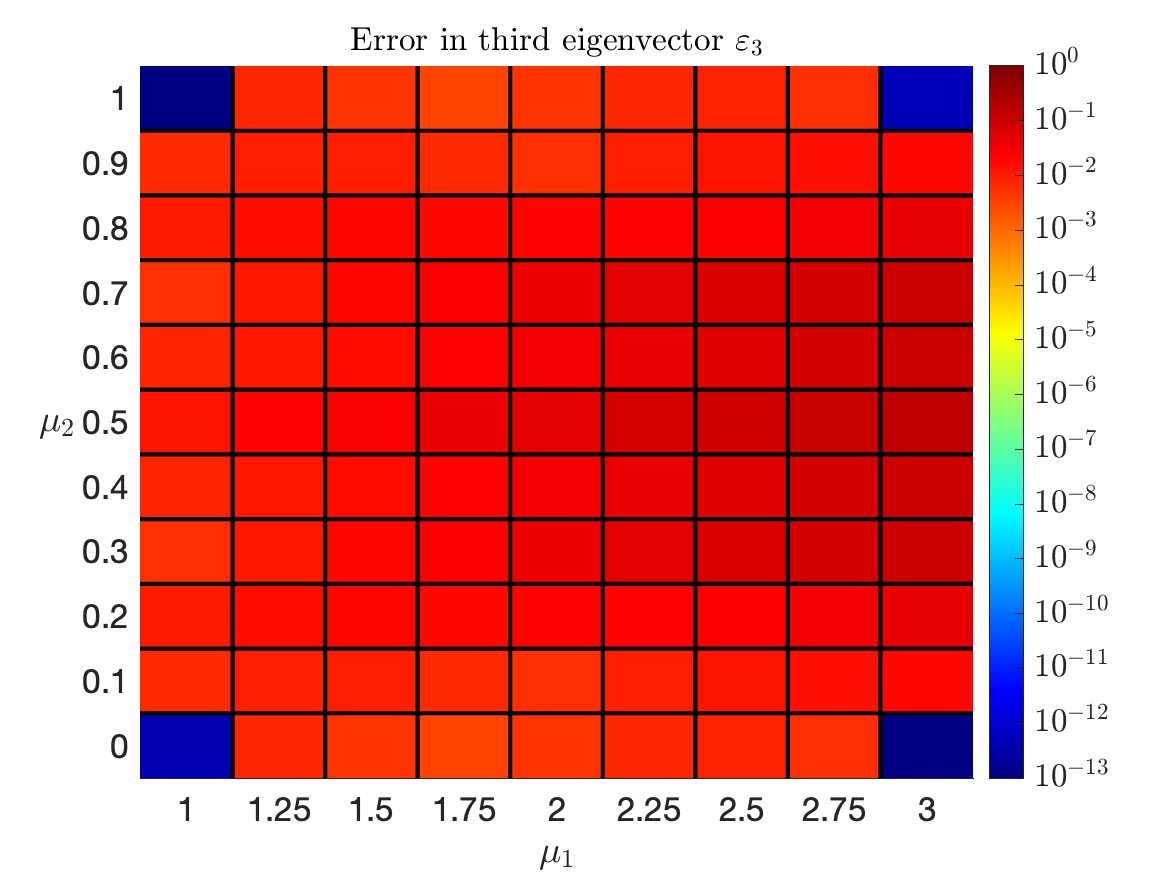}
\includegraphics[width=0.48\linewidth]{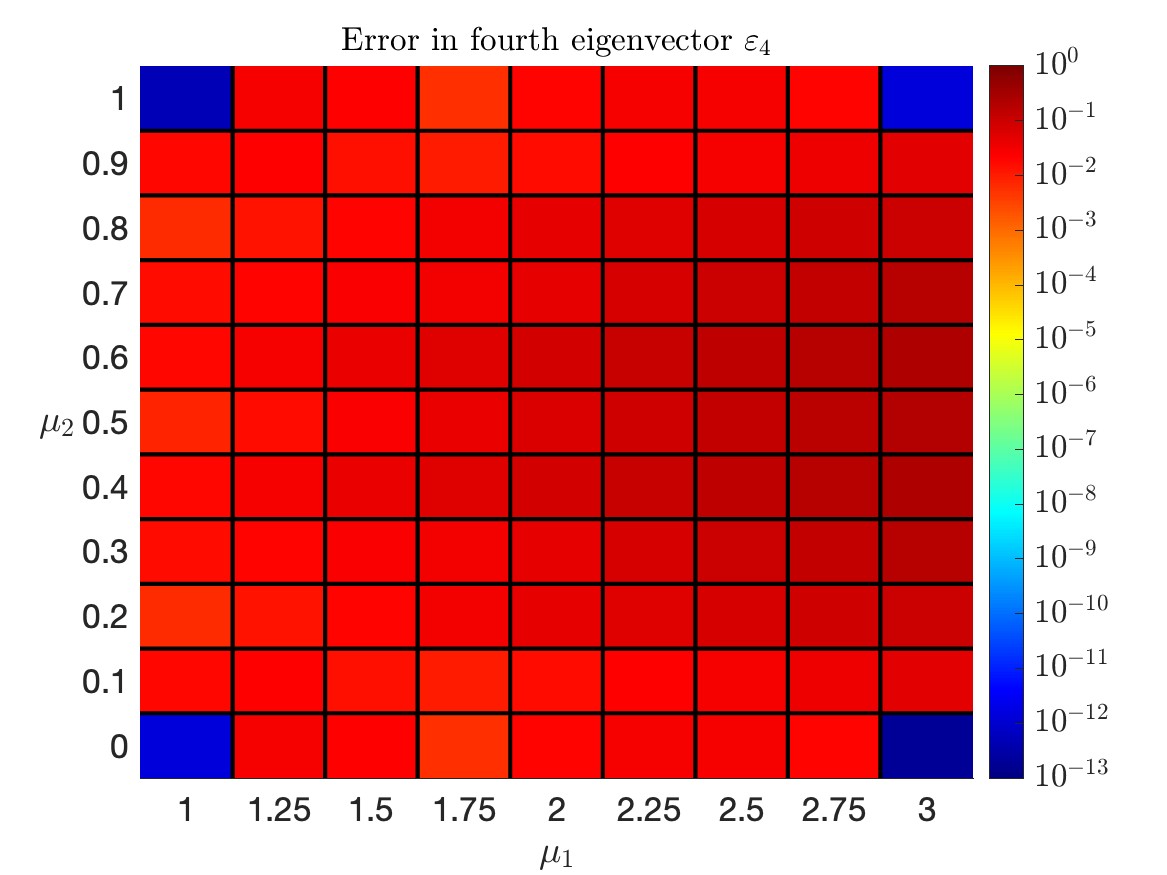}\\
\caption{Error in first 4 eigenvectors at different testing parameters 
\rev{$\boldsymbol{\mu} \in \textsf{D}_{\text{test}}$} 
with coarse sampling of training parameters $\textsf{D}_{\text{train}} = \{ 1, 3 \} \times \{ 0, 1 \}$
in the example of parametric harmonic oscillator.}
\label{fig:schrodinger_harmonic_eigenvector_error_coarse}
\end{figure}

\revv{
We investigate the effects of sampling for collecting the sample snapshots 
by comparing results obtained using two different sets of sample parameters: 
\begin{itemize}
\item $\textsf{D}_{\text{train}} = \{ 1, 3 \} \times \{ 0, 1 \}$ consisting of two sample parameters with a density of $(2, 1)$, and 
\item $\textsf{D}_{\text{train}} = \{ 1, 2, 3 \} \times \{ 0, 0.5, 1 \}$ consisting of three sample parameters with a density of $(1, 0.5)$.
\end{itemize}
}
\rev{
As the first trial, 
the eigenvectors $\{ \boldsymbol{\phi}_k(\boldsymbol{\mu}) \}_{k=1}^{4}$ at the sample parameters 
$\boldsymbol{\mu} \in \textsf{D}_{\text{train}} = \{ 1, 3 \} \times \{ 0, 1 \}$ are used as snapshots 
to construct the basis matrix $\mathbf{Q}$ of column size $r = 16$. 
The projected system \eqref{eq:rom-evp} is solved with the built-in function \texttt{eig} in MATLAB.
Since all eigenvalues are \rev{simple}, direct comparison can be made to the FOM and ROM eigenvectors of the same index. 
Figure~\ref{fig:schrodinger_harmonic_error_bound_coarse} illustrates the behavior of the ROM approximated $k$-th eigenvalue $\widetilde{\lambda}_k$ and 
the ROM approximation error of the $k$-th eigenvector $\varepsilon_k$, for $1 \leq k \leq r = 36$ at the parameter $\boldsymbol{\mu} = (3, 0.8)$, 
which again verifies the validity of our theoretical error estimates. 
Figure~\ref{fig:schrodinger_harmonic_eigenvalue_error_coarse} and Figure~\ref{fig:schrodinger_harmonic_eigenvector_error_coarse} 
show the error in the first 4 eigenvalues and eigenvectors at different testing parameters 
$\boldsymbol{\mu} \in \textsf{D}_{\text{test}} = \{ (1 + 0.25 t_1, 0.1 t_2) : t_1 \in [0,8] \cap \mathbb{Z}, t_2 \in [0, 10] \cap \mathbb{Z} \}$. 
In this example, both approximations are extremely accurate in the reproductive case, i.e. at the parameters 
$\boldsymbol{\mu} \in \textsf{D}_{\text{train}}$, 
where the errors are of order $10^{-12}$. 
However, in the interpolation cases, where $\boldsymbol{\mu} \in \textsf{D} \setminus \textsf{D}_{\text{train}}$, 
the accuracy becomes the worse. 
The error in the $4$-th eigenvector is of order $10^{0}$. 
}

\begin{figure}[ht!]
\centering
\includegraphics[width=0.48\linewidth]{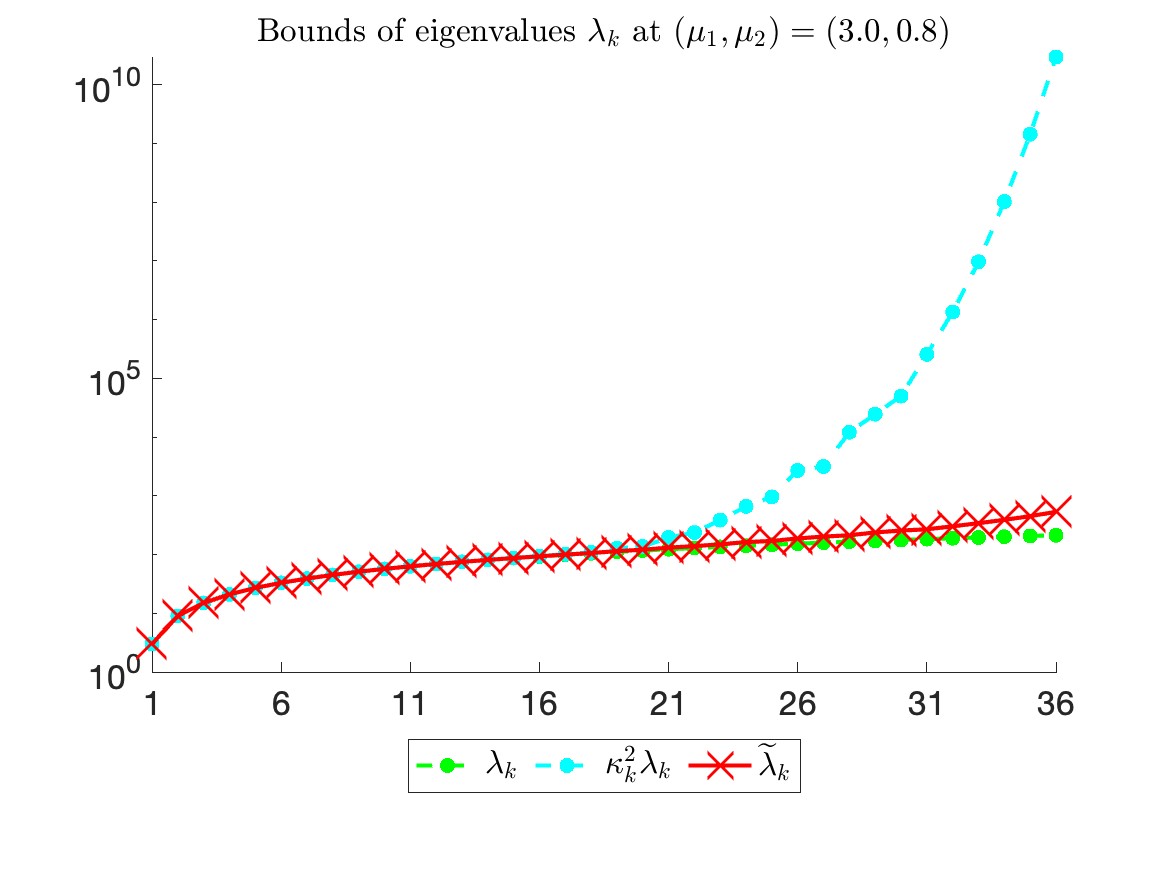}
\includegraphics[width=0.48\linewidth]{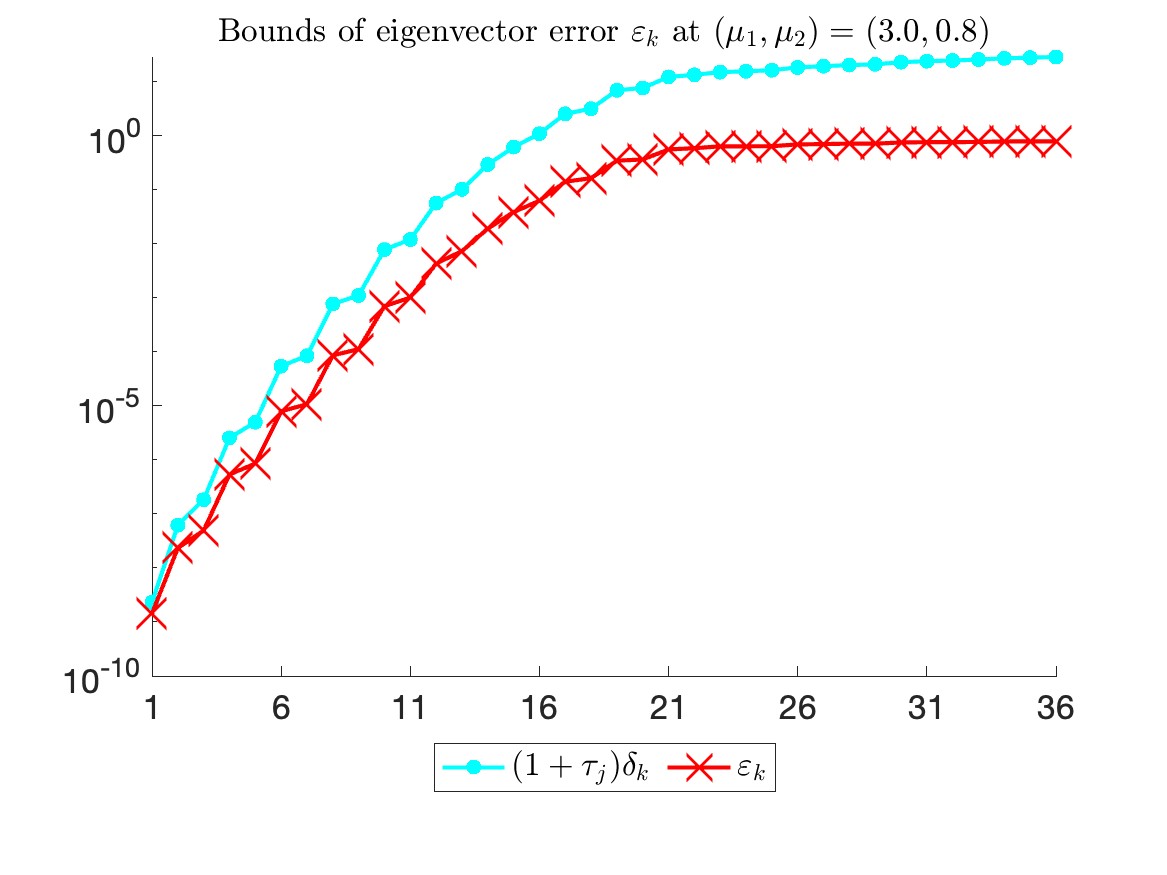}
\caption{Bounds of all 36 eigenvalues (left) and eigenvectors (right) at the parameter  
$\boldsymbol{\mu} = (3, 0.8)$ 
with fine sampling of training parameters $\textsf{D}_{\text{train}} = \{ 1, 2, 3 \} \times \{ 0, 0.5, 1 \}$
in the example of parametric harmonic oscillator.}
\label{fig:schrodinger_harmonic_error_bound_fine}
\end{figure}

\begin{figure}[ht!]
\centering
\includegraphics[width=0.48\linewidth]{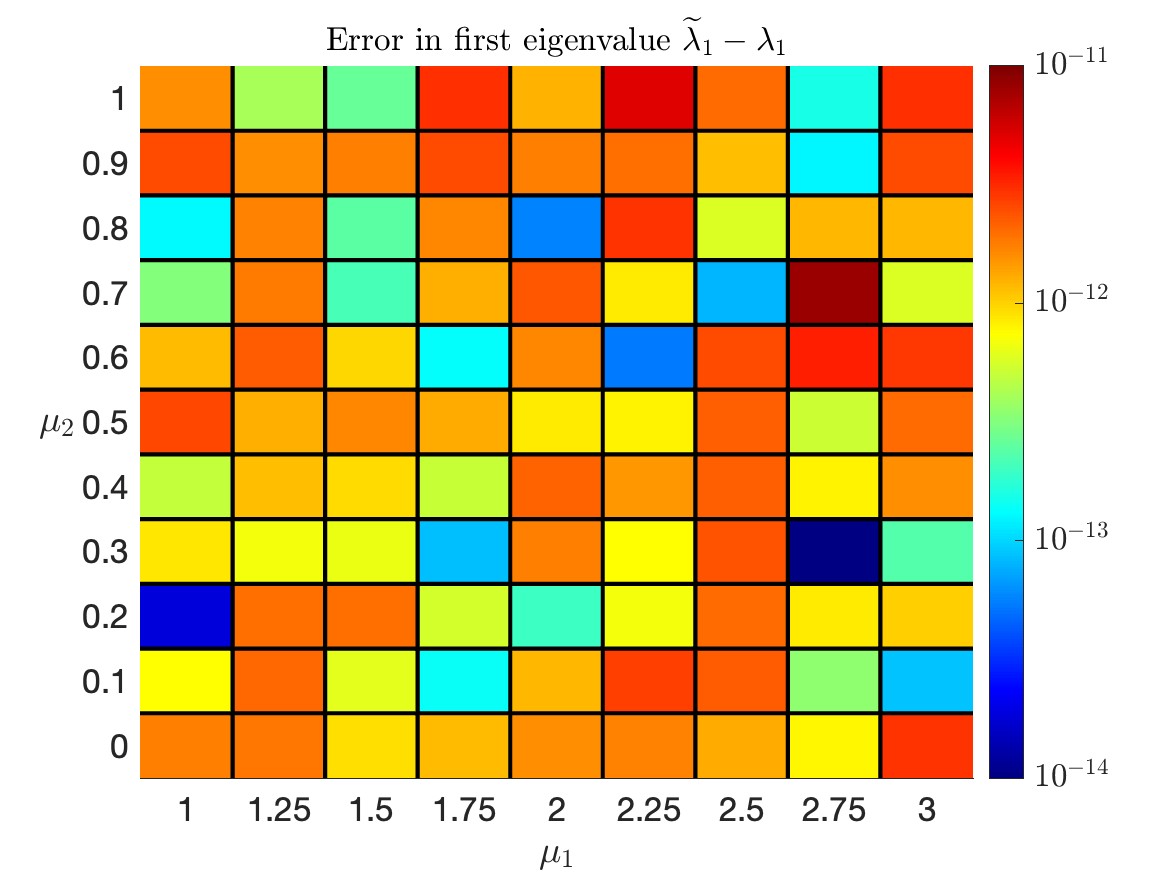}
\includegraphics[width=0.48\linewidth]{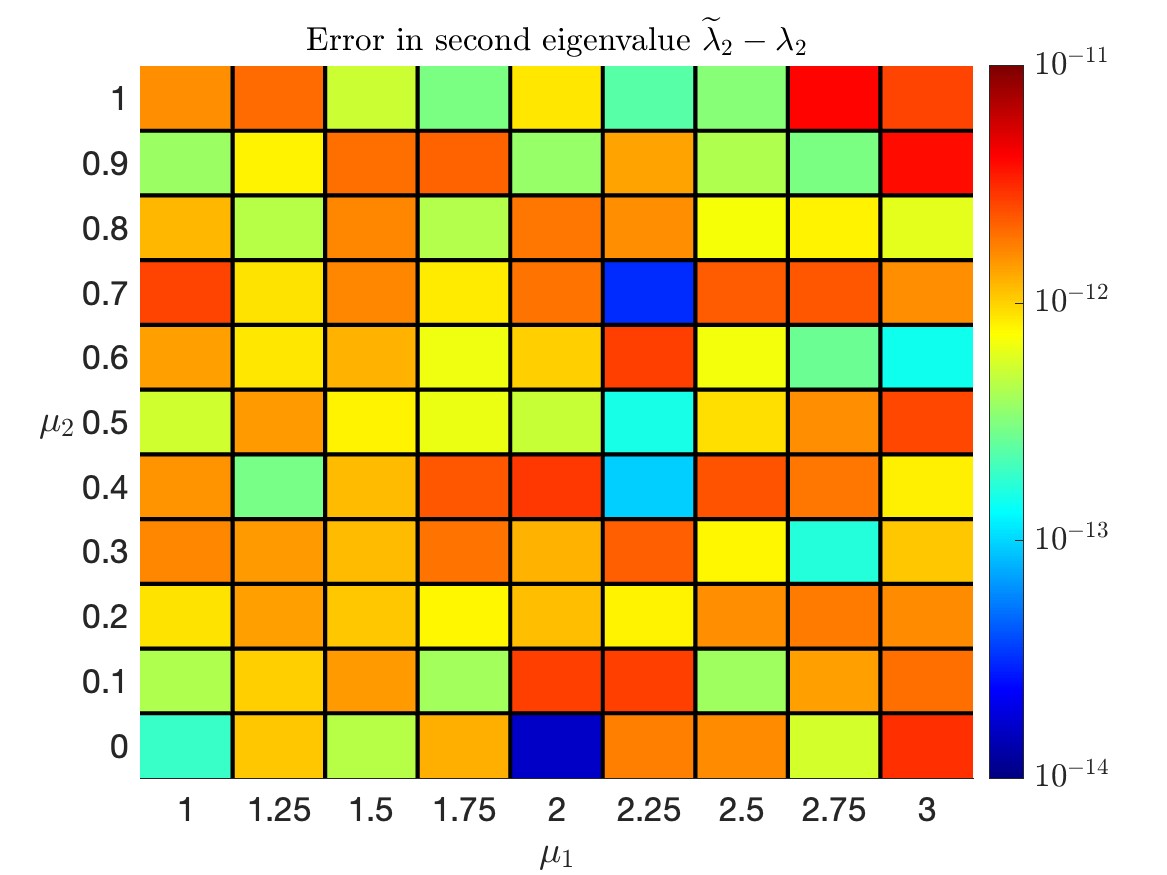}\\
\includegraphics[width=0.48\linewidth]{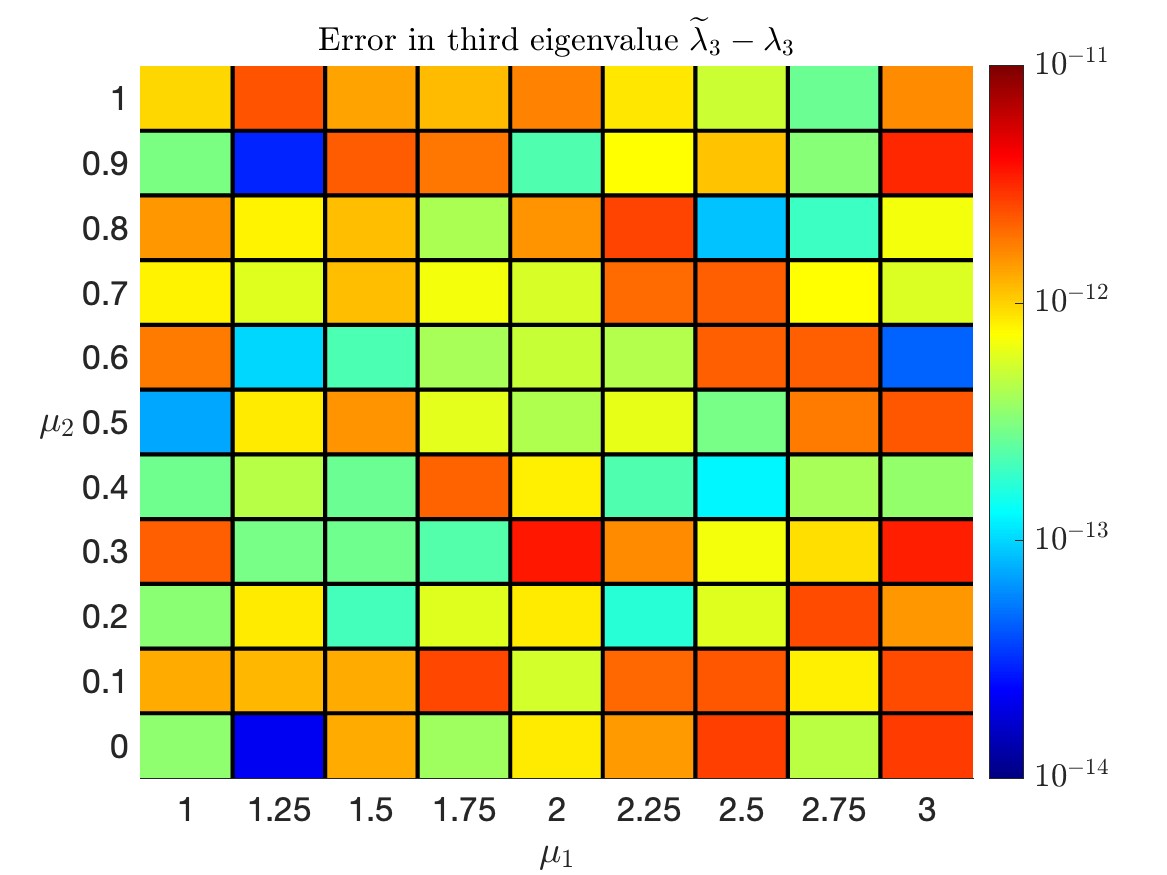}
\includegraphics[width=0.48\linewidth]{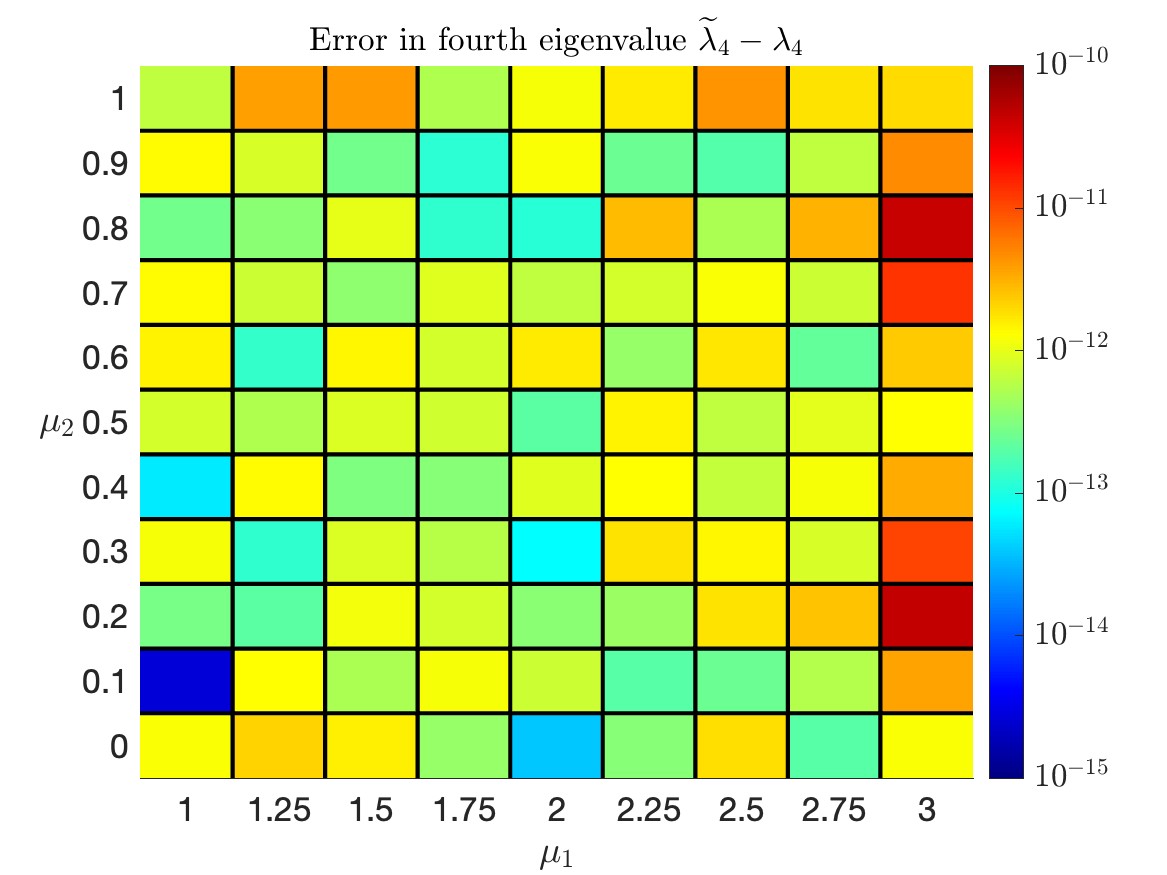}\\
\caption{Error in first 4 eigenvalues at different testing parameters 
\rev{$\boldsymbol{\mu} \in \textsf{D}_{\text{test}}$} 
with fine sampling of training parameters $\textsf{D}_{\text{train}} = \{ 1, 2, 3 \} \times \{ 0, 0.5, 1 \}$
in the example of parametric harmonic oscillator.}
\label{fig:schrodinger_harmonic_eigenvalue_error_fine}
\end{figure}

\begin{figure}[ht!]
\centering
\includegraphics[width=0.48\linewidth]{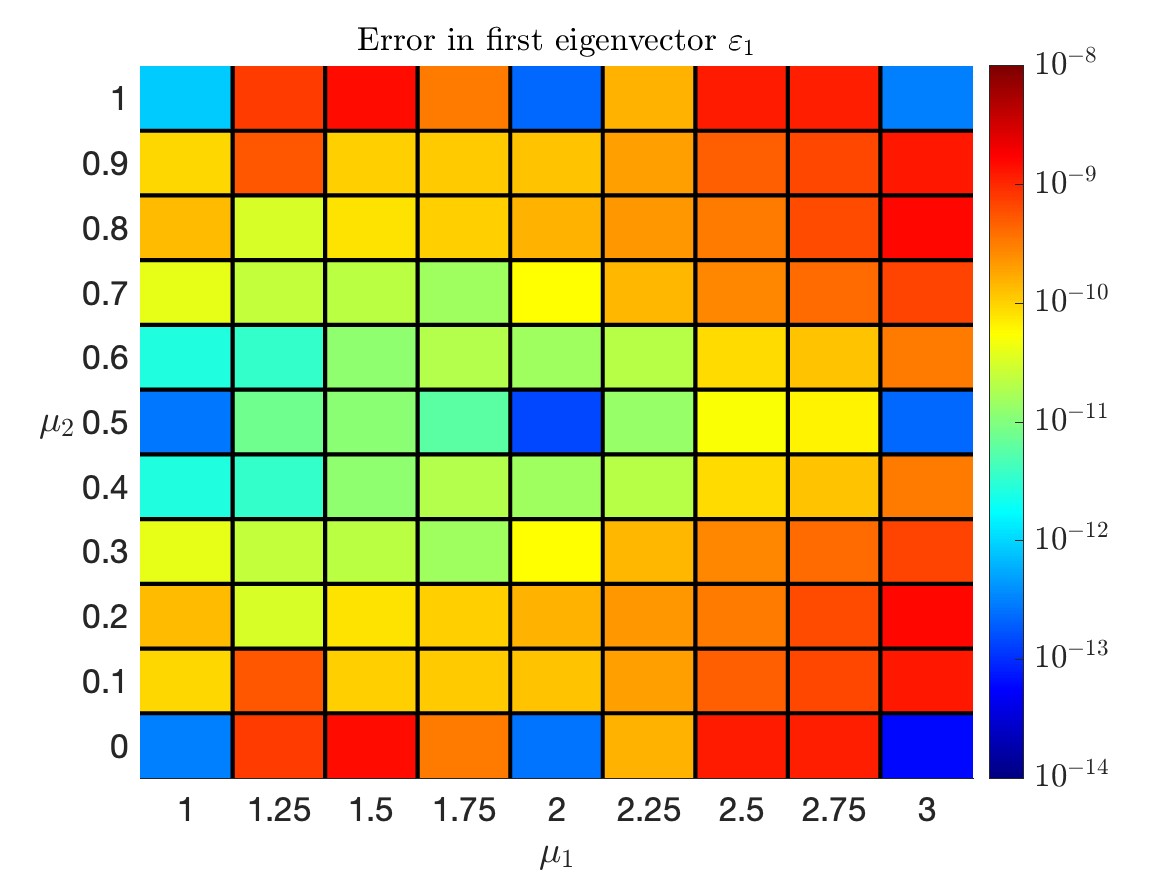}
\includegraphics[width=0.48\linewidth]{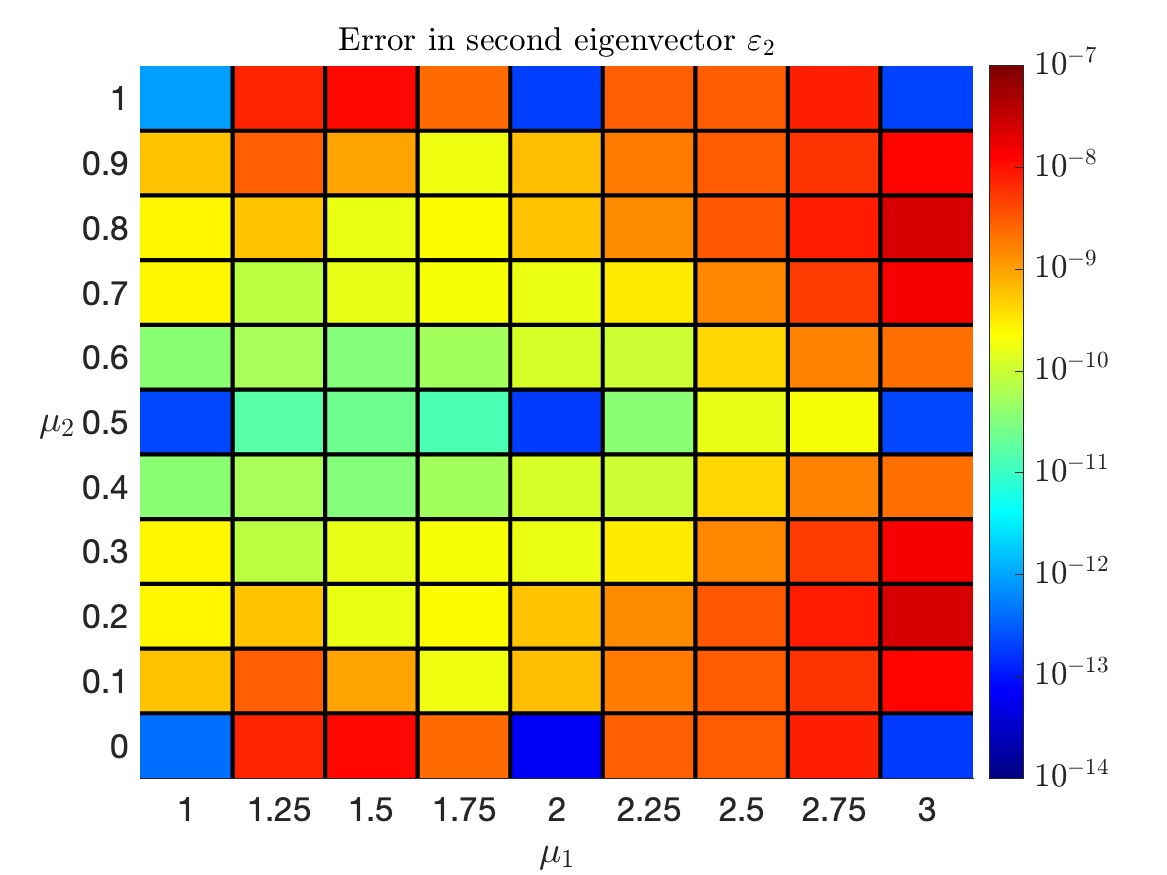}\\
\includegraphics[width=0.48\linewidth]{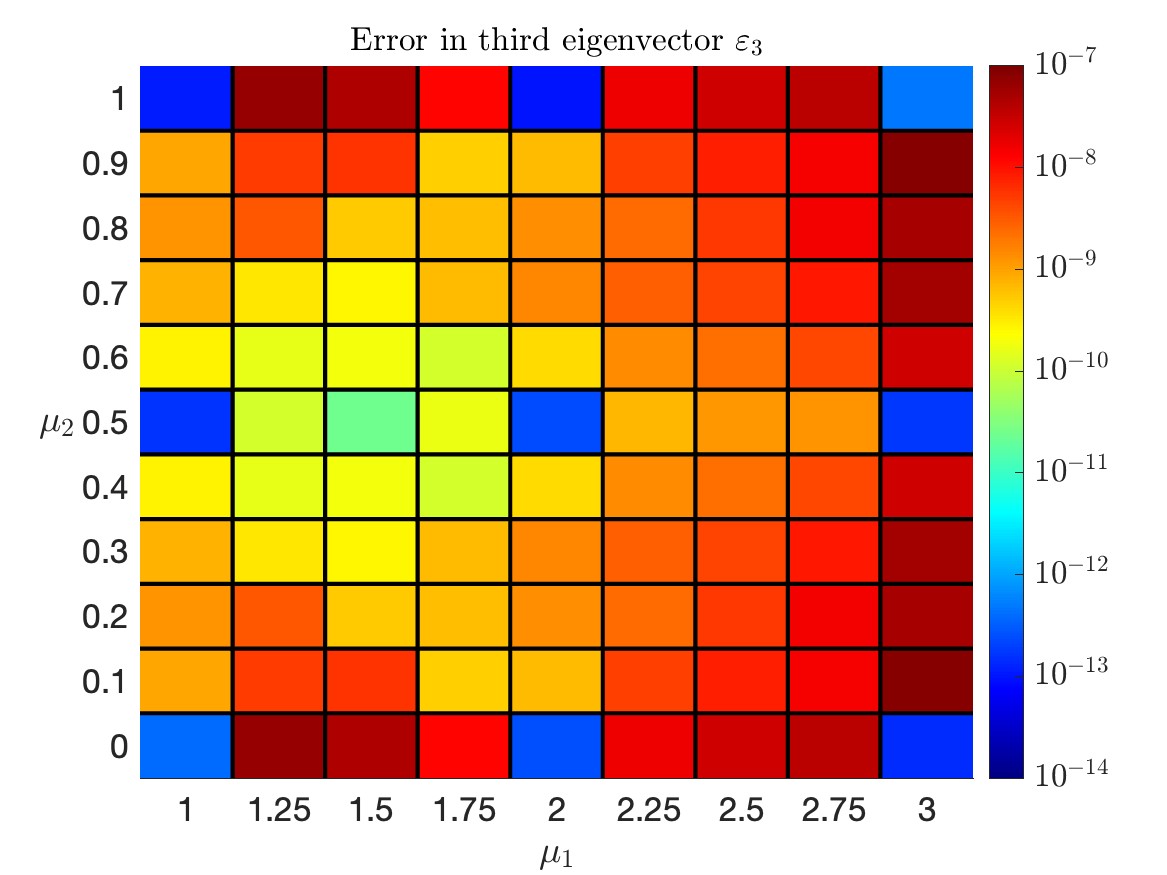}
\includegraphics[width=0.48\linewidth]{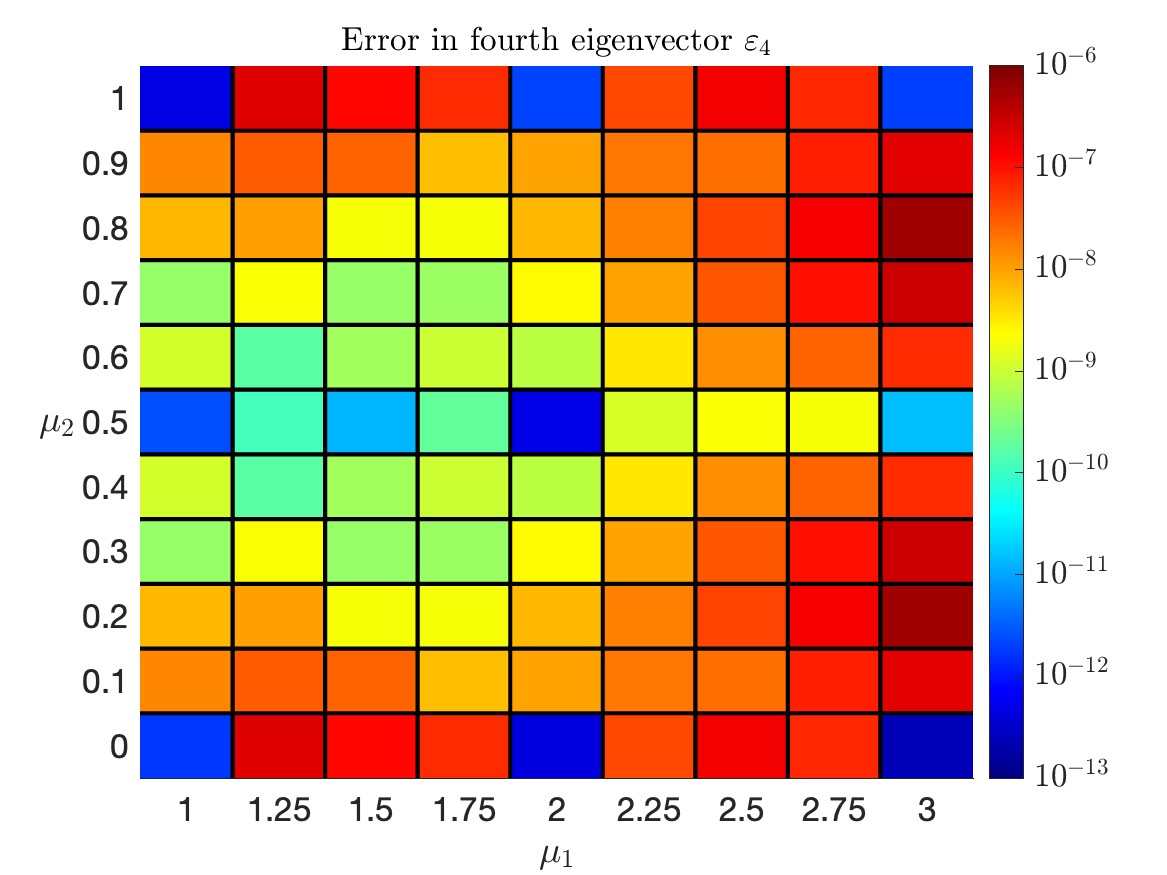}\\
\caption{Error in first 4 eigenvectors at different testing parameters 
\rev{$\boldsymbol{\mu} \in \textsf{D}_{\text{test}}$} 
with fine sampling of training parameters $\textsf{D}_{\text{train}} = \{ 1, 2, 3 \} \times \{ 0, 0.5, 1 \}$
in the example of parametric harmonic oscillator.}
\label{fig:schrodinger_harmonic_eigenvector_error_fine}
\end{figure}

\rev{
To improve the results, we increase the density of the sample parameters used for training. 
The testing procedure is repeated with the sample parameters 
$\boldsymbol{\mu} \in \textsf{D}_{\text{train}} = \{ 1, 2, 3 \} \times \{ 0, 0.5, 1 \}$,
leading to a basis matrix $\mathbf{Q}$ of column size $r = 36$. 
Figure~\ref{fig:schrodinger_harmonic_error_bound_fine} shows that while the validity of our theoretical error estimates remains, 
the improved quality of the basis matrix $\mathbf{Q}$ leads to a 
the more moderate growth of the upper bounds than that 
in Figure~\ref{fig:schrodinger_harmonic_error_bound_coarse} using the coarse sampling. 
Figure~\ref{fig:schrodinger_harmonic_eigenvalue_error_fine} and Figure~\ref{fig:schrodinger_harmonic_eigenvector_error_fine} 
show the error in the first 4 eigenvalues and eigenvectors at different testing parameters 
$\boldsymbol{\mu} \in \textsf{D}_{\text{test}}$. 
In this example, both approximations are extremely accurate at all the testing parameters, 
where the errors are of order $10^{-6}$ or below. 
}

\begin{figure}[ht!]
\centering
\includegraphics[width=0.48\linewidth]{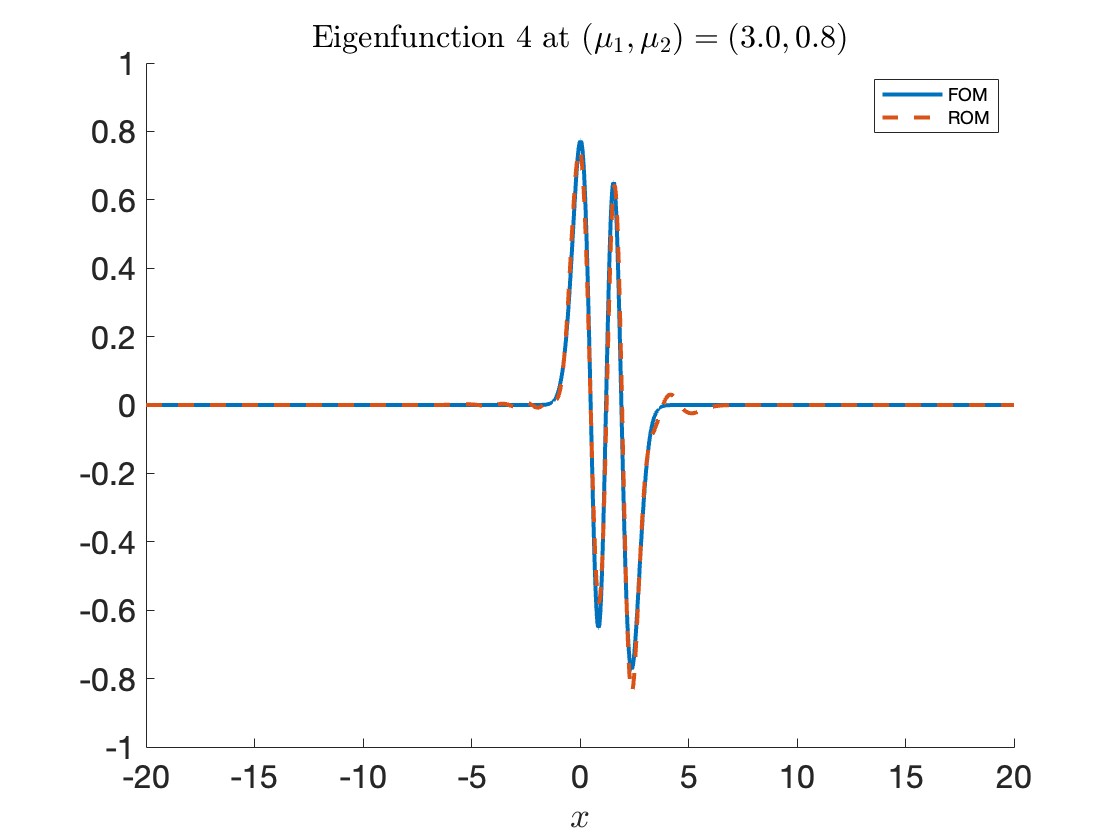}
\includegraphics[width=0.48\linewidth]{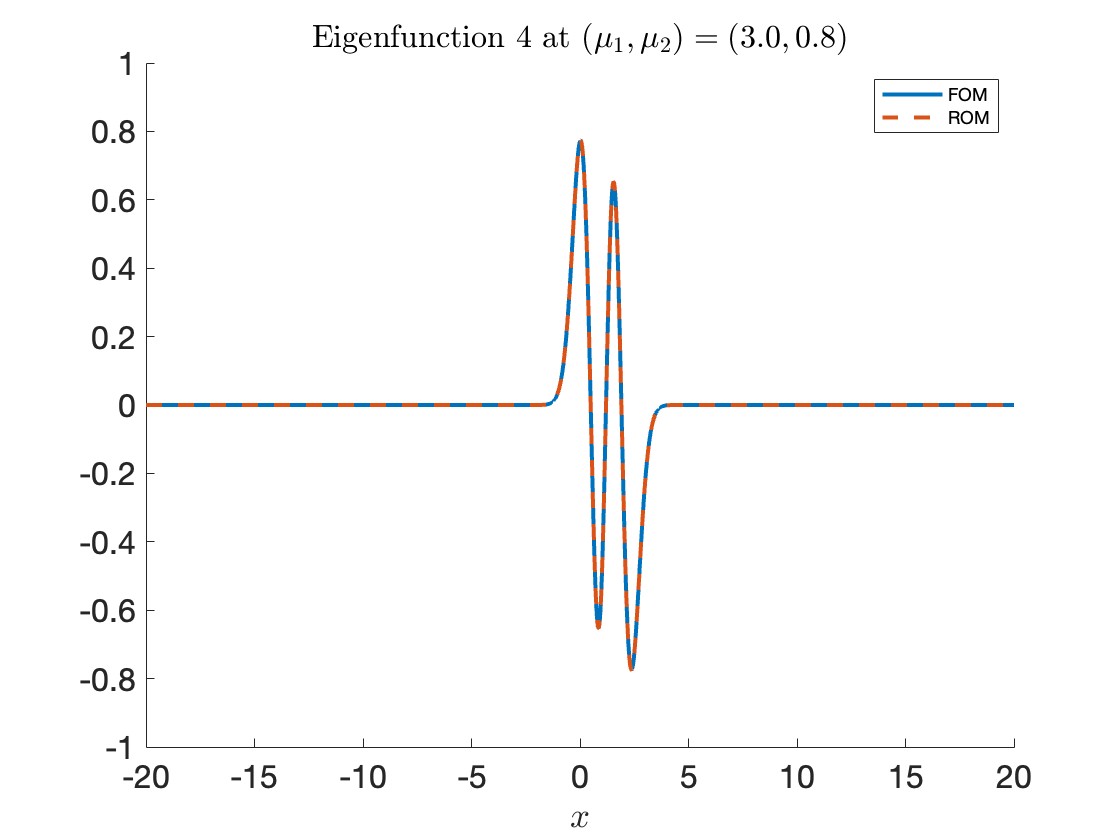}
\caption{Finite element representations of the FOM solution (blue solid line) and the ROM solution (red dashed line) of 
the 4-th eigenvectors at the parameter $\boldsymbol{\mu} = (3, 0.8)$ with
coarse sampling $\textsf{D}_{\text{train}} = \{ 1, 3 \} \times \{ 0, 1 \}$ (left) and 
fine sampling $\textsf{D}_{\text{train}} = \{ 1, 2, 3 \} \times \{ 0, 0.5, 1 \}$ (right) 
in the example of parametric harmonic oscillator.}
\label{fig:schrodinger_harmonic_eigenfunction}
\end{figure}

\rev{
Figure~\ref{fig:schrodinger_harmonic_eigenfunction} depicts the finite element representations of the FOM solution and the ROM solution of 
the 4-th eigenvectors at the parameter $\boldsymbol{\mu} = (3, 0.8)$ with
the two densities of the sample parameters used for training. 
With coarse sampling $\textsf{D}_{\text{train}} = \{ 1, 3 \} \times \{ 0, 1 \}$, 
the ROM approximation exhibits overshooting and undershooting at the extrema, 
as well as noticeable fluctuations around $x = 4$.
These inaccuracy was resolved with fine sampling $\textsf{D}_{\text{train}} = \{ 1, 2, 3 \} \times \{ 0, 0.5, 1 \}$. 
}

\subsection{2D parametric Gaussian well potential}

\begin{figure}[ht!]
\centering
\includegraphics[width=0.48\textwidth]{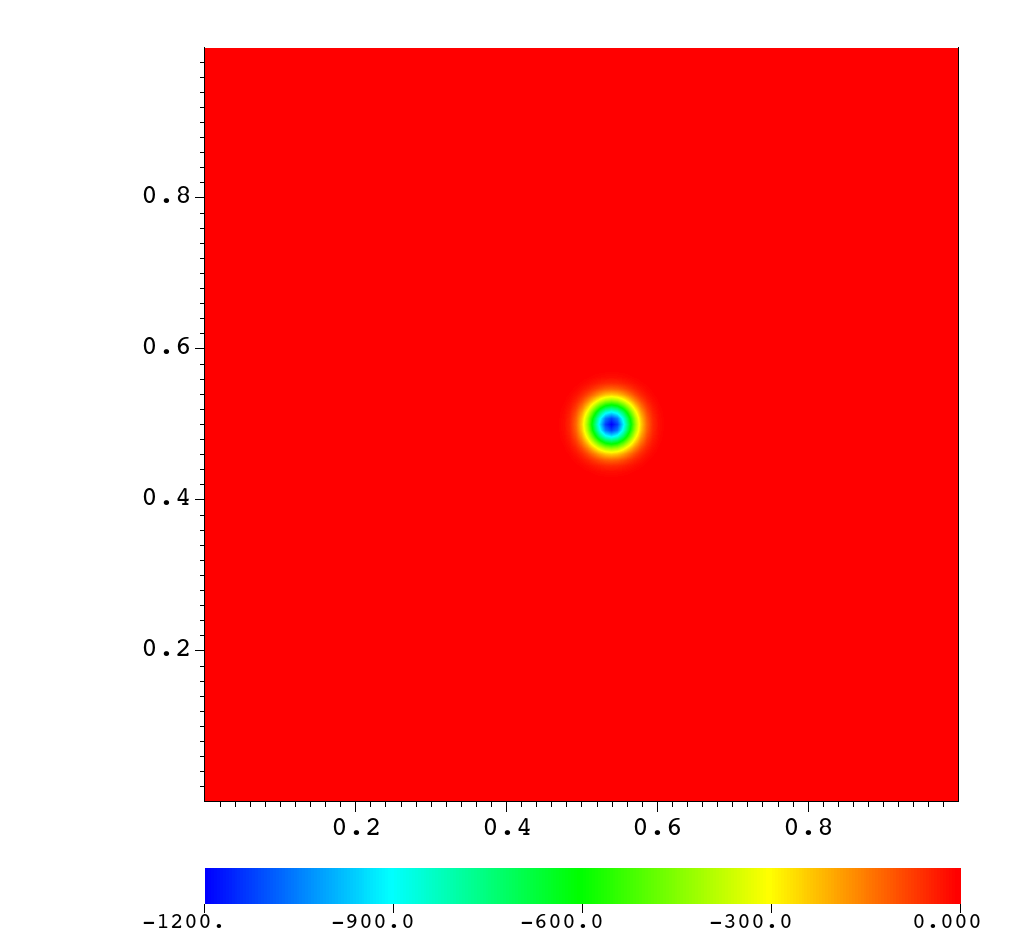}
\caption{Potential coefficient $\rho$ at the parameter $\mu = 5$
 in the example of parametric Gaussian well potential.}
\label{fig:schrodinger_square_potential}
\end{figure}

In this example, the computational domain is taken as the unit square 
$\Omega = (0,1)^2 \subset \mathbb{R}^2$, and the parameter domain is 
$\textsf{D} = [-5, 5]$ 
\rev{which reduces the parameter $\boldsymbol{\mu}$ to a scalar $\mu$}.
For $\mu \in \textsf{D}$, 
the coefficient fields are given by 
\begin{equation}
\begin{aligned}
\sigma(x; \mu) & = 1 && \text{for } x \in \Omega, \\
\rho(x; \mu) & = -1200 \exp\left( -1024 (x - x_\circ(\mu))^2 \right) && \text{for } x \in \Omega, 
\end{aligned}
\end{equation}
which corresponds to the Schr\"{o}dinger operator with a parametric Gaussian well potential
with $x_\circ(\mu) = (0.5 + \mu / 128, 0.5)$ being the center of the Gaussian well, and 
the Neumann boundary condition is prescribed, i.e., 
\begin{equation}
\alpha(x; \mu) = 0, \quad \beta(x; \mu) = 1, \quad \text{for } x \in \partial\Omega.
\end{equation}
Figure~\ref{fig:schrodinger_square_potential} depicts the potential coefficient $\rho$ 
in the domain $\Omega$ at the parameter $\mu = 5$.

The domain $\Omega$ is divided into a uniform mesh of size $h = 1/128$, 
and $Q^1$ Lagrange finite elements is used in the finite element discretization, 
which results in a system size of $n = 16641$. 
We are interested in the first $p = 1$ eigenpair, with the eigenvalue being negative.
The eigenvalue problem \eqref{eq:fom-evp} is solved by LOBPCG.

\begin{figure}[ht!]
\centering
\includegraphics[width=0.48\linewidth]{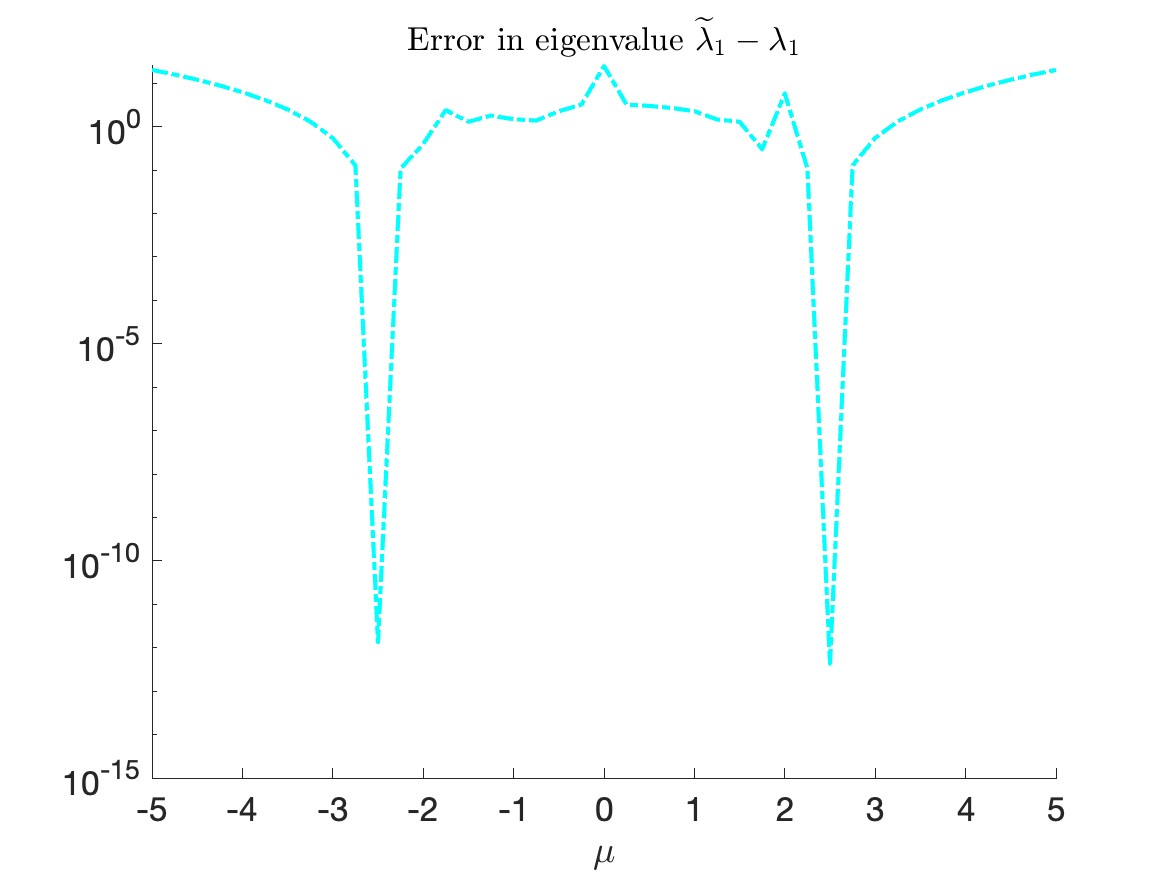}
\includegraphics[width=0.48\linewidth]{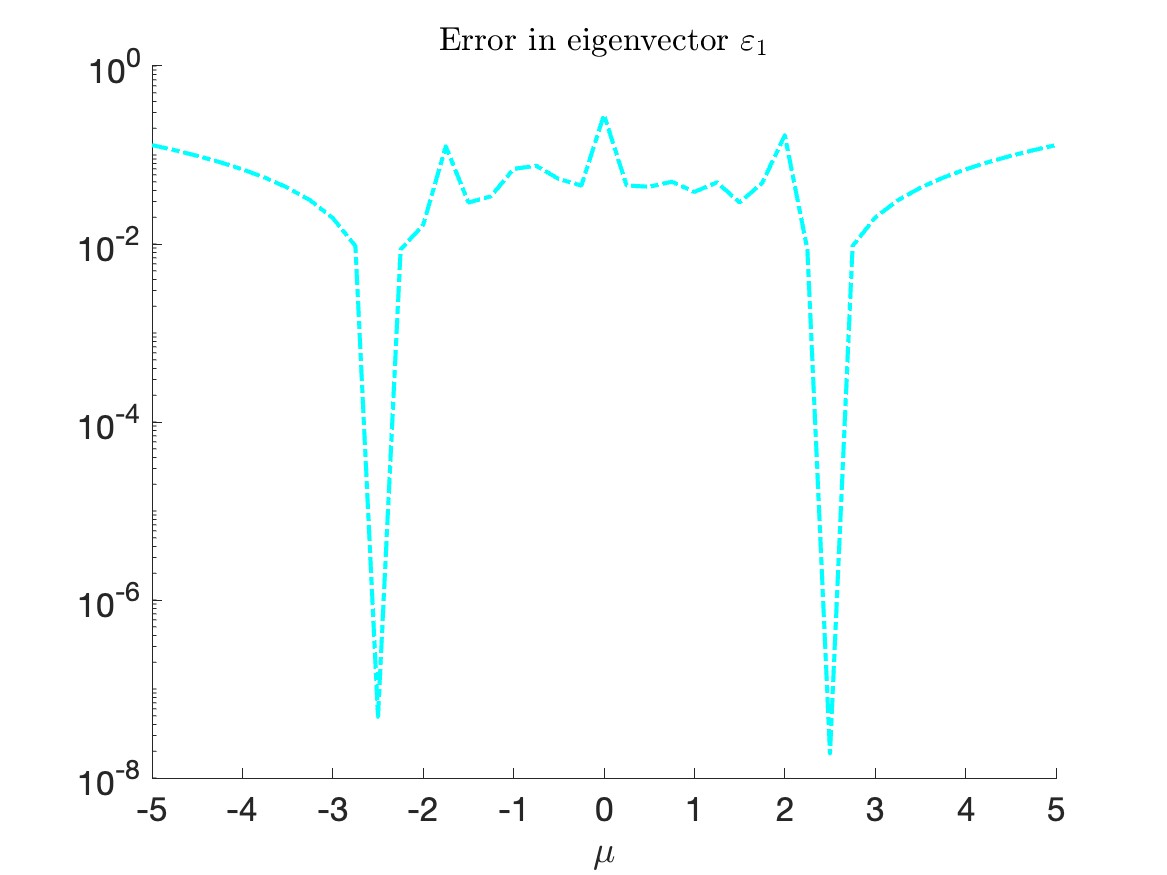}
\caption{Error in lowest eigenvalue (left) and eigenvector (right) at different testing parameters 
\rev{$\mu \in \textsf{D}_{\text{test}}$} in the example of parametric Gaussian well potential.}
\label{fig:schrodinger_square_error}
\end{figure}
\begin{figure}[ht!]
\centering
\includegraphics[width=0.48\textwidth]{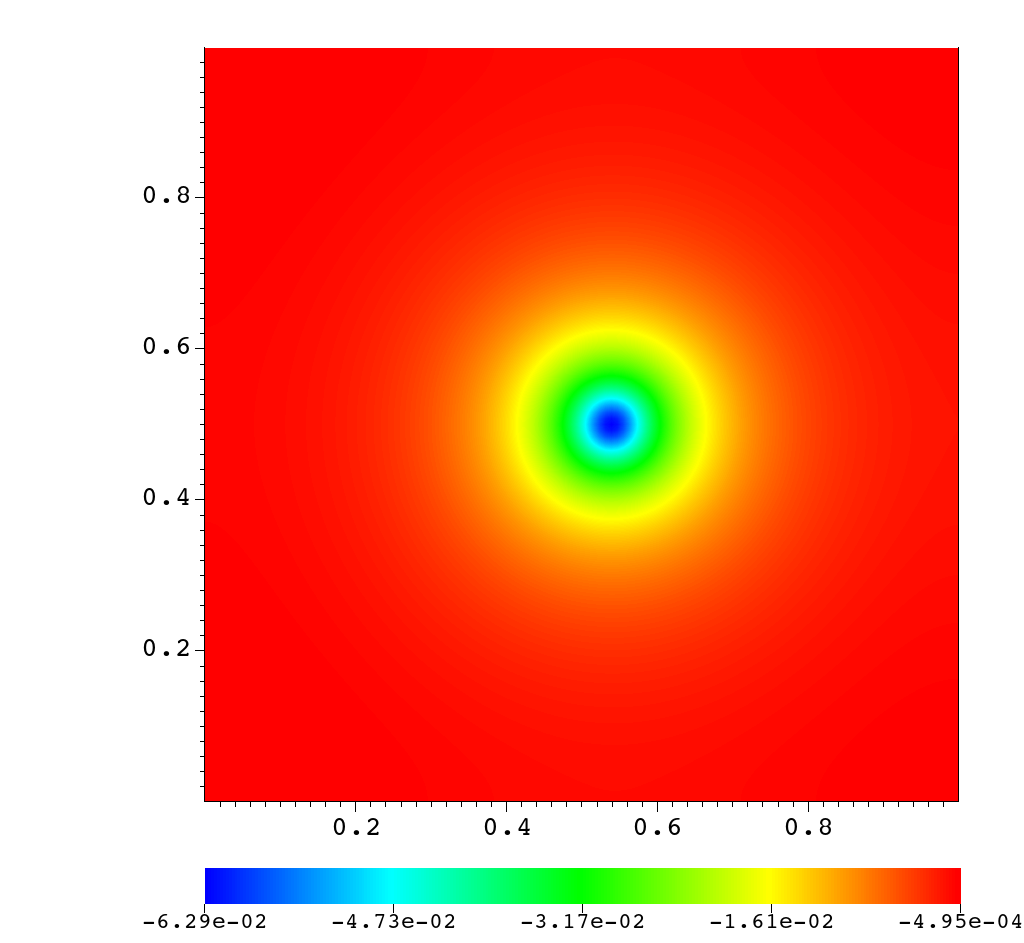}
\includegraphics[width=0.48\textwidth]{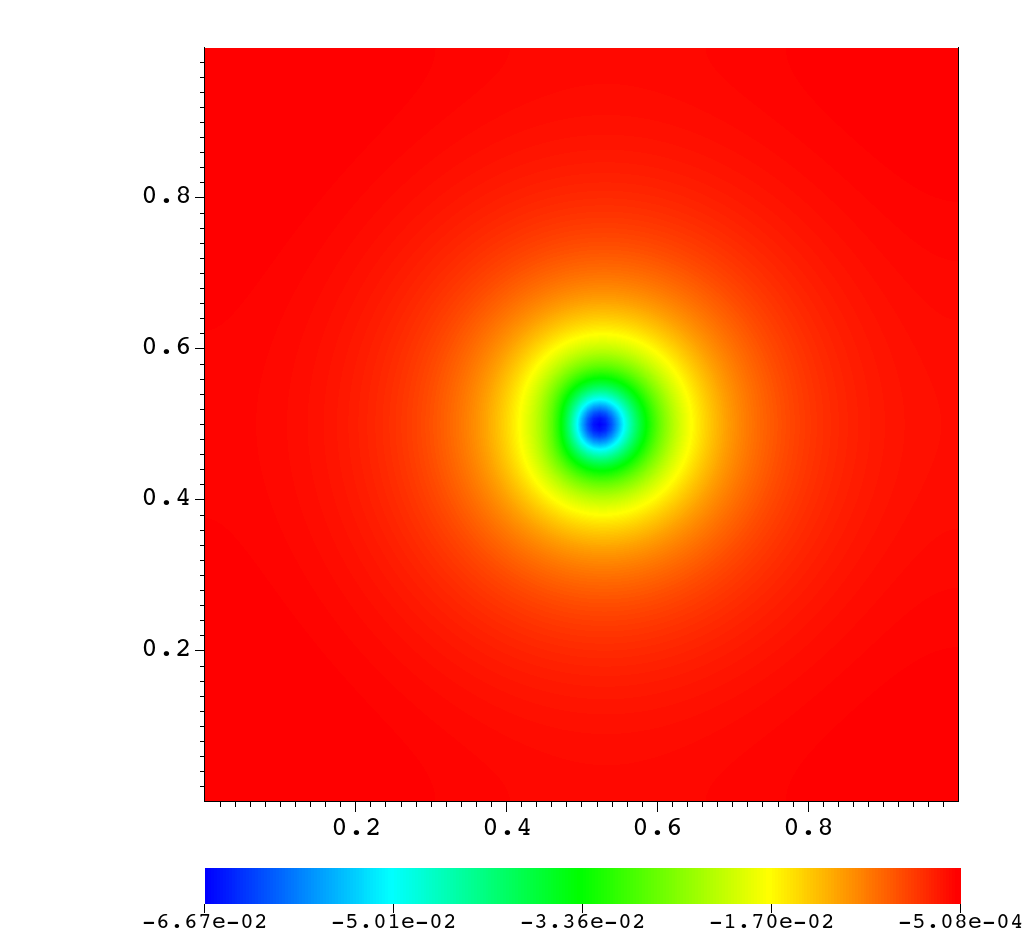}
\caption{Finite element representations of the FOM solution (left) and the ROM solution (right) of 
the first eigenvector at the parameter $\mu = 5$  
in the example of parametric Gaussian well potential.}
\label{fig:schrodinger_square_eigenfunction_1}
\end{figure}

The eigenvector $\boldsymbol{\phi}_1(\mu)$ at the sample parameters $\mu \in \textsf{D}_{\text{train}} = \{ -2.5, 2.5 \}$ are used as snapshots 
to construct the basis matrix $\mathbf{Q}$ of column size $r = 2$.  
The projected system \eqref{eq:rom-evp} is solved with the routine \texttt{dsygv} in LAPACK.
Figure~\ref{fig:schrodinger_square_error} shows the error in the first eigenvalue and eigenvector at different testing parameters 
\rev{$\boldsymbol{\mu} \in \textsf{D}_{\text{test}} = \{ -5 + 0.25 t : t \in [0,40] \cap \mathbb{Z} \}$}.
Figure~\ref{fig:schrodinger_square_eigenfunction_1} depicts the finite element representations of the FOM solution and the ROM solution of 
the first eigenvector with the parameter $\mu = 5$. 
It can be seen that the peak of the ROM solution has the peak shifted to the left  
due to the snapshots sampled, and the decay of the ROM solution becomes asymmetric. 

\subsection{3D parametric diatomic well potential}

In this example, the computational domain is taken as the unit cube 
$\Omega = (0,1)^3 \subset \mathbb{R}^3$, and the parameter domain is 
$\textsf{D} = [-2,2]$ 
\rev{which reduces the parameter $\boldsymbol{\mu}$ to a scalar $\mu$}.
For $\mu \in \textsf{D}$, 
the coefficient fields are given by 
\begin{equation}
\begin{aligned}
\sigma(x; \mu) & = 1/324 && \text{for } x \in \Omega, \\
\rho(x; \mu) & = \sum_{k=1}^2 \sum_{p=1}^2 A_k \exp\left( -\dfrac{(x - x^{(p)}_\circ(\mu))^2}{2 R_k^2} \right)  && \text{for } x \in \Omega, 
\end{aligned}
\end{equation}
with $A_1 = -28.9$, $A_2 = -3.6$, 
$R_1 = 7/450$, $R_2 = 3/50$,
$x^{(1)}_\circ(\mu) = (13/36 - \mu / 128, 1/2, 1/2)$ 
and $x^{(2)}_\circ(\mu) = (23/36 + \mu / 128, 1/2, 1/2)$, 
which corresponds to a Kohn-Sham operator mimicking a diatomic well potential, and 
the Neumann boundary condition is prescribed, i.e., 
\begin{equation}
\alpha(x; \mu) = 0, \quad \beta(x; \mu) = 1, \quad \text{for } x \in \partial\Omega.
\end{equation}

\begin{figure}[ht!]
\centering
\includegraphics[width=0.48\linewidth]{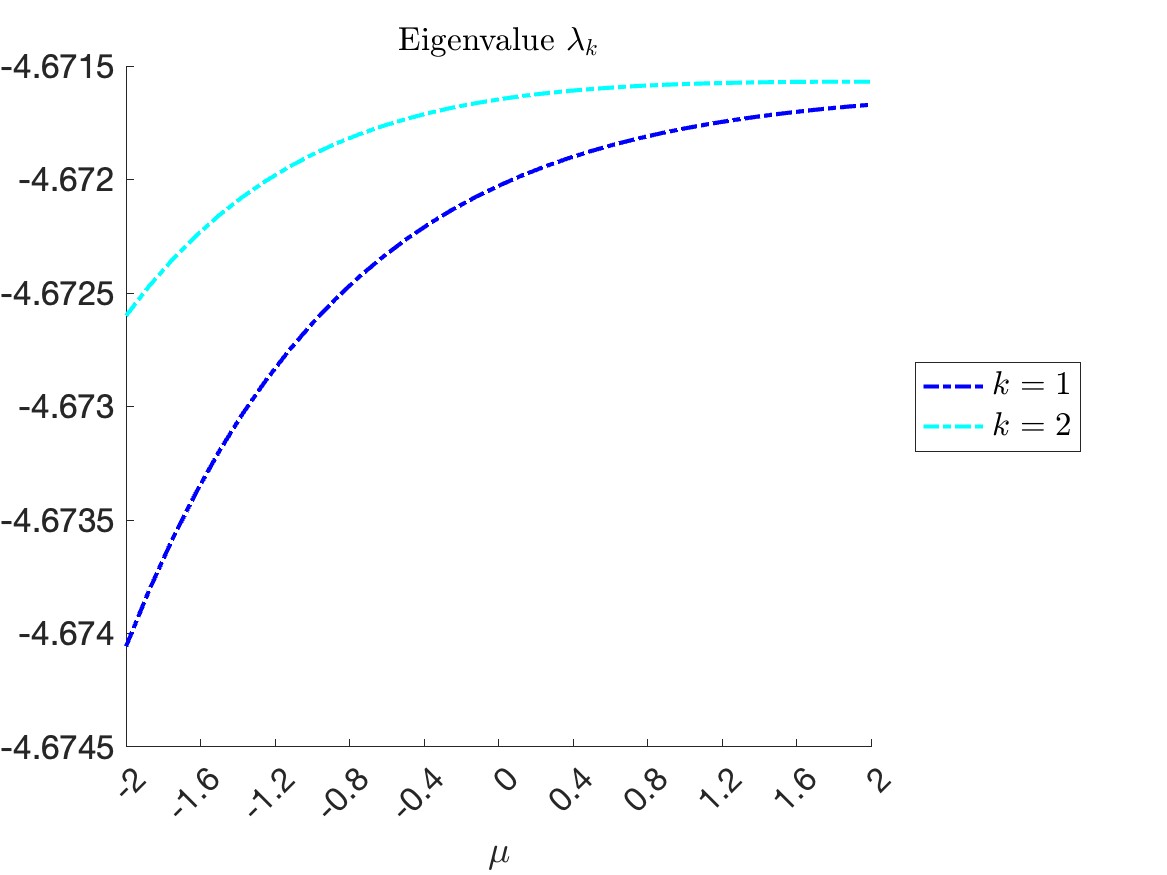}
\caption{First 2 eigenvalues at different parameters $\mu$ in the example of parametric parametric diatomic well potential.}
\label{fig:schrodinger_cube_eigenvalue}
\end{figure}

The domain $\Omega$ is divided into a uniform mesh of size $h = 1/128$, 
and $Q^1$ Lagrange finite elements is used in the finite element discretization, 
which results in a system size of $n = 2146689$. 
We are interested in the first $p = 2$ eigenpairs, with the eigenvalues being negative. 
The eigenvalue problem \eqref{eq:fom-evp} is solved by LOBPCG.
\rev{Figure~\ref{fig:schrodinger_cube_eigenvalue} shows the first 2 eigenvalues $\lambda_k(\mu)$ 
of \eqref{eq:fom-evp} at different parameters $\mu \in \textsf{D}$.}
In this case, the first 2 eigenvalues are both approximately $-4.67$ for $\mu \in \textsf{D}$. 

\begin{figure}[ht!]
\centering
\includegraphics[width=0.48\linewidth]{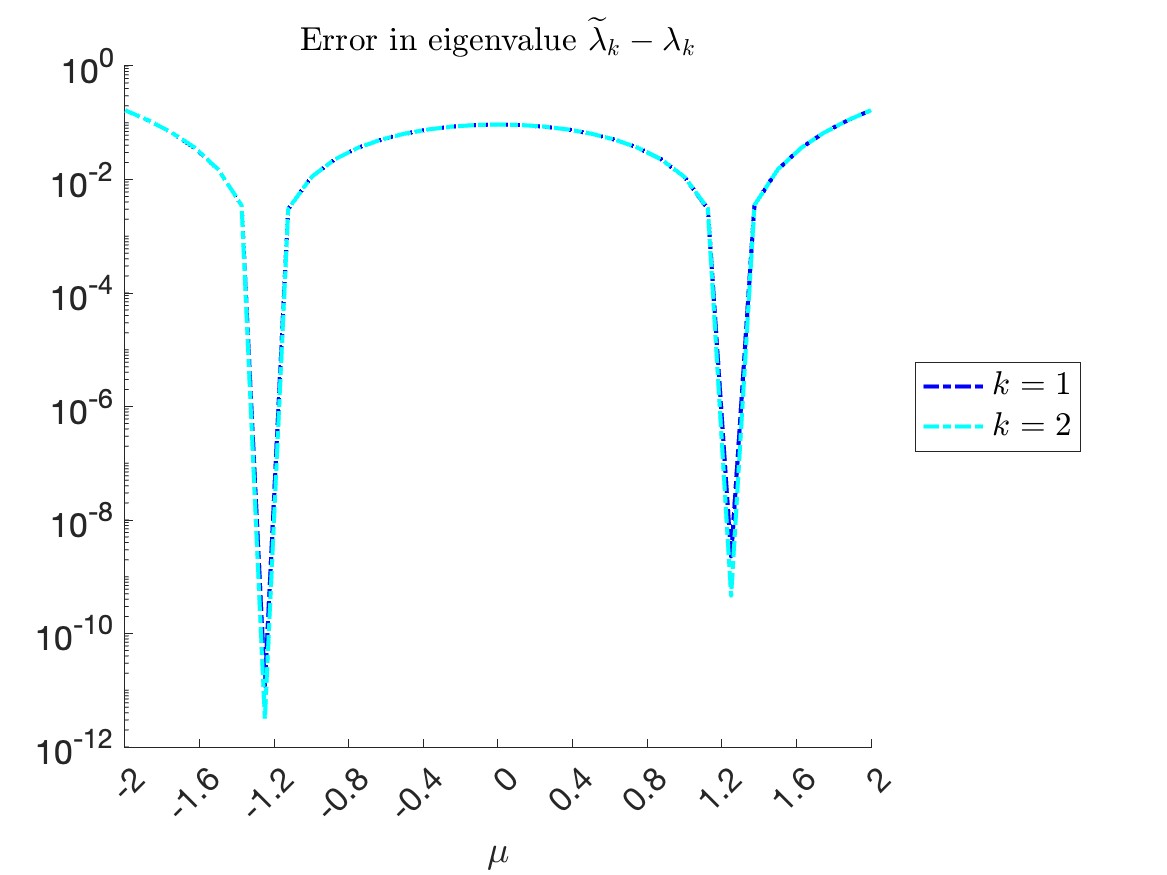}
\includegraphics[width=0.48\linewidth]{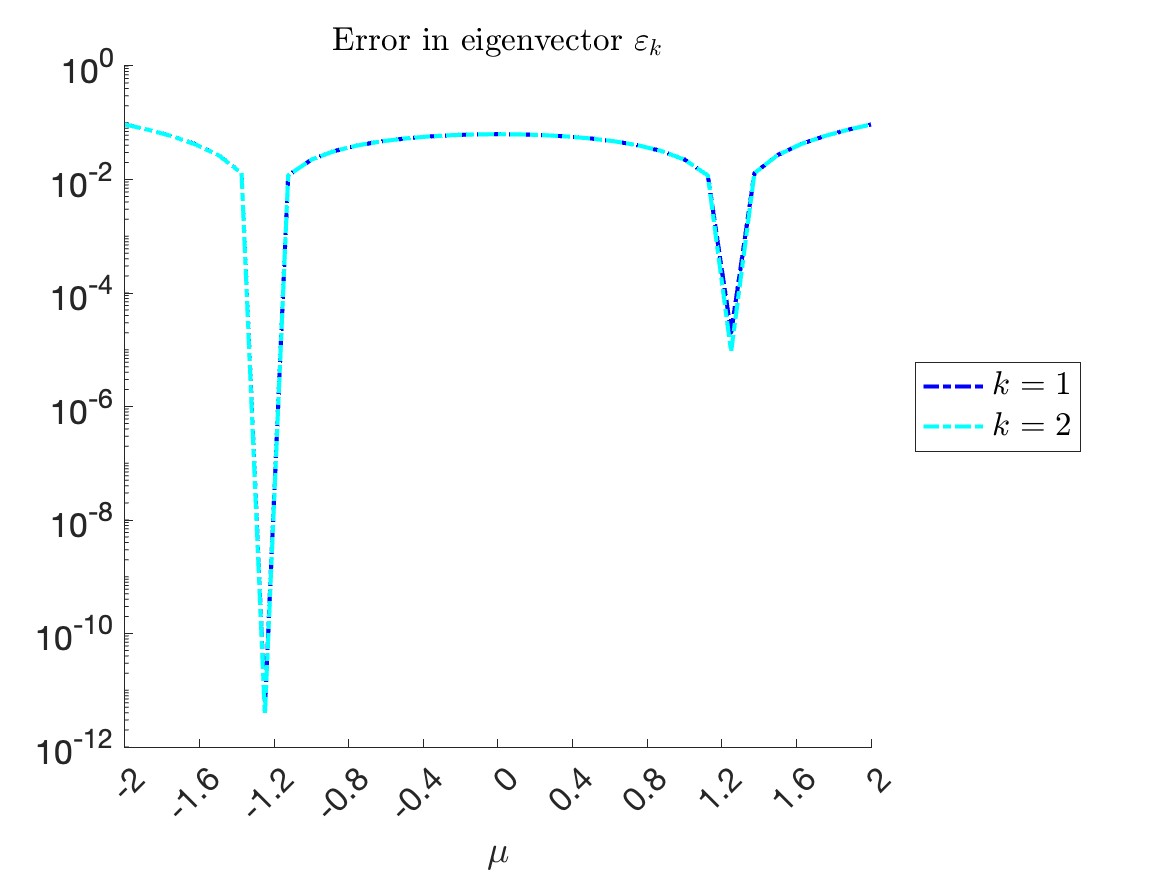}
\caption{Error in first \rev{2} eigenvalues (left) and eigenvectors (right) at different testing parameters 
\rev{$\mu \in \textsf{D}_{\text{test}}$ 
with coarse sampling of training parameters $\mu \in \textsf{D}_{\text{train}} = \{-1.25, 1.25\}$} 
in the example of parametric diatomic well potential.}
\label{fig:schrodinger_cube_error_coarse}
\end{figure}

\begin{figure}[ht!]
\centering
\includegraphics[width=0.48\linewidth]{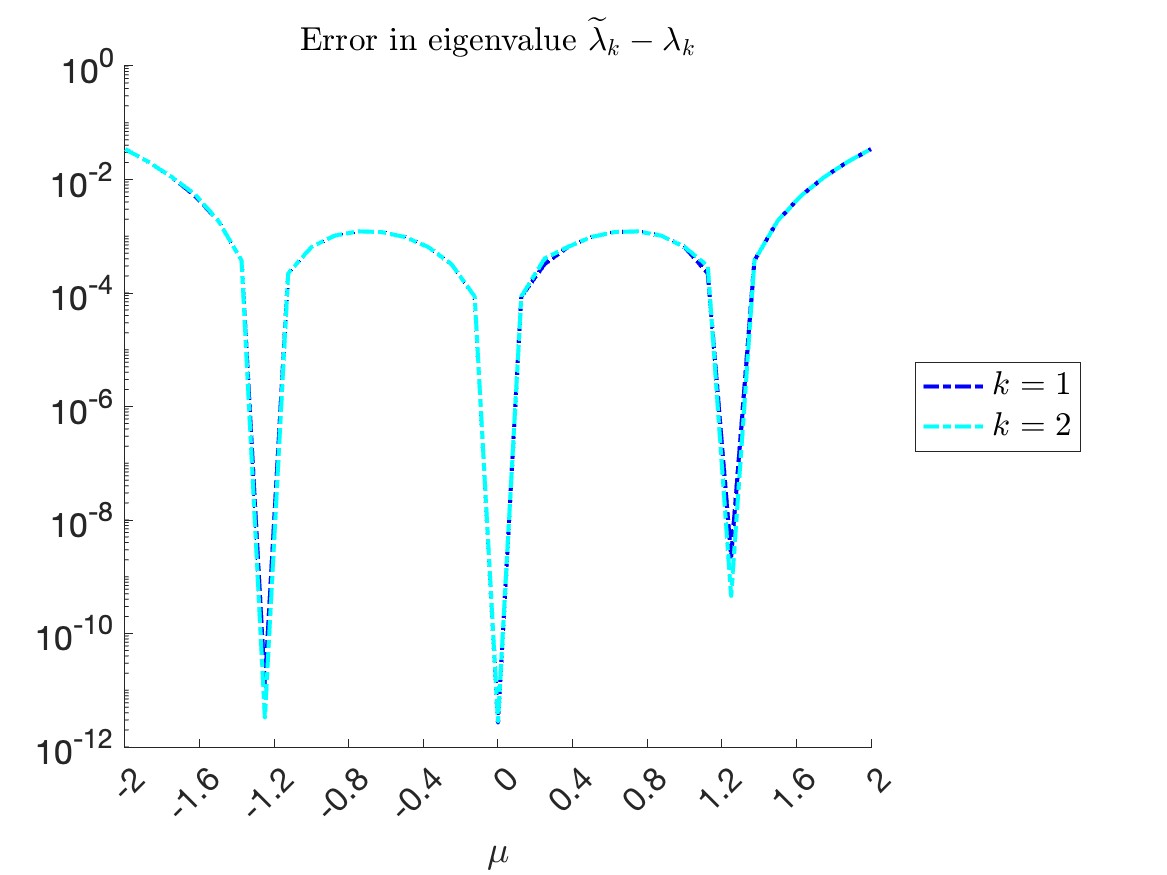}
\includegraphics[width=0.48\linewidth]{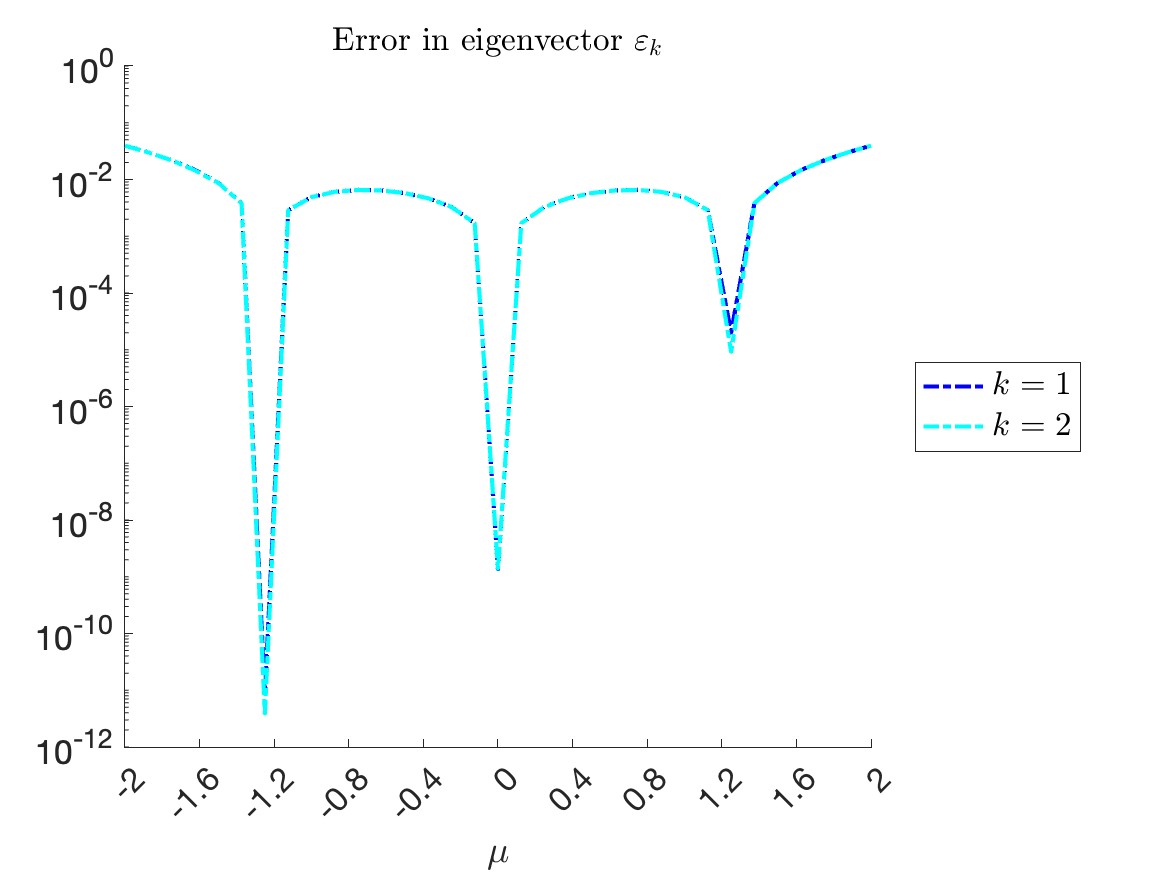}
\caption{Error in first \rev{2} eigenvalues (left) and eigenvectors (right) at different testing parameters 
\rev{$\mu \in \textsf{D}_{\text{test}}$ 
with fine sampling of training parameters $\mu \in \textsf{D}_{\text{train}} = \{-1.25, 0, 1.25\}$} 
in the example of parametric diatomic well potential.}
\label{fig:schrodinger_cube_error_fine}
\end{figure}

We investigate the effects of sampling for collecting the sample snapshots. 
\rev{
We compare results obtained using two different sets of sample parameters: 
$\textsf{D}_{\text{train}} = \{ -1.25, 1.25 \}$ consisting of two sample parameters with a density of 2.5, and 
$\textsf{D}_{\text{train}} = \{ -1.25, 0, 1.25 \}$ consisting of three sample parameters with a density of 1.25.
}
The eigenvectors $\{ \boldsymbol{\phi}_k(\boldsymbol{\mu}) \}_{k=1}^{2}$ at these sample parameters $\boldsymbol{\mu} \in \textsf{D}_{\text{train}}$ are used as snapshots
to construct the basis matrix $\mathbf{Q}$ of column size $r = 4$ for the two-parameter set and $r = 6$ for the three-parameter set, respectively.
The projected system \eqref{eq:rom-evp} is solved with the routine \texttt{dsygv} in LAPACK.
Figure~\ref{fig:schrodinger_cube_error_coarse} and Figure~\ref{fig:schrodinger_cube_error_fine} show 
the error in the first 2 eigenvalues and eigenvectors at different testing parameters 
\rev{$\boldsymbol{\mu} \in \textsf{D}_{\text{test}} = \{ -2 + 0.125 t : t \in [0,32] \cap \mathbb{Z} \}$} for these two sampling configurations.
It can be observed that using a denser sampling does not only make the approximations of the eigenvalues and eigenvectors
become extremely accurate in the additional reproductive cases $\mu = 0$,
but also significantly increases the accuracy all the cases.
At the extrapolation case $\boldsymbol{\mu} = 2$, the ROM approximation error of both eigenvalues is reduced from 0.09 to 0.03, and
the ROM approximation error of both eigenvectors is reduced from 0.17 to 0.04.
We remark that very similar results are observed for the second eigenpair.

\subsection{3D parametric contrast heterogeneous diffusion}

\begin{figure}[ht!]
\centering
\includegraphics[width=0.48\textwidth]{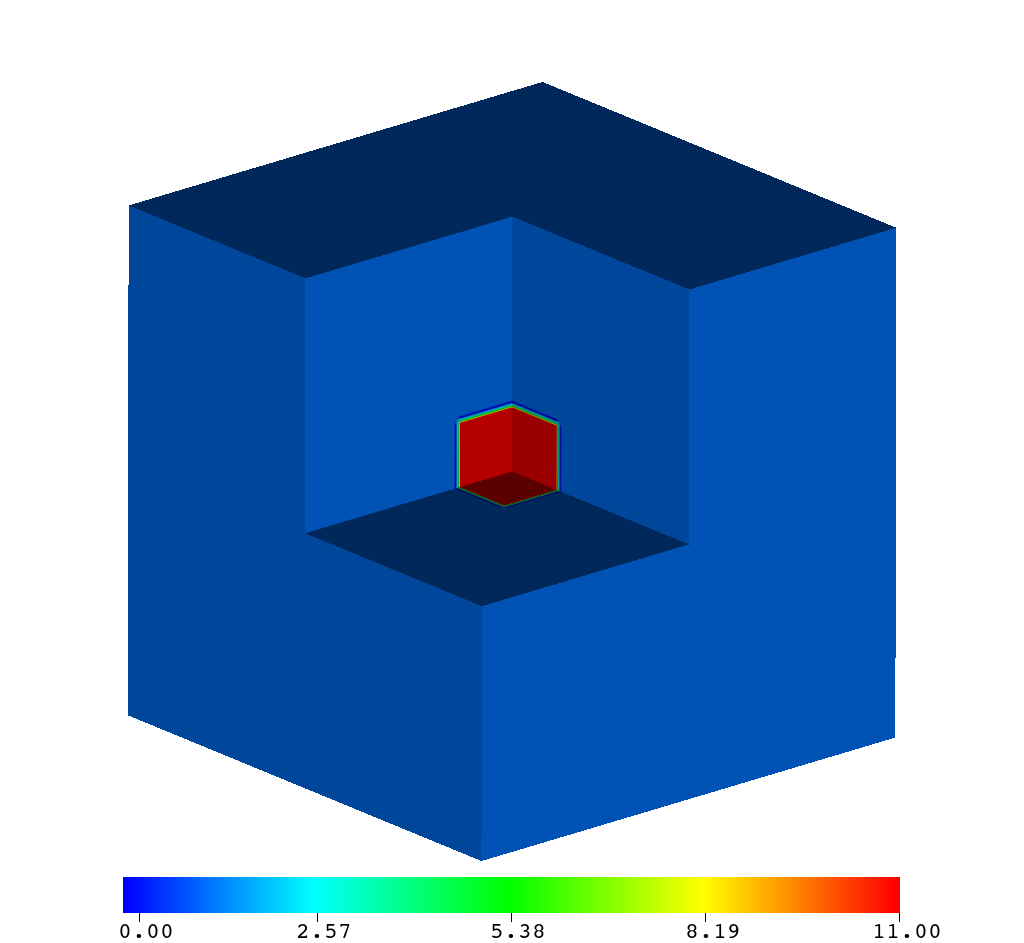}
\caption{Conductivity coefficient $\sigma$ at the parameter $\mu = 10$
 in the example of parametric contrast heterogeneous diffusion.}
\label{fig:diffusion_fichera_conductivity}
\end{figure}

In this example, the computational domain is taken as the Fichera domain  
$\Omega = (-1,1)^3 \setminus (-1,0]^3 \subset \mathbb{R}^3$, and the parameter domain is 
$\textsf{D} = \rev{[0,20]}$ 
\rev{which reduces the parameter $\boldsymbol{\mu}$ to a scalar $\mu$}.
For $\mu \in \textsf{D}$, 
the coefficient fields are given by 
\begin{equation}
\begin{aligned}
\sigma(x; \mu) & = 1 + \mu \chi_B(x) && \text{for } x \in \Omega, \\
\rho(x; \mu) & = 0 && \text{for } x \in \Omega, 
\end{aligned}
\end{equation}
which corresponds to a heterogeneous diffusion operator with parametric contrast between the box 
$B = \{x \in \Omega: \vert x \vert_\infty \leq 0.25\}$ and the background $\Omega \setminus B$, 
and the Dirichlet boundary condition is prescribed, i.e., 
\begin{equation}
\alpha(x; \mu) = 1, \quad \beta(x; \mu) = 0, \quad \text{for } x \in \partial\Omega.
\end{equation}
Figure~\ref{fig:diffusion_fichera_conductivity} depicts the conductivity coefficient $\sigma$ 
in the domain $\Omega$ at the parameter $\mu = 10$.

\begin{figure}[ht!]
\centering
\includegraphics[width=0.48\linewidth]{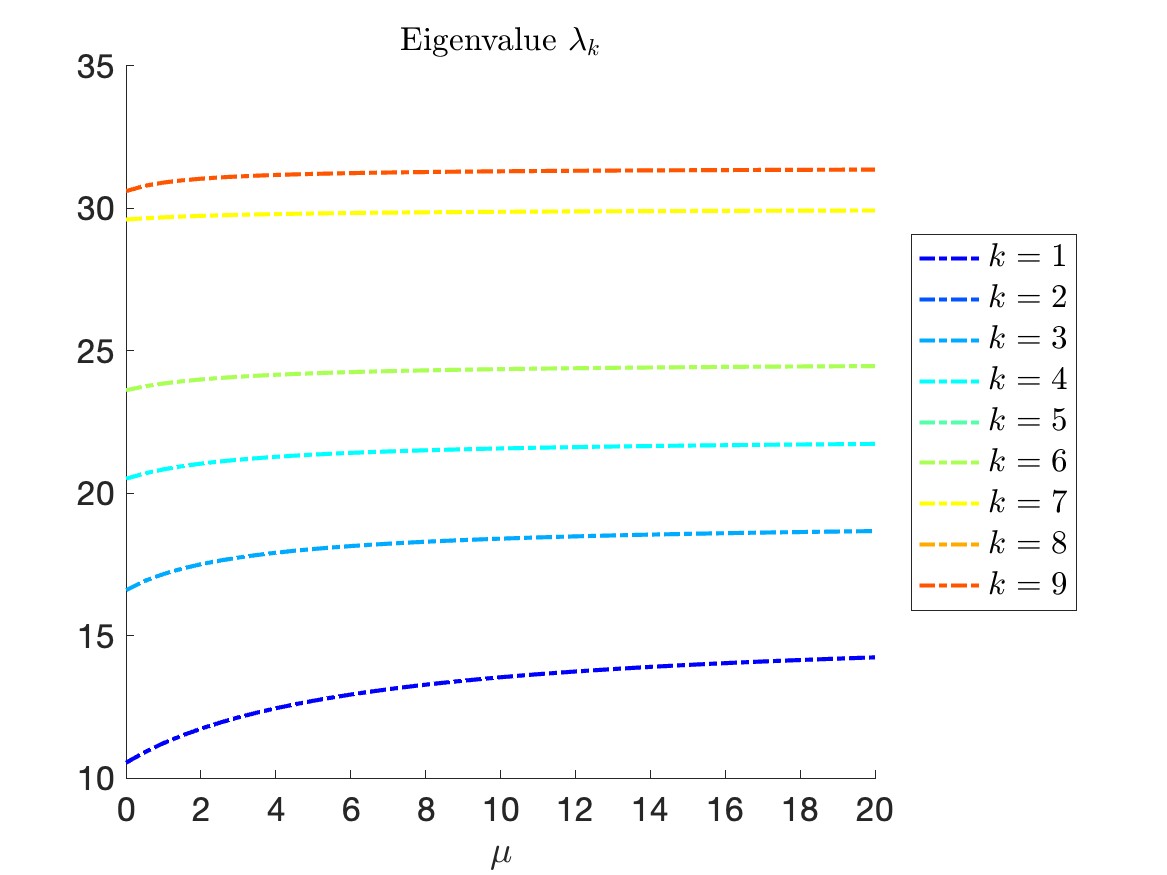}
\caption{First \rev{9} eigenvalues at different parameters $\mu$ in the example of parametric contrast heterogeneous diffusion.}
\label{fig:diffusion_fichera_eigenvalue}
\end{figure}

The domain $\Omega$ is divided into a uniform mesh of size $h = 1/64$, 
and $Q^2$ Lagrange finite elements is used in the finite element discretization, 
which results in a system size of $n = 1884545$. 
We are interested in the first $p = \rev{9}$ eigenpairs. 
The eigenvalue problem \eqref{eq:fom-evp} is solved by LOBPCG.
Figure~\ref{fig:diffusion_fichera_eigenvalue} shows the first \rev{9} eigenvalues $\{\lambda_k(\mu)\}_{k=1}^{\rev{9}}$ 
of \eqref{eq:fom-evp} at different \rev{testing} parameters \rev{$\mu \in \textsf{D}_{\text{test}} = \{ 0.5 t: t \in [0, 40] \cap \mathbb{Z} \}$}. 
In this case, there are 3 sets of double eigenvalue, namely $\lambda_2 = \lambda_3$, 
$\lambda_5 = \lambda_6$, and $\lambda_8 = \lambda_9$, uniformly in $\mu \in [0, \rev{20}]$. 

\begin{figure}[ht!]
\centering
\includegraphics[width=0.48\linewidth]{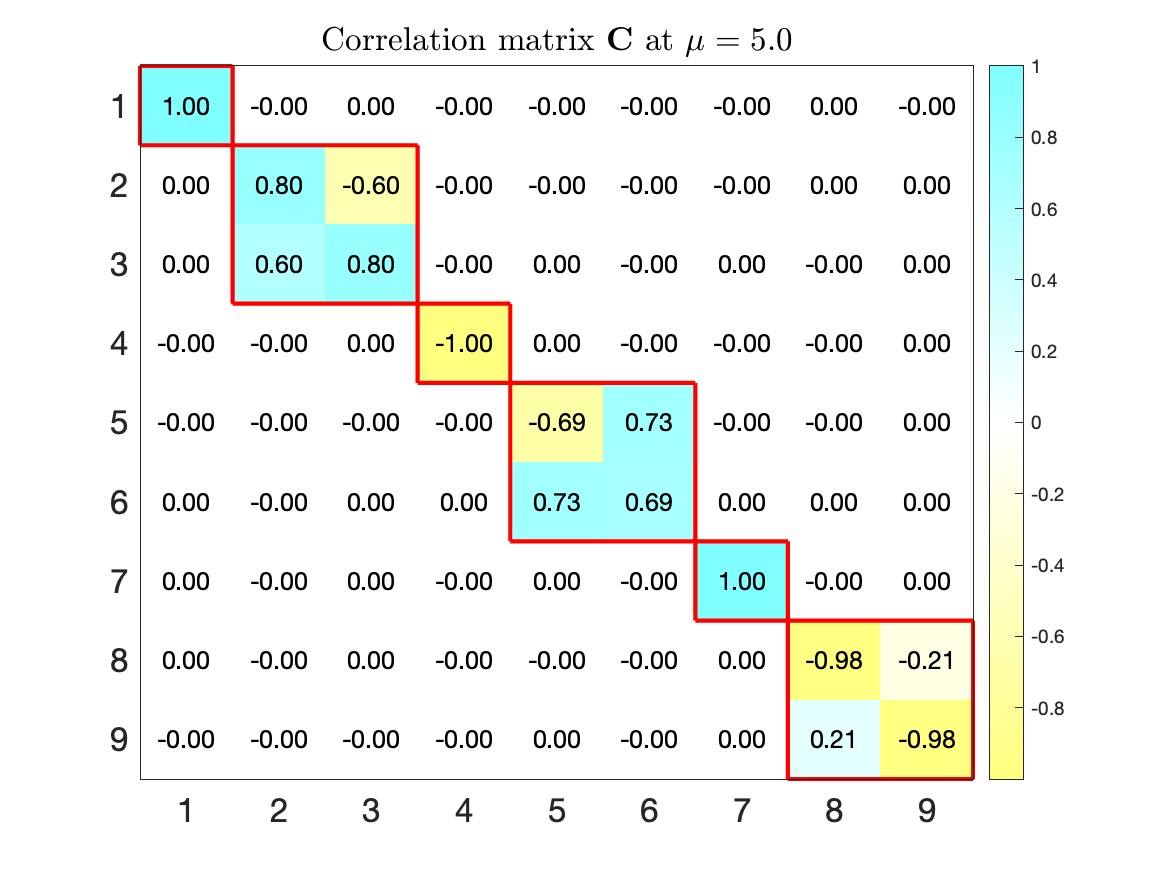}
\includegraphics[width=0.48\linewidth]{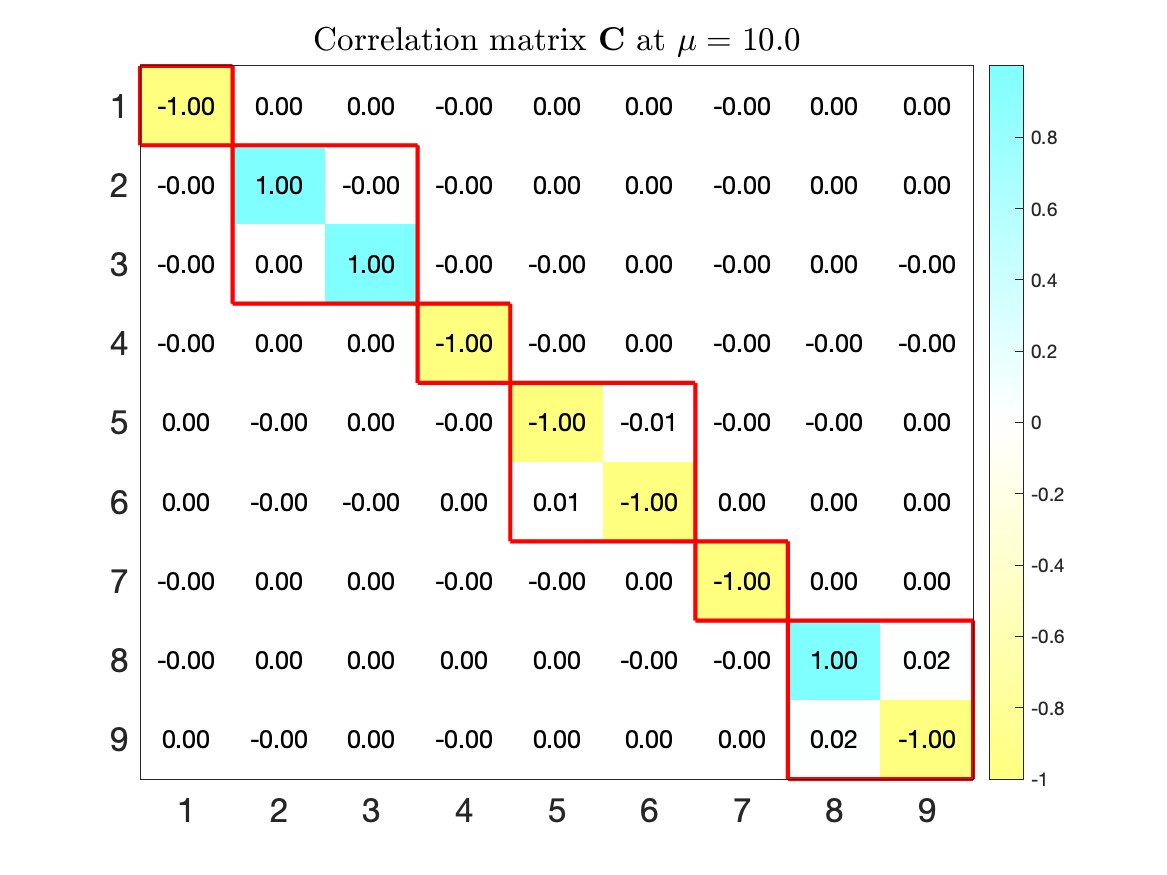}
\caption{Correlation matrix $\mathbf{C}$ at the \rev{parameters $\mu = 5$ (left) and $\mu = 10$ (right)} 
in the example of parametric contrast heterogeneous diffusion.}
\label{fig:diffusion_fichera_correlation_matrix}
\end{figure}

\begin{figure}[ht!]
\centering
\includegraphics[width=0.48\linewidth]{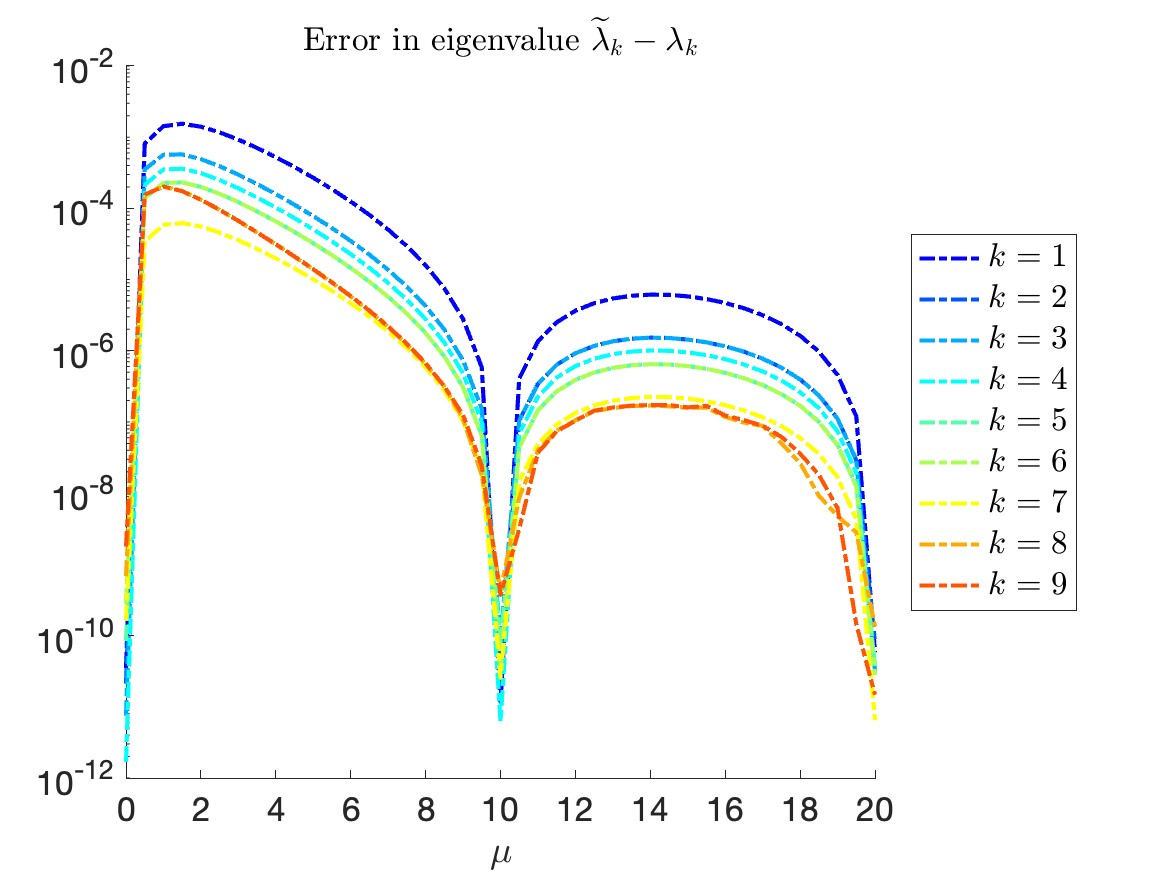}
\includegraphics[width=0.48\linewidth]{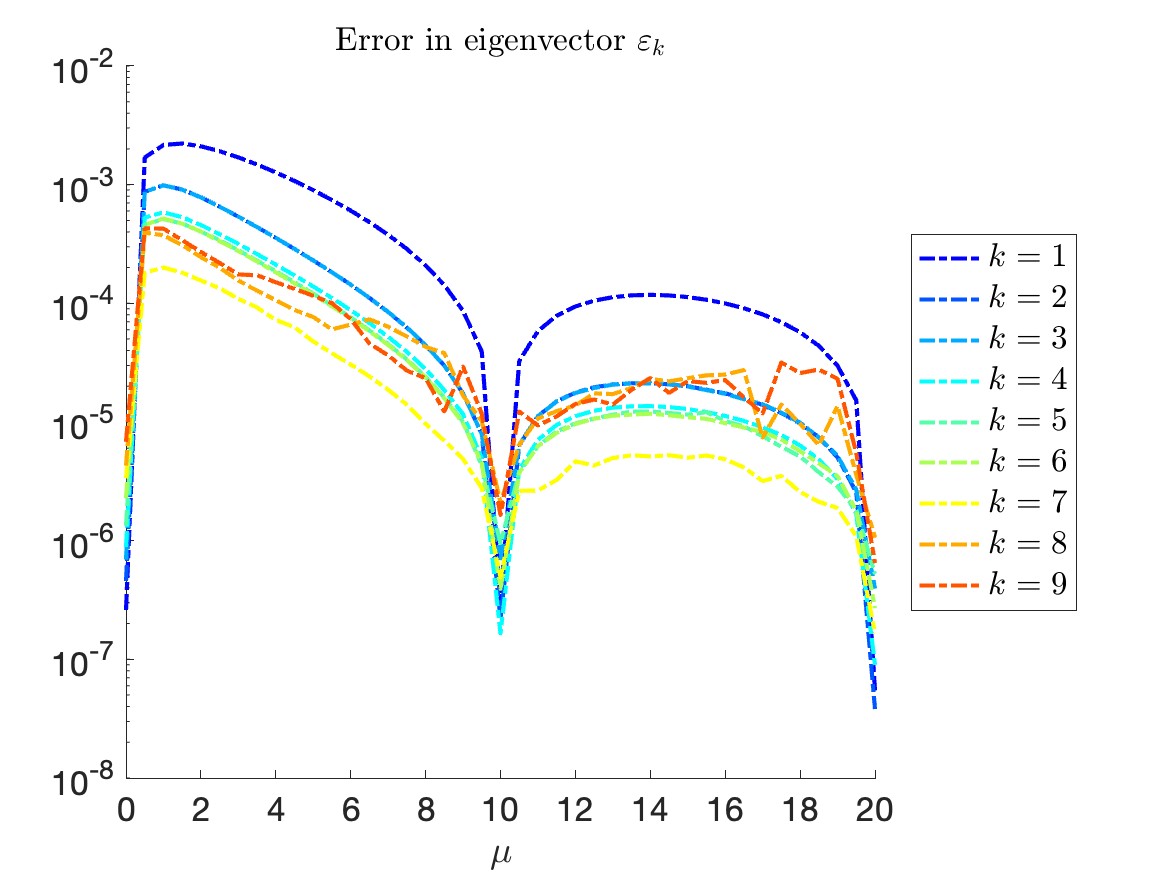}
\caption{Error in first \rev{9} eigenvalues (left) and eigenvectors (right)
 at different testing parameters \rev{$\mu \in \textsf{D}_{\text{test}}$} in the example of parametric contrast heterogeneous diffusion.}
\label{fig:diffusion_fichera_error}
\end{figure}

The eigenvectors $\{ \boldsymbol{\phi}_k(\mu) \}_{k=1}^{\rev{9}}$ at the sample parameters $\mu \in \textsf{D}_{\text{train}} = \{ 0, 10, 20 \}$ are used as snapshots 
to construct the basis matrix $\mathbf{Q}$ of column size $r = \rev{27}$.  
The projected system \eqref{eq:rom-evp} is solved with the routine \texttt{dsygv} in LAPACK.
Due to the repeated eigenvalues, in the calculation of the error in the $k$-th eigenvector, 
we need to identify the subspace $\mathcal{E}_j$ where $\lambda_k = \nu_j$, and 
compute the projection $\widetilde{\mathbf{P}}_{\mathbf{M}}^{\mathcal{S}_j}$, 
which depends on the correlation matrix $\mathbf{C} \in \mathbb{R}^{\rev{9 \times 9}}$, defined by 
\begin{equation}
\mathbf{C}_{km} = \widetilde{\boldsymbol{\phi}}_m^\top \mathbf{M} \boldsymbol{\phi}_k.
\end{equation} 
The ROM approximation of the $k$-th eigenvector can be expressed as 
\begin{equation}
\widetilde{\mathbf{P}}_{\mathbf{M}}^{\mathcal{S}_j} \boldsymbol{\phi}_k 
= \sum_{m: \lambda_m = \lambda_k} \mathbf{C}_{km} \widetilde{\boldsymbol{\phi}}_m.
\end{equation}
Figure~\ref{fig:diffusion_fichera_correlation_matrix} shows all the entries of the correlation matrix $\mathbf{C}$ at the parameter 
\rev{$\mu = 5$ and $\mu = 10$ respectively}. 
The calculation of the ROM approximations of the eigenvectors only utilizes the block matrices surrounded by the red boxes, 
within which the indices of the eigenpairs corresponding to the same eigenspace. 
Moreover, it can be observed that the off-diagonal blocks are approximately zero (up to 2 decimal places), 
meaning the FOM solutions $\boldsymbol{\phi}_k$ are almost orthogonal to all the ROM approximations 
$\widetilde{\boldsymbol{\phi}}_m$ outside their own eigenspace where $\lambda_m \neq \lambda_k$. 
Figure~\ref{fig:diffusion_fichera_error} shows the error in the first \rev{9} eigenvalues and eigenvectors at different parameters $\mu \in [0,20]$. 
Similar to the previous examples, both approximations are extremely accurate in the reproductive case, 
and reasonably accurate in the interpolation cases. 

%% file: conclusion.tex
\section{Conclusion}

In this paper, we presented theoretical and numerical studies of reduced basis approximation for solving parametric linear eigenvalue problems. 
We derived general a-priori error estimates of the eigenvalues and the eigenvectors through the min-max principles. 
The error bounds are controlled by several quantities, including the oblique projectors and the spectral gaps.
We provided various numerical examples of finite element discretization of Laplace, Kohn-Sham, and diffusion eigenvalue problems,
with parametric boundary conditions and coefficient fields. 
We found that, through efficient sampling of snapshot data for constructing the reduced basis, 
the projected eigenvalue problem can provide accurate approximations and be solved thousands times faster. 
The research findings will pave the way of further studies on reduced order modeling of nonlinear eigenvalue problems, 
and more complicated applications such as quantum molecular dynamics.

\section*{Declaration of Competing Interest}
The authors declare that they have no known competing financial interests or personal relationships that could have appeared to influence the work reported in this paper.

\section*{Acknowledgments}

This work was supported by Laboratory Directed Research and Development (LDRD) Program by the U.S. Department of Energy (24-ERD-035). 
Lawrence Livermore National Laboratory is operated by Lawrence Livermore National Security, LLC, for the U.S. Department of Energy, National Nuclear Security Administration under Contract DE-AC52-07NA27344. IM release number: LLNL-JRNL-867049.
This manuscript has been co-authored by UT-Battelle, LLC, under contract DE-AC05–00OR22725 with the US Department of Energy (DOE).

\section*{Disclaimer}
This document was prepared as an account of work sponsored by an agency of the United States government. Neither the United States government nor Lawrence Livermore National Security, LLC, nor any of their employees makes any warranty, expressed or implied, or assumes any legal liability or responsibility for the accuracy, completeness, or usefulness of any information, apparatus, product, or process disclosed, or represents that its use would not infringe privately owned rights. Reference herein to any specific commercial product, process, or service by trade name, trademark, manufacturer, or otherwise does not necessarily constitute or imply its endorsement, recommendation, or favoring by the United States government or Lawrence Livermore National Security, LLC. The views and opinions of authors expressed herein do not necessarily state or reflect those of the United States government or Lawrence Livermore National Security, LLC, and shall not be used for advertising or product endorsement purposes.

%% file: appendix.tex
\rev{
\section{Evaluation of upper bound for $\widetilde{\lambda}_k$}
\label{sec:compute_kappa}
We discuss the procedure of computing $\kappa_k$ in the eigenvalue upper bound in Theorem~\ref{thm:eigval-error}, 
from the matrix of generalized eigenvectors $\boldsymbol{\Phi}$. 
First, we define the projected matrix $\boldsymbol{\Psi}^{(k)} = \mathbf{P}_{\mathbf{A}} \boldsymbol{\Phi}^{(k)}$. This matrix represents the action of the oblique projection on the true eigenvectors. It can be computed explicitly as:
\begin{equation}
    \boldsymbol{\Psi}^{(k)} = \mathbf{Q} (\mathbf{Q}^\top \mathbf{A} \mathbf{Q})^{-1} \mathbf{Q}^\top \mathbf{A} \boldsymbol{\Phi}^{(k)}
\end{equation}
We also define $\mathbf{M}_{\boldsymbol{\Phi}}^{(k)} = (\boldsymbol{\Phi}^{(k)})^\top \mathbf{M} \boldsymbol{\Phi}^{(k)} \in \mathbb{S}_{++}^k$ and $\mathbf{M}_{\boldsymbol{\Psi}}^{(k)} = (\boldsymbol{\Psi}^{(k)})^\top \mathbf{M} \boldsymbol{\Psi}^{(k)} \in \mathbb{S}_{++}^k$. 
We express any vector $\mathbf{y}$ in $\mathcal{R}(\boldsymbol{\Phi}^{(k)}) \setminus \{ \mathbf{0} \}$ as a linear combination of the columns of $\boldsymbol{\Phi}^{(k)}$, i.e., $\mathbf{y} = \boldsymbol{\Phi}^{(k)} \widehat{\mathbf{x}}_k$ for some non-zero vector $\widehat{\mathbf{x}}_k \in \mathbb{R}^k \setminus \{ \mathbf{0} \}$.
Substituting this into the definition of $\kappa_k$, we get:
\begin{equation}
\begin{split}
    \kappa_k & = \sup_{\widehat{\mathbf{x}}_k \in \mathbb{R}^k \setminus \{\mathbf{0}\}} \dfrac{\| \boldsymbol{\Phi}^{(k)} \widehat{\mathbf{x}}_k \|_{\mathbf{M}}}{\| \mathbf{P}_{\mathbf{A}} \boldsymbol{\Phi}^{(k)} \widehat{\mathbf{x}}_k \|_{\mathbf{M}}} \\
    &= \sup_{\widehat{\mathbf{x}}_k \in \mathbb{R}^k \setminus \{\mathbf{0}\}} \dfrac{(\boldsymbol{\Phi}^{(k)} \widehat{\mathbf{x}}_k)^\top \mathbf{M} (\boldsymbol{\Phi}^{(k)} \widehat{\mathbf{x}}_k)}{(\boldsymbol{\Psi}^{(k)} \widehat{\mathbf{x}}_k)^\top \mathbf{M} (\boldsymbol{\Psi}^{(k)} \widehat{\mathbf{x}}_k)} \\
    &= \sup_{\widehat{\mathbf{x}}_k \in \mathbb{R}^k \setminus \{\mathbf{0}\}} \dfrac{\widehat{\mathbf{x}}_k^\top ((\boldsymbol{\Phi}^{(k)})^\top \mathbf{M} \boldsymbol{\Phi}^{(k)}) \widehat{\mathbf{x}}_k}{\widehat{\mathbf{x}}_k^\top ((\boldsymbol{\Psi}^{(k)})^\top \mathbf{M} \boldsymbol{\Psi}^{(k)}) \widehat{\mathbf{x}}_k} \\
    &= \sup_{\widehat{\mathbf{x}}_k \in \mathbb{R}^k \setminus \{\mathbf{0}\}} \dfrac{\widehat{\mathbf{x}}_k^\top \mathbf{M}_{\boldsymbol{\Phi}}^{(k)} \widehat{\mathbf{x}}_k}{\widehat{\mathbf{x}}_k^\top \mathbf{M}_{\boldsymbol{\Psi}}^{(k)} \widehat{\mathbf{x}}_k} \\
    & = \sup_{\widehat{\mathbf{x}}_k \in \mathbb{R}^k \setminus \{\mathbf{0}\}} \mathbf{R}_{\mathbf{M}_{\boldsymbol{\Phi}}^{(k)}, \mathbf{M}_{\boldsymbol{\Psi}}^{(k)}}(\widehat{\mathbf{x}}_k). 
\end{split}
\end{equation}
As a direct consequence of Lemma~\ref{lemma:rayleigh-quotient-fom}, the supremum of this generalized Rayleigh quotient 
is the largest eigenvalue of the generalized eigenvalue problem
\begin{equation}
    \mathbf{M}_{\boldsymbol{\Phi}}^{(k)} \widehat{\mathbf{x}}_k = \widehat{\lambda}_k \mathbf{M}_{\boldsymbol{\Psi}}^{(k)} \widehat{\mathbf{x}}_k.
\end{equation}
}